# On exchangeable random variables and the statistics of large graphs and hypergraphs[*]

Tim Austin

*Department of Mathematics
University of California, Los Angeles
Los Angeles, CA 90095, USA
e-mail:* timaustin@math.ucla.edu
*url:* www.math.ucla.edu/~timaustin

**Abstract:** De Finetti's classical result of [18] identifying the law of an exchangeable family of random variables as a mixture of i.i.d. laws was extended to structure theorems for more complex notions of exchangeability by Aldous [1, 2, 3], Hoover [41, 42], Kallenberg [44] and Kingman [47]. On the other hand, such exchangeable laws were first related to questions from combinatorics in an independent analysis by Fremlin and Talagrand [29], and again more recently in Tao [62], where they appear as a natural proxy for the 'leading order statistics' of colourings of large graphs or hypergraphs. Moreover, this relation appears implicitly in the study of various more bespoke formalisms for handling 'limit objects' of sequences of dense graphs or hypergraphs in a number of recent works, including Lovász and Szegedy [52], Borgs, Chayes, Lovász, Sós, Szegedy and Vesztergombi [17], Elek and Szegedy [24] and Razborov [54, 55]. However, the connection between these works and the earlier probabilistic structural results seems to have gone largely unappreciated.

In this survey we recall the basic results of the theory of exchangeable laws, and then explain the probabilistic versions of various interesting questions from graph and hypergraph theory that their connection motivates (particularly extremal questions on the testability of properties for graphs and hypergraphs).

We also locate the notions of exchangeability of interest to us in the context of other classes of probability measures subject to various symmetries, in particular contrasting the methods employed to analyze exchangeable laws with related structural results in ergodic theory, particular the Furstenberg-Zimmer structure theorem for probability-preserving $\mathbb{Z}$-systems, which underpins Furstenberg's ergodic-theoretic proof of Szemerédi's Theorem.

The forthcoming paper [10] will make a much more elaborate appeal to the link between exchangeable laws and dense (directed) hypergraphs to establish various results in property testing.

Received January 2008.

---

[*]This is an original survey paper

**Contents**

## 1. Introduction

This survey paper is about the laws of random colourings of the complete $k$-uniform hypergraph $\binom{S}{k}$ on a countably infinite vertex set $S$ that are invariant under permutations of the vertex set. When $k = 1$ a complete description of these laws follows from a classical result of de Finetti [18, 19] for $\{0,1\}$-valued exchangeable random variables. More generally, suppose that $(K, \Sigma_K)$ is a standard Borel space (which will serve as our space of 'colours') and that $k \geq 1$. We shall be concerned with the structure of those probability measures $\mu$ on the measurable space $(K^{\binom{S}{k}}, \Sigma_K^{\otimes \binom{S}{k}})$ (the set of all $K$-coloured complete $k$-uniform hypergraphs on $S$) that are invariant under the natural vertex-permuting action $\mathrm{Sym}_0(S) \curvearrowright \binom{S}{k}$, where $\mathrm{Sym}_0(S)$ is the group of finitely-supported permutations of $S$.

Measures (or, equivalently, the associated canonical processes of coordinate-projections onto $K$) enjoying such symmetries were subject to a number of studies during the 1970's and 80's, culminating in the first complete analyses by Kingman [47], Hoover [41, 42], Aldous [1, 2, 3] and Kallenberg [44] for increasingly general classes of process.

More recently, a similar structural description has emerged independently in the work of a group of researchers on 'limit objects' for sequences of large finite graphs or hypergraphs, whose structure can often serve as a 'proxy' for the 'leading order statistics' of such graphs or hypergraphs (see, for example, Lovász and Szegedy [52], Borgs, Chayes, Lovász, Sós, Szegedy and Vesztergombi [17], Elek and Szegedy [24]).

We shall survey the former area (giving a description close in spirit to those of Aldous [3] and Kallenberg [44], where the picture is more complete), and then describe how these two strands of research are actually closely related.

The link between them areas arises because exchangeable random colourings can themselves serve as such limit objects, and because once this identification is made many of the results of the more recent formalisms simply follow from the older structure theorems for exchangeable laws. This basic identification seems to appear first in Tao [62], where parts of the older structure theory are then implicitly re-proved, but without a development of the full formalism.

In addition, similar structural results were already obtained by Fremlin and Talagrand in [29] (a paper that seems to have gone unnoticed by many more recent researchers) for a class of random graphs that are subject only to a rather weaker symmetry than full exchangeability: in the terminology we will adopt below their random graphs on $\mathbb{N}$ are 'spreadable', according to which all induced finite random subgraphs on a fixed number of vertices have the same law so long as the order of those vertices is respected. This requirement of order-preservation demands a more subtle analysis than in the exchangeable case. Spreadability was studied by Ryll-Nardzewski [58] and Kallenberg [44] as a natural weakening of the hypothesis of exchangeability, but the work of Fremlin and Talagrand [29] also includes an analysis of a related extremal problem (about critical densities for finite subgraphs of certain infinite random graphs), in an early precursor to more recent work relating such questions in finitary combinatorics to the analysis of these random graphs.

Tao's use of exchangeable random hypergraph colourings as such proxies is motivated by a Furstenberg-like correspondence principle between properties of finite hypergraphs and those of exchangeable random hypergraphs. He goes on to give an infinitary analysis of versions of the graph and hypergraph removal lemmas. This makes concrete certain analogies between Furstenberg's ergodic-theoretic work and hypergraph-based approaches to proving Szemerédi's Theorem, and in many respects the present paper is a continuation of this program. In [10] we will extend the infinitary methods of [62] to give an infinitary account of general hypergraph property testing, calling on the structural result described in the present paper.

Another abstract approach to the asymptotic statistics of large dense graphs and hypergraphs has recently appeared in work of Razborov [54], [55]. His construction rests on the notion of a 'flag algebra', constructed from collections of combinatorial structures in a more abstract algebraic manner. Here, too, there is a close parallel between the analysis of the infinitary structures that result and the earlier works pertaining to exchangeable random hypergraph colourings.

In this survey we will first recall versions of the basic results of the theory of exchangeable laws, and will then examine how various purely combinatorial questions admit a parallel version in the setting of these laws, and can occasionally shed light on the original finitary versions through a 'limit object' analysis. In the process of describing this link, we will show how the various other infinitary formalisms described above all recover essentially the same structure as the study of exchangeable processes. We will finish by locating these basic underlying structural results in the broader context of ergodic theory, where a related but necessarily less complete analysis underlies the fundamental Furstenberg-Zimmer structure theory for probability-preserving systems, and so — through a correspondence principle with finitary combinatorics similar to that mentioned above — enables Furstenberg's proof of Szemerédi's Theorem.

A much more thorough account of the theory of exchangeability and related symmetries for stochastic processes, as well as its historical development, can be found in the book [46] of Kallenberg. The treatment of this theory in the present survey will also be skewed to better exhibit its relationship with the more recent work in combinatorics, since the versions of the probabilistic results most central to this relationship (our Corollary 3.5 and Theorems 2.9 and 3.21) are not quite the most general known to probabilists (which are more closely related to the setting of partite hypergraphs that we examine in Subsection 3.7).

**Remark.** As this paper neared completion, many of the main relations between exchangeability and hypergraph theory that it was written to advertise were also independently reported in work of Diaconis and Janson [20].

## 2. Preliminaries

### 2.1. Background notation and definitions

*Combinatorics*

In this paper we will often be concerned with uniform hypergraphs over some countably infinite vertex set $S$. We shall write $\binom{S}{k}$ for the set of $k$-subsets of $S$ and $\binom{S}{\leq k}$ (resp. $\binom{S}{<k}$) for $\bigcup_{j\leq k}\binom{S}{j}$ (resp. $\bigcup_{j<k}\binom{S}{j}$), including $\binom{S}{<\infty}$ for the set of all finite subsets of $S$.

We shall denote the subset of tuples $(s_1, s_2, \ldots, s_i) \in S^i$ with all $s_j$ distinct by $\mathrm{Inj}([i], S)$ (a notation that will re-appear in [10], where injections are conveniently treated as the morphisms of a category of index sets).

In addition to ordinary $k$-uniform hypergraphs, we will need to consider multicoloured hypergraphs.

**Definition 2.1** (Finite palettes). By a **finite $k$-palette** we understand a sequence $K = (K_i)_{i=0}^{k}$ of finite sets. We shall refer to $k$ as the **rank** of the palette.

**Remark.** We shall later have cause to extend the above definition to general $k$-palettes, whose individual spaces are endowed with a standard Borel structure but need not be finite. ◁

**Definition 2.2** (Coloured hypergraphs). Given a vertex set $S$ and a finite set $K_k$, we define a $K$**-coloured $k$-uniform hypergraph on** $S$ to be a map $H : \binom{S}{k} \to K_k$, and a $K$**-coloured $k$-uniform hypergraph on $S$ with loops** to be a map $H : V^k \to K_k$ that is invariant under coordinate permutations of $V^k$.

More generally, given a finite $k$-palette $K = (K_i)_{i=0}^k$ we define a $K$**-coloured hypergraph on** $S$ to be a sequence $H_i$ of $K_i$-coloured $i$-uniform hypergraphs on $S$ for $i = 0, 2, \ldots, k$.

We will also need the analog of the above definition in the directed case.

**Definition 2.3** (Coloured directed hypergraphs). Given a vertex set $S$ and a finite set $K_k$, we define a $K$**-coloured directed $k$-uniform hypergraph on** $S$ to be a map $H : \mathrm{Inj}([k], S) \to K_k$, and a $K$**-coloured $k$-uniform hypergraph on $S$ with loops** to be a map $H : V^k \to K_k$, and more generally given a finite $k$-palette $K = (K_i)_{i=0}^k$ we define a $K$**-coloured directed hypergraph on** $S$ to be a sequence $H_i$ of $K_i$-coloured directed $i$-uniform hypergraphs on $S$.

Given a $K$-coloured (directed) hypergraph $H$ on $V$ and a subset of vertices $W \subseteq V$ we write $H|_W$ for the induced (directed) hypergraph on $W$ given by the sequence of maps $H_i|_{W^i}$.

If $S$ is endowed with a total order $<$, we shall sometimes write $\{s_1 < s_2 < \cdots < s_n\}$ for a subset $\{s_1, s_2, \ldots, s_n\}$ whose members are understood to have been listed in increasing order.

We shall write $\mathrm{Sym}_0(S)$ or just $\mathrm{Sym}_0$ for the group of finitely-supported permutations of $S$, and given a subset $I \subseteq S$ and some $g \in \mathrm{Sym}_0(S)$ we shall write $g(I)$ for $\{g(x) : x \in I\}$.

*Measure theory and probability*

We assume familiarity with basic measure theory and probability; a suitable reference for all the background we will need is Kallenberg [45]. We recall here only some slightly more specialized notions.

We shall generally work with measurable spaces that are **standard Borel spaces**: that is, that are isomorphic as measurable spaces to a Polish space with its Borel $\sigma$-algebra. Occasionally we shall use a proper $\sigma$-subalgebra of the Borel $\sigma$-algebra of such a space. The basic properties of these spaces are described in Appendix A of [45] (under the shorter name 'Borel spaces'). When we refer to a space $K$ as standard Borel, its $\sigma$-algebra is to be understood, and will be denoted by $\Sigma_K$ when necessary; similarly, we shall usually refer to a probability measure 'on $K$' in place of 'on $\Sigma_K$'. The standard Borel assumption has various important technical consequences for the management of $\sigma$-subalgebras and measurable functions: for example, given a probability measure $\mu$ on such a space $X$, we may adopt the common convention from ergodic theory of identifying a $\sigma$-subalgebra T of $\Sigma_X$ with a factor map $\phi : X \to Y$ into some other standard Borel space $Y$, defined uniquely up to equality $\mu$-almost everywhere. When endowed with a Borel probability measure, we will refer to the resulting measure space as a **Lebesgue space**; note that point masses are allowed in this convention, so these spaces may have atoms.

The reader who is uncomfortable with these definitions will lose nothing by assuming throughout that our underlying spaces are compact metrizable and carrying their Borel measurable structure.

Suppose that $Y$ and $X$ are standard Borel spaces. Then by a **probability kernel from $Y$ to $X$** we understand a function $P : Y \times \Sigma_X \to [0, 1]$ such that

- the map $y \mapsto P(y, A)$ is $\Sigma_Y$-measurable for every $A \in \Sigma_X$;
- the map $A \mapsto P(y, A)$ is a probability measure on $\Sigma_X$ for every $y \in Y$.

One instructive intuition about such a kernel is that it serves as a 'randomized map' from $Y$ to $X$: rather than specify a unique image in $X$ for each point $y \in Y$, it specifies only a probability distribution $P(y, \,\cdot\,)$ from which we should choose a point of $X$. The first of the above conditions is then the natural sense in which this assignment of a probability distribution is measurable in $y$; and, indeed, a popular alternative definition of probability kernel is as a measurable function from $Y$ to the set $\Pr X$ of probability measures on $X$. We will adopt the standard notation $P(y, \mathrm{d}x)$ for the measure associated by such a kernel to the point $y$, and will write $P : Y \rightsquigarrow X$ when $P$ is as above.

Given a kernel $P : Y \rightsquigarrow X$ and a probability measure $\nu$ on $Y$, we define the measure $P_\# \nu$ on $X$ by

$$P_\# \nu(A) := \int_Y P(y, A) \, \nu(\mathrm{d}y);$$

this measure on $X$ can be interpreted as the law of a member of $X$ selected randomly by first selecting a member of $Y$ with law $\nu$ and then selecting a member of $X$ with law $P(y, \,\cdot\,)$. By analogy with the case of a function between measurable spaces, we will refer to this as the **pushforward** of $\nu$ by $P$. This extends standard deterministic notation: given a measurable function $\phi : Y \to X$, we may associate to it the deterministic probability kernel given by $P(y, \,\cdot\,) = \delta_{\phi(y)}$ (the point mass at the image of $y$ under $\phi$), and now $P_\# \nu$ is the usual pushforward measure $\phi_\# \nu$.

Certain special probability kernels naturally serve as adjoints to factor maps, in the sense of the following theorem.

**Theorem 2.4.** *Suppose that $Y$ and $X$ are as above, that $\mu$ is a probability measure on $X$ and that $\phi : X \to Y$ is a measurable factor map. Then, denoting the push-forward $\phi_\# \mu$ by $\nu$, there is a $\nu$-almost surely unique probability kernel $P : Y \rightsquigarrow X$ such that $\mu = P_\# \nu$ and which* **represents the conditional expectation with respect to $\phi$**: *for any $f \in L^1(\mu)$, the function*

$$x_1 \mapsto \int_X f(x) \, P(\phi(x_1), \mathrm{d}x)$$

*is a version of the $\mu$-conditional expectation of $f$ with respect to $\phi^{-1}(\Sigma_Y)$. We call this $P$ the **disintegration of $\mu$ over $\phi$**.*

*Proof.* See Theorem 6.3 in Kallenberg [45]. □

Motivated by this result, we will sometimes refer to data $(Y, \nu, P)$ for which $\mu = P_\#\nu$ as a **quasifactor** of $(X, \mu)$, even when it does not arise as the adjoint of some factor map. An alternative convention is to interpret a quasifactor as a probability distribution on the set $\Pr X$ of probability measures on $X$, where $\Pr X$ is given its 'evaluation $\sigma$-algebra'. This corresponds to our definition by identifying $\nu \in \Pr Y$ and $P : Y \rightsquigarrow X$ with the *law* of $P$ under $\nu$ regarded as a measurable function $Y \to \Pr X$. A more detailed discussion of quasifactors in ergodic theory can be found in Chapter 8 of Glasner [35], where this alternative convention is used.

It is clear that in general a probability kernel need not correspond to a disintegration over some factor map: indeed, if $P$ is the disintegration of $\mu$ over a map $\phi$ then it must obey the extra condition that its fibre measures $P(y, \cdot)$ are mutually singular and concentrated on the fibres $\phi^{-1}\{y\}$. Thus, for example, letting $p : [0,1]^2 \to (0,\infty)$ be some strictly positive Borel measurable function for which

$$\int_0^1 p(x,y)\,\mathrm{d}y = 1$$

for all $x \in [0,1]$ and defining $P : [0,1] \rightsquigarrow [0,1]$ by

$$P(x, A) = \int_A p(x,y)\,\mathrm{d}y$$

gives a probability kernel with all fibre measures mutually absolutely continuous, and so this $P$ is not a disintegration of any measure over some measurable map $[0,1] \to [0,1]$. If the probability kernel $P$ does arise as a disintegration over some map $\phi : X \to Y$, we will write that $\phi$ **recovers** $P$; if no such $\phi$ exists we call $P$ **unrecoverable**. It turns out that under reasonably general hypotheses the mutual singularity of the measures $P(y, \cdot)$ is sufficient, as well as necessary, for recoverability, but we will not make use of this fact; see, for example, Theorem 8.3 in Glasner [35]. We will be interested in certain instances of unrecoverability later in this paper. It is worth noting, however, that any quasifactor $P : Y \rightsquigarrow X$ with $\mu = P_\#\nu$ defines a joining $\lambda$ on $(Y \times X, \Sigma_Y \otimes \Sigma_X)$ of $(Y, \nu)$ and $(X, \mu)$ by setting

$$\lambda := \int_Y \delta_y \otimes P(y, \cdot)\,\nu(\mathrm{d}y),$$

and now we can identify both $(X, \mu)$ and $(Y, \nu)$ as true factors of the larger probability space $(Y \times X, \Sigma_Y \otimes \Sigma_X, \lambda)$ so that $P$ gives the conditional distribution of the second coordinate relative to the first: that is, if $Q : Y \rightsquigarrow Y \times X$ is the disintegration of $\lambda$ over the first coordinate and $\pi_X : Y \times X \to X$ is the projections onto the second coordinate then $P = \pi_X \circ Q$.

Various natural operations on either measures or functions have analogs for probability kernels, and we will need some of these. If $P : Y \rightsquigarrow X$ and $S$ is some index set then we will write $P^{\otimes S}$ for the probability kernel $Y \rightsquigarrow X^S$ defined by $P^{\otimes S}(y, \cdot) := P(y, \cdot)^{\otimes S}$, the simple $S$-indexed product of copies of the probability measure $P(y, \cdot)$ (it is routine to check that this retains the measurability properties required of a probability kernel). On the other hand, if $Z$ is a third standard Borel space and $\phi : Z \to Y$ and $Q : Z \rightsquigarrow Y$, we define the **compositions** $P \circ \phi : Z \rightsquigarrow X$ and $P \circ Q : Z \rightsquigarrow X$ by

$$P \circ \phi(z, \cdot) := P(\phi(z), \cdot), \quad \text{and} \quad P \circ Q(z, \cdot) := \int_Y P(y, \cdot) Q(z, \mathrm{d}y).$$

In the remainder of this paper we shall be interested in various classes of **stochastic process**: a family of random variables $\pi_t$ all defined on a single 'background' probability space, indexed by some set $T$ and taking values in standard Borel target spaces $K^{(t)}$. In many of our examples these target spaces $K^{(t)}$ will all agree; more generally there will be a short list of different target spaces $K_0, K_2, \ldots, K_k$, with $\pi_t$ taking values in $K_i$ according as $t$ falls in the $i^{\text{th}}$ cell of some partition $T = T_0 \cup T_1 \cup \cdots \cup T_k$.

It is a common and convenient practice in probability to regard the background probability space that supports these random variables as 'hidden', and to perform all desired constructions in terms of the particular functions $\pi_t$; however, in this survey we shall instead work largely with **canonical processes**: that is, those for which the underlying probability space is itself the product $X = \prod_{t \in T} K^{(t)} = K_0^{T_0} \times K_1^{T_1} \times \cdots K_k^{T_k}$ with its product measurable structure and some probability measure, and $\pi_t : X \to K^{(t)}$ is just the projection onto the coordinate indexed by $t$ for each $t \in T$. Given this projection interpretation, we will sometimes extend our notation by writing $\pi_I$ for the projection of $\prod_{t \in T} K^{(t)}$ onto $\prod_{t \in I} K^{(t)}$ for any subset $I \subseteq T$. The underlying probability measure on $X$ is referred to as the **joint law** of the process.

Not only is this picture quite concrete and intuitive for the arguments that we will need, but it will also emerge naturally when we turn to the 'correspondence principles' relating our work to structures in combinatorics.

### 2.2. Exchangeable families of random variables

Our main objects of study will be families $(\pi_t)_{t \in T}$ of random variables with values in some spaces of 'colours' $K_0, K_1, \ldots$ or $K_k$, and subject to the additional hypothesis that their joint law is invariant under some class of permutations of the (countably infinite) index set $T$. Clearly these permutations must preserve the partition $T_0 \cup T_1 \cup \cdots \cup T_k$ in order for this to make sense. We will discuss later how such processes can serve as 'proxies' for the statistics of some large combinatorial structure (such as a dense graph or hypergraph), and for these examples such a symmetry assumption will arise owing to an equivalent symmetry in the statistics we wish to count: for example, the number of induced four-simplices in a three-uniform hypergraph is unaffected by vertex-permutations of that hypergraph.

We will consider different theories corresponding to different classes of permutations. If $\Gamma$ is a group of permutations of $T$, preserving the partition $T_0 \cup T_1 \cup \cdots \cup T_k$ where relevant, then following Aldous [3] we refer to the joint law $\mu$ of the family $(\pi_t)_{t \in T}$ as $(T, \Gamma)$-**exchangeable** if it is invariant under the coordinate-permuting action of $\Gamma$: that is, writing $\tau^g : K_0^{T_0} \times K_1^{T_1} \times \cdots K_k^{T_k} \to K_0^{T_0} \times K_1^{T_1} \times \cdots K_k^{T_k}$ for the coordinate-permuting action

$$\tau^g\big((\omega_t)_{t \in T}\big) = (\omega_{g(t)})_{t \in T} \qquad g \in \Gamma, \ (\omega_t)_{t \in T} \in K_0^{T_0} \times K_1^{T_1} \times \cdots K_k^{T_k},$$

we ask that $\mu = (\tau^g)_\# \mu$ for all $g \in \Gamma$. We shall extend this notation for the action $\tau$ to apply to maps between different index sets: if $\psi : T \to T'$ is a bijection that respects corresponding partitions $T_0 \cup T_1 \cup \cdots \cup T_k$ and $T'_0 \cup T'_1 \cup \cdots \cup T'_k$ then define $\tau^\psi : K_0^{T'_0} \times K_1^{T'_1} \times \cdots K_k^{T'_k} \to K_0^{T_0} \times K_1^{T_1} \times \cdots K_k^{T_k}$ (notice the contravariance) by

$$\tau^\psi\big((\omega_{t'})_{t' \in T'}\big) = (\omega_{\psi(t)})_{t \in T}.$$

We choose to make both $T$ and $\Gamma$ explicit in our nomenclature since many of our examples will be obtained from actions of the same group, but on different sets. In particular, the leading examples of the theory of exchangeability are the following:

1. Let $S$ be a countably infinite set, whose members we shall refer to as 'vertices', and let $T = S$ and $\Gamma = \mathrm{Sym}_0(S)$, the group of all finitely-supported permutations of $S$.[1] Now we may interpret $K$ as a set (possibly a continuum) of 'colours', and a point of $K^S$ as a $K$-colouring of the vertices in $S$. A joint law on $K^S$ may be interpreted as a random such colouring (with the individual random variables $\pi_s$ giving the colours of individual vertices), and the exchangeability of a random $K$-colouring asserts that its law does not depend on the ordering of the vertices. This example (albeit not with this terminology) is the oldest and simplest instance of exchangeability to have been studied, leading to the complete structural description given by work of de Finetti [18, 19], Dynkin [21] and Hewitt and Savage [40], and was the genesis of most subsequent work; we shall review it in Subsection 3.1 below.
2. Next let $S$ and $\Gamma$ be as above, but take for $T$ the set $\binom{S}{k}$ of all $k$-*hyperedges* of the complete graph on $S$. Now our joint law corresponds to a random $K$-colouring of the complete $k$-uniform hypergraph on $S$ that is invariant under all finitely-supported vertex-set permutations. In Subsection 3.2 we shall recount the complete description of these in case $k = 2$ that follows from the work of Kingman [47], Hoover [42] and Aldous [2], and shall then extend those results to general hypergraph colourings in Subsection 3.3 following Kallenberg [44]. (Although we note that their work actually corresponds to the more general setting of *partite* graphs and hypergraphs, to which we turn in item 5 below.)

---

[1] We make the assumption of finite support only because working with all permutations introduces the additional technicalities of working with an uncountable group; however, with the right conventions these are routinely surmountable, and the resulting theory is easily seen to be equivalent.

3. However, in the course of analyzing the above examples we shall meet a natural need to generalize them further, in particular in order to formulate an inductive argument that can be closed. Thus we shall be mostly concerned with random $K$-colourings on $S$ for some general $k$-palette $K = (K_i)_{i=0}^k$, with the above example corresponding to the case in which all $K_i$ for $i \leq k - 1$ are singletons. It is this setting that motivates the above introduction of a possibly nontrivial index set partition, since we shall want to write $T = \bigcup_{i \leq k} T_i$ with $T_i := \binom{S}{i}$.
4. A different extension of example 2 can be obtained by allowed coloured directed hypergraphs: these correspond simply to processes indexed by the injections in $\mathrm{Inj}([k], S)$ rather than the subsets in $\binom{S}{k}$. As in the undirected case, we shall also then need to allow such injections for several different values of $k$, and so obtain a further theory that amounts to a coupling of several $\mathrm{Inj}([j], S)$-indexed processes for different $j \leq k$.
5. Viewing our above examples as random colourings of some combinatorial structure that enjoy a certain related symmetry, it is easy to cook up more in the same way. The range of new ideas that are needed for their analysis seems to be limited, but two in particular are worthy of discussion.

    Firstly, we will consider random colourings of partite hypergraphs. In fact, it is this combinatorial setup that corresponds most closely to the original study of 'exchangeable random arrays' in the probabilistic literature: consider, in particular, the presentation in Section 14 of Aldous [3]. In the simplest case, that of bipartite graphs, we require two countable infinite vertex classes $S_1, S_2$, and then have $T = S_1 \times S_2$ and $\Gamma$ the group of permutations of $T$ obtained by a separate permutation of each $S_i$. A more complicated partite example is of interest as it is the natural infinitary proxy for the simplex-counting that underlies modern hypergraph proofs of Szemerédi's Theorem (see, in particular, Nagle, Rödl and Schacht [53] and Gowers [36]).

    Secondly, within a space of 'colourings' such as $K^T$ we may obtain certain $\Gamma$-invariant subspaces by demanding various relationships among the colours of different points of $T$, provided these demands are themselves permutation-invariant. As a refinement of our earlier directed hypergraph examples, these can be naturally formulated in great generality as random models of certain kinds of theory. In fact, these constituted some of the first instances of exchangeability for higher ranks to appear in the literature, for example in the model theoretic work of Gaifman [34] and Krauss [51], and in [42] Hoover cites these works as motivation for their further study. More recently, they have appeared appeared implicitly in the work of Razborov [54, 55], and in Subsection 4.3 we shall briefly relate his formalism to the more classical probabilistic analysis.

Note that from the above list we have omitted the most heavily-studied instance of exchangeability of all: stationary processes, for which $T = \Gamma = \mathbb{Z}$ (or, more generally, some other group) acting on itself by (right-) translation. As is standard, these stationary processes include (isomorphic copies of) all probability-preserving Γ-systems on some abstract standard Borel probability space: the basic objects of ergodic theory.

However, that setting is, of course, quite wild: such general systems can exhibit hugely varying kinds of behaviour, most of which seems beyond a sensible classification. What is more surprising is that the more special examples in the above list all *do* admit some quite precise structure theorem for the exchangeable laws, proved using rather different, more probabilistic arguments than those common in general ergodic theory (which lean more towards representation theory and harmonic analysis). The closest parallel results in general ergodic theory are perhaps those of Furstenberg-Zimmer structure theory ([30, 31, 69, 68]), and we shall make a rough comparison between our exchangeable families and that theory in Section 4.7.

Let us mention also another symmetry condition for processes closely related to exchangeability. Given a law $\mu$ on the product space $K^{\binom{\mathbb{N}}{k}}$, we refer to it as **spreadable** if it is invariant under any passage to subsequences: that is, if the law $\mu$ of $(\pi_e)_{e \in \binom{\mathbb{N}}{k}}$ is the same as that of $(\pi_{\phi(e)})_{e \in \binom{\mathbb{N}}{k}}$ for any strictly increasing function $\phi : \mathbb{N} \to \mathbb{N}$. Spreadable processes formed an end-point of the original vein of probabilistic research into exchangeability and related symmetries. In [44], Kallenberg showed that for any spreadable set-indexed process $(\pi_\alpha)_{\alpha \in \binom{\mathbb{N}}{k}}$ there is some exchangeable process $(\tilde{\pi}_\phi)_{\phi \in \mathrm{Inj}([k],\mathbb{N})}$ indexed by distinct $k$-tuples in $\mathbb{N}$ such that $(\pi_\alpha)$ has the same distribution as the subprocess of $(\tilde{\pi}_\phi)$ indexed by tuples in increasing order, where we identify $\{i_1 < i_2 < \cdots < i_k\} \in \binom{\mathbb{N}}{k}$ with $(i_1, i_2, \ldots, i_k) \in \mathrm{Inj}([k], \mathbb{N})$; he then gave a more-or-less complete structure theorem for both classes. He also extended these conclusions to processes indexed by the set $\binom{\mathbb{N}}{<\infty}$ of all finite subsets of $\mathbb{N}$; this amounts to a spreadable coupling of the examples indexed by sets of different fixed sizes.

Spreadable doubleton-indexed processes also appear, under the name of 'deletion-invariant random graphs', in the independent work [29] of Fremlin and Talagrand; in particular, theirs may be the first use of this kind of analysis to answer a more purely 'combinatorial' question, albeit still about an class of infinitary structures. Let us also note in passing that several other names for spreadability have also appeared in the literature, such as 'spreading-invariant' (as in Kingman [48]) and 'contractible' in Kallenberg's recent book [46].

We shall discuss spreadability briefly in Subsection 3.9, but will not examine in detail the modifications it requires to the exchangeability theory, since these lie further from the comparison we wish to make with combinatorics. Kallenberg's paper [44] contains an essentially complete account of this subject.

In the rest of this introductory subsection we will give some flavour of the structural results that follow by considering cases 1 and 2 above.

In case 1, a complete description of exchangeable $\mu$ is a classical consequence of a theorem of de Finetti:

**Theorem 2.5** (Structure of exchangeable random colouring). *Let $S$ be an infinite set, $K$ a standard Borel space and $\mu$ an exchangeable probability on $K^S$. Then $\mu$ can be represented as a mixture of product measures: there exist a*

Lebesgue space $(Z_0, \mu_0)$ and a probability kernel $P_1 : Z_0 \rightsquigarrow K$ such that

$$\mu = \int_{Z_0} P_1^{\otimes S}(z, \cdot) \mu_0(\mathrm{d}z).$$

This result serves as an important and instructive precursor to more general analyses, so we will recall a full proof in Section 3.1 below.

We can rephrase Theorem 2.5 as a 'recipe' for producing all possible exchangeable measures on $K^S$. We first select some Lebesgue space $(Z_0, \mu_0)$ and some probability kernel $P_1 : Z_0 \rightsquigarrow K$, and now a $\mu$-distributed random member of $K^S$ is obtained thus:

- First, choose some $z$ at random from $Z_0$, according to $\mu_0$;
- Then, from the law $P_1(z, \cdot)$ that corresponds (measurably) to this choice of $z$, choose each coordinate $\omega_s$ of $\omega \in K^S$ independently with distribution $P_1(z, \cdot)$.

How might such a recipe for building an exchangeable measure in case $k = 1$ extend to $k = 2$? Certainly, one possibility is to choose each edge-colour in $K$ independently at random from some fixed probability distribution $\nu_0$ on $K$, so $\mu = \nu_0^{\otimes \binom{S}{2}}$: for example, if $K = \{0, 1\}$ and $\nu_0\{1\} = p$ then this gives simply the infinite Erdős-Renyi random graph $\mathbf{G}(S, p)$. Next, just as for the case $k = 1$, we could form a mixture of such measures, so that $\mu = P_\#^{\otimes \binom{S}{2}} \mu_0$ for some Lebesgue space $(Z_0, \mu_0)$ and some $P : Z_0 \rightsquigarrow K$.

However, we are now able also to introduce a third possibility. For simplicity let us first consider a special case. If $G = (V, E)$ is some fixed finite graph, then we may form an exchangeable measure $\mu$ on $\{0, 1\}^{\binom{S}{2}}$ by first picking some vertex $v_s \in V$ for each $s \in S$ independently and uniformly at random, and then including each edge $\{s, t\} \in \binom{S}{2}$ according to whether $\{v_s, v_t\} \in E$. Thus, we take infinitely many samples of the vertices in $V$, indexed by the set $S$, and then join members of $S$ by edges according to the behaviour in $G$ of the corresponding sample vertices, and so form an infinite random expansion of the graph $G$. This certainly has an exchangeable law. This construction gives a fairly simple way to convert any fixed finite graph into an infinite exchangeable random graph by treating the original graph as a source of sample vertices.

We may extend to the case of arbitrary $K$ by now ignoring the edge set $E$ of $G$, and introducing instead a $K$-coloured graph $G : V^k \to K$ with loops. This construction also extends directly to higher-rank hypergraphs; let us make it formal with the following definition.

**Definition 2.6** (Sampling random hypergraph). Suppose that $H$ is a fixed $k$-uniform $K$-coloured hypergraph with loops on a vertex set $V$. Then we define the **$H$-sampling random hypergraph** $\mu_H$ on $K^{\binom{S}{k}}$ to be the law of the random $k$-uniform $K$-coloured hypergraph on $S$ obtained by sampling for each $s \in S$ independently and uniformly at random a vertex $v_s \in V$ and then letting the colour of an edge $e = \{s_1, s_2, \ldots, s_k\} \in \binom{S}{k}$ be $H(v_{s_1}, v_{s_2}, \ldots, v_{s_k})$.

Note that it was in order to allow for collisions in this list $v_{s_1}$, $v_{s_2}$, ..., $v_{s_k}$ of vertices that we required the input coloured hypergraph $H$ to have loops.

It is this construction of an exchangeable random hypergraph from a fixed finite hypergraph that underlies the correspondence principle between statistical properties of fixed finite hypergraph colourings and random infinite hypergraph colourings. We will take up this principle and some of the parallels that it supports between the two settings in Subsections 4.5 and 4.6.

However, the sampling random hypergraphs do not capture all the possible behaviours of an exchangeable random hypergraph: we need to allow for more different sources of randomness. For example, if $k = 3$, then in addition to sampling vertices of a fixed 3-uniform $K$-coloured hypergraph $H$, we could consider an alternative sampling procedure based on a finite set $W$ of possible *edge*-colours and a symmetric function $H : W^3 \to K$, by defining the measure $\theta_H$ to be the law of the random hypergraph obtained by first sampling independently and uniformly from $W$ a colour $w_e$ for each edge $e \in \binom{S}{2}$, and then setting the $K$-colour of $u \in \binom{S}{3}$ to be $H(w_{e_1}, w_{e_2}, w_{e_3})$ when $\binom{u}{2} = \{w_1, w_2, w_3\}$.

For higher ranks $k$, the obvious extension of the above idea yields a random-sampling construction of exchangeable random coloured hypergraphs corresponding to each intermediate rank $i \leq k$. It turns out that these different sampling procedures really do, in general, lead to distinct random hypergraphs — not every 'edge-sampling' random 3-uniform hypergraph as constructed above can be recovered by an alternative 'vertex-sampling' construction — but we postpone giving an example that witnesses this until Subsection 3.6, after the proof of the main Structure Theorem 2.9, as various elements of that proof will ease the analysis of the example.

Any general structure theorem for exchangeable random hypergraphs must allow for all of these different sampling procedures in a single recipe; our main result is that, once formulated appropriately, *this* is essentially all we need. The final structure that emerges extends the above sampling construction in essentially two ways: firstly by introducing a whole tower of sampling procedures, first of vertices, then edges, then 3-hyperedges, and so on up to the full rank $k$; and secondly by allowing each of these sampling procedures its own infusion of 'additional randomness', by replacing deterministic colourings with probability kernels into the relevant space of colours. It is most naturally phrased as an extension of the 'random recipe' we gave above as a formulation of de Finetti's Theorem.

**Definition 2.7.** By a sequence of **ingredients** we understand a sequence $Z_0$, $Z_1$, ..., $Z_{k-1}$ of standard Borel spaces, a probability measure $\mu_0$ on $Z_0$ and (setting $Z_k := K$), for each $i = 1, 2, \ldots, k$, a probability kernel

$$P_i : Z_0 \times Z_1^i \times Z_2^{\binom{i}{2}} \times \cdots \times Z_{i-1}^{\binom{i}{i-1}} \rightsquigarrow Z_i$$

that is symmetric under the natural coordinate-permuting action of $\mathrm{Sym}([i])$ on the domain.

**Definition 2.8.** Given a sequence of ingredients

$$(Z_0, \mu_0), (Z_1, P_1), \ldots, (Z_{k-1}, P_{k-1}), (K, P_k),$$

by the **standard recipe** we understand the following procedure for picking a random member of $K^{\binom{S}{k}}$:

- Choose $z_\emptyset \in Z_0$ at random according to the law $\mu_0$;
- Colour each vertex $s \in S$ by some $z_s \in Z_1$ chosen independently according to $P_1(z_\emptyset, \cdot)$;
- Colour each edge $a = \{s, t\} \in \binom{S}{2}$ by some $z_a \in Z_2$ chosen independently according to $P_2(z_\emptyset, z_s, z_t, \cdot)$;

$$\vdots$$

- Colour each $(k-1)$-hyperedge $u \in \binom{S}{k-1}$ by some $z_u \in Z_{k-1}$ chosen independently according to $P_{k-1}(z_0, (z_s)_{s \in u}, \ldots, (z_v)_{v \in \binom{u}{k-2}}, \cdot)$;
- Colour each $k$-hyperedge $e \in \binom{S}{k}$ by some colour in $K$ chosen independently according to $P_k(z_0, (z_s)_{s \in e}, \ldots, (z_u)_{u \in \binom{e}{k-1}}, \cdot)$.

If $\mu$ is the law of the resulting random $K$-coloured $k$-uniform hypergraph, we will say that the ingredients **yield** $\mu$ upon following the standard recipe.

Symbolically the standard recipe may be represented as the following recursive construction of a family of exchangeable random hypergraphs terminating in $\mu$: $\mu_0$ on $Z_0$ is given, and then we define

$$\mu_1 := (\mathrm{id}_{Z_0}, P_1^{\otimes S})_\# \mu_0 \quad \text{on } Z_0 \times Z_1^S,$$

$$\mu_2 := \left( \mathrm{id}_{Z_0}, \mathrm{id}_{Z_1^S}, \bigotimes_{a \in \binom{S}{2}} P_2 \circ (\mathrm{id}_{Z_0}, \pi_{\binom{a}{1}}) \right)_\# \mu_1 \quad \text{on } Z_0 \times Z_1^S \times Z_2^{\binom{S}{2}},$$

$$\vdots$$

$$\mu_{k-1} := \left( \mathrm{id}_{Z_0}, \mathrm{id}_{Z_1^S}, \ldots, \mathrm{id}_{Z_{k-2}^{\binom{S}{k-2}}}, \bigotimes_{u \in \binom{S}{k-1}} P_{k-1} \circ (\mathrm{id}_{Z_0}, \pi_{\binom{u}{1}}, \ldots, \pi_{\binom{u}{k-2}}) \right)_\# \mu_{k-2}$$

$$\text{on } Z_0 \times Z_1^S \times \cdots \times Z_{k-1}^{\binom{S}{k-1}}$$

$$\mu = \left( \bigotimes_{e \in \binom{S}{k}} P_k \circ (\mathrm{id}_{Z_0}, \pi_{\binom{e}{1}}, \ldots, \pi_{\binom{e}{k-1}}) \right)_\# \mu_{k-1} \quad \text{on } K^{\binom{S}{k}}.$$

**Theorem 2.9** (Structure theorem for uniform exchangeable random hypergraph colourings)**.** *For any $k$-uniform exchangeable random hypergraph $\mu$ there is some sequence of ingredients*

$$(Z_0, \mu_0), (Z_1, P_1), \ldots, (Z_{k-1}, P_{k-1}), (K, P_k)$$

*which yields $\mu$ upon following the standard recipe.*

**Example.** In terms of the structure theorem the $H$-sampling random $\{0,1\}$-coloured hypergraph $\mu_H$ for $H$ a $k$-uniform hypergraph with loops on a finite vertex set $V$ has the following ingredients

$$(\{*\}), (V, P_1^H), (\{*\}, P_2^H), \ldots, (\{*\}, P_{k-1}^H), (\{0,1\}, P_k^H)$$

where we write $\{*\}$ for a one-point space, $P_1^H$ is the uniform distribution on $V$, $P_{i+1}^H$ is identically $\delta_*$ for $i = 0, 1, \ldots, k-2$ and

$$P_k^H(*, (v_1, v_2, \ldots, v_k), *, \ldots, *, \cdot) = \delta_{H(v_1, v_2, \ldots, v_k)}.$$

◁

Let us note here that the forms of the structure theorems for exchangeable random colourings that we have adopted above are not quite the same as those favoured in the earlier probabilistic literature on the subject, such as Aldous [3] and Kallenberg [44]. In particular, while we have allowed the use of rather abstract spaces and probability kernels as structural ingredients for our exchangeable laws $\mu$, it has in the past been popular to exhibit a stochastic process with a representation in terms of simple ingredients of a different kind whose law is $\mu$. For example, the version of Theorem 2.9 that follows most easily from the formalism of Kallenberg [44] reads as follows:

**Theorem.** *For any $k$-uniform exchangeable random $K$-coloured hypergraph $\mu$ there is some measurable function $f : [0,1]^{\binom{[k]}{\leq k}} \to K$ such that, if $(\xi_a)_{a \in \binom{S}{\leq k}}$ is a collection of independent uniform $[0,1]$-valued random variables then the stochastic process $\left(f(\xi_a)_{a \subseteq e}\right)_{e \in \binom{S}{k}}$ has joint law $\mu$.*

We shall discuss these older 'representation' versions of our results for graphs in Subsection 3.2 and for hypergraphs in Subsection 3.3, but in fact the original probabilistic theorems correspond most closely to the setting of random colourings of partite hypergraphs, and so we shall defer most of our discussion of them until we address this setting in Subsection 3.7. This older formalism will also re-appear in Subsection 4.2, where we shall find it paralleled very closely by the results obtained recently by Elek and Szegedy on 'limit object' representations for the statistics of large graphs and hypergraphs: their 'limit objects' can simply be identified with function $f$ above. ast

The forms given for the structure theorems above suggest various questions on what further refinements might be possible. We will show in Subsection 3.6 that in general a simple exchangeable random hypergraph does require type spaces of several ranks in its recipe: it cannot be recovered using only a mixture of vertex-sampling random graphs for some standard Borel probability space of sample vertices.

A slight subtlety in our means for obtaining our structural ingredients is the following. Given an exchangeable random hypergraph colouring $\mu$, we might ask whether it can be that there not only exists some suitable collection of ingredients $(Z_i, P_i)$ yielding $\mu$, but that the probability kernels of these ingredients can be recovered as factors of the original probability space $(K^{\binom{S}{k}}, \mu)$? That is, can we obtain our products of the type spaces $Z_i$ and the probability kernels $P_i$ from some tower of factors

$$(K^{\binom{S}{k}}, \mu) \xrightarrow{\phi_{k-1}} \left(Z_{k-1}^{\binom{S}{k-1}}, \mu_{k-1}\right) \xrightarrow{\phi_{k-2}} \left(Z_{k-2}^{\binom{S}{k-2}}, \mu_{k-2}\right) \xrightarrow{\phi_{k-3}} \cdots$$
$$\xrightarrow{\phi_1} \left(Z_1^S, \mu_1\right) \xrightarrow{\phi_0} (Z_0, \mu_0),$$

where (setting $\mu_k := \mu$) the image measure

$$\mu_i := (\phi_1)_\# \cdots (\phi_{k-2})_\# (\phi_{k-1})_\# \mu$$

is precisely the exchangeable random $i$-uniform hypergraph $Z_i$-colouring obtained by stopping our recipe at level $i$?

This is possible in the case of de Finetti's Theorem 3.1 (that is, $k = 1$), where one option is to obtain relative independence over the $(\mathrm{Sym}_0(S)$-invariant) 'tail $\sigma$-algebra' of the original system. However, it turns out to be impossible in general if $k \geq 2$: in Subsection 3.6 we will also present an example of an exchangeable random graph for which we must use unrecoverable quasifactors.

## 2.3. Relations to combinatorics: correspondence principles and limit objects

There are certain theorems of combinatorics that can be conveniently approached by first converting the original combinatorial data into a related kind of stochastic process, and then applying to that process the more analytic methods of probability or ergodic theory. The successes of this approach are perhaps most striking in arithmetic combinatorics, and in particular in the ergodic-theoretic approach to Szemerédi's Theorem discovered by Furstenberg in [30] (for an overview of these relationships, see also Bergelson [15] and Tao and Vu [64]). Indeed, several more general results of density Ramsey theory in the arithmetic or related settings are still known only by such methods.

The machinery of exchangeable random hypergraph colourings stands in a similar relation to finitary hypergraph theory as does the ergodic theory of $\mathbb{Z}^d$-actions to arithmetic combinatorics. While it is not clear that this 'correspondence principle' can have such powerful consequences for finitary hypergraph theory as for its arithmetic counterpart, there are certain kinds of combinatorial question for which parallel versions can be set up for exchangeable random hypergraph colourings which may then either be instructive or interesting in their own right. We shall give a brief overview of these parallels here, and will then examine two particular examples, coming from property testing and extremal combinatorics, in Subsections 4.5 and 4.6.

Let us first describe informally how this 'correspondence principle' might arise. Suppose we are given a combinatorial structure $T$ of very large finite size and a dense subset $E \subseteq T$. Suppose further that the structure $T$ has a natural group $\Gamma$ of symmetries, and that we are interested in the 'statistics' of the subset $E$ that count the number of copies of a fixed (small) pattern $L \subseteq T$ lying inside $E$: that is, the number of $g \in \Gamma$ with $g(L) \subseteq E$. Then we can form from $E$ an associated random subset $\mathbf{E}$ of $X$, by setting $\mathbf{E} := \mathbf{g}(E)$ for a symmetry $\mathbf{g} \in \Gamma$ chosen uniformly at random; letting $\mathsf{P}$ be the law of $\mathbf{E}$ in $\binom{T}{\leq |T|}$, it is clear that this is $\Gamma$-invariant for the canonical action of $\Gamma$ on $\binom{T}{\leq |T|}$. For each $t \in T$, we can consider the indicator function of the event $A_t := \{\mathbf{E} \ni t\}$; evaluation of all products of these indicator functions for different $t \in T$ describes completely the law $\mathsf{P}$.

The count we wish to make of those $g$ such that $g \in \Gamma$, upon renormalizing by $|\Gamma|$, is now just the probability $\mathsf{P}(\mathbf{E} \supseteq L)$. It may be that some of the properties of these probabilities can be established by studying the possible structure of general such $\Gamma$-invariant laws $\mathsf{P}$ on $\binom{T}{\leq |T|}$; in particular, sometimes the asymptotics of these properties may be more easily analyzed by considering a nested sequence $T_1 \subseteq T_2 \subseteq \ldots$ of larger and larger index sets, and then defining from the associated invariant laws $\mathsf{P}_n$ some vague limit $\mathsf{P}$ on $\binom{T_\infty}{<\infty}$ for a suitable countably infinite index set $T_\infty$.

The correspondence principle of interest to us in the coloured hypergraph setting corresponds to the above outline in the case that $T$ is the complete $k$-uniform hypergraph on some large finite vertex set $S_n$ and $\Gamma$ is the group of automorphisms of $T$ given by vertex set permutations; after taking a vague limit as above, we are naturally left with an exchangeable random $k$-uniform hypergraph (equivalently, a $\{0,1\}$-coloured hypergraph) on the countably infinite vertex set $S := \bigcup_{n \geq 1} S_n$.

It is this relation that leads to the overlap between counting substructures of finite coloured hypergraph and exchangeable random coloured hypergraphs. In view of this it is only to be expected that various properties relating to the first-order statistics of these finite objects have analogs for their infinite random counterparts, and that arguments relating to these properties are often skeletally similar in the two categories.

In the setting of arithmetic combinatorics, a slightly more delicate such correspondence principle underlies Furstenberg's analysis of densities of arithmetic progressions in a dense subset of $\{1, 2, \ldots, N\}$. These are translated into probabilities of multiple recurrence in probability-preserving systems, the analysis of which then leads to the ergodic-theoretic proof of Szemerédi's Theorem. Of course, in Furstenberg's case the translation to an infinitary result made the proof of the finitary result much more manageable through efficiently hiding a large number of intricate estimates required in a finitary version of the argument (consider, in this relation, Tao's finitarization of the proof in [63]). This situation is surely rather rare; in the hypergraph setting, it seems more common that where parallel versions of theorems are available, their proofs are of comparable sophistication on either side of the divide.

We should stress that many finitary questions do not have a sensible infinitary analog, since in order to do so they must take a sensible 'asymptotic' form as the orders of the hypergraphs in question tend to $\infty$. Among those that do exhibit this form, the infinitary version of the question typically bears only on the finitary version up to 'leading-order'. For such questions, we are often able to sharpen the finitary version to ask for higher-order information — perhaps by asking for estimates on rates of convergence, as opposed to the mere fact of convergence — that has no clear infinitary analog: $\mathrm{O}(\sqrt{n})$ is more than just $\mathrm{o}(n)$. Indeed, in some cases it is only after raising the stakes in this way that the finitary question seems interesting, and so the infinitary version may appear somewhat degenerate. For example, it is hard to imagine a modification of Theorem 2.9 that can be brought to bear on the finer Turán problem for bipartite graphs (for which the leading-order density is 0, and this is already the answer to which the infinitary approach pertains) or on the problem of Zarankiewicz (see, for example, Chapter IV of Bollobás [16]).

Two specific arenas that do allow a parallel study of finite hypergraphs and exchangeable random hypergraphs are questions of efficient property testing (to be considered in Subsection 4.5, and at much greater length in the forthcoming paper [10]) and extremal questions (to be considered in Subsection 4.6). To date, there are some early indications that this parallelism may occasionally bear fruit in the study of finitary questions, such as recent results of Razborov [54, 55] on an instance of the graph-copy problem, although it is surely too early to judge whether a trend may emerge.

Insofar as the structure theory for exchangeable random coloured hypergraphs is of use in proving results in the above settings, in principle the ways in which it can be used should be replicable purely in terms of counting substructures of the original finitary graphs or hypergraphs. The benefit to be derived from an examination of the infinitary latter machinery may lie largely in the more convenient organization of such finitary arguments in the future, rather than in forging a new way actually to write out proofs. In place of the infinitary structural analysis of exchangeable random hypergraphs, purely finitary arguments tend to rest heavily on Szemerédi's graph regularity lemma (see, for example, Sections IV.5 and IV.6 of Bollobás [16]) and its later hypergraph versions developed by Nagle, Rödl and Schacht [53] and by Gowers [36]; these serve as a kind of 'structure theory' in the finitary category, and show some interesting similarities with the exchangeability theory (although deriving either from the other does not seem quite straightforward). We shall discuss this comparison a little further in 4.4.

## 3. Exchangeable families of random variables

We give an analysis of many of the instances of exchangeability listed above, and provide some examples to illustrate different phenomena that can arise in the process. We have tried to arrange our particular examples in roughly increasing order of generality and sophistication.

### 3.1. Warmup: de Finetti's Theorem

We first recall a version of de Finetti's Theorem ([18, 19]) and its use to deduce a structure theorem for the special case $k = 1$ of the setting of random hyper-graph colourings. This serves as helpful motivation for later, more complicated instances.

Since our vertex set $S$ is infinite, it can be partitioned into two infinite subsets, say as $S = S_1 \cup S_2$. Now observe that within the canonical family of random variables $(\pi_s)_{s \in S}$, the subfamily $(\pi_s)_{s \in S_1}$ has the same joint law under any bijective identification of $S$ with $S_1$, say $\psi : S_1 \to S$. Indeed, it suffices to check this assertion for the finite-dimensional projections: if $s_1, s_2, \ldots, s_k \in S_1$ and $A_1, A_2, \ldots, A_k \in \Sigma_K$ then

$$\begin{aligned}
&(\pi_{S_1})_\# \mu \{\pi_{s_1} \in A_1, \pi_{s_2} \in A_2, \ldots, \pi_{s_k} \in A_k\} \\
&\stackrel{\text{def}}{=} \mu\{\pi_{s_1} \in A_1, \pi_{s_2} \in A_2, \ldots, \pi_{s_k} \in A_k\} \\
&= \mu\{\pi_{\psi(s_1)} \in A_1, \pi_{\psi(s_2)} \in A_2, \ldots, \pi_{\psi(s_k)} \in A_k\} \quad \text{by exchangeability} \\
&= (\tau^\psi)_\# \mu \{\pi_{s_1} \in A_1, \pi_{s_2} \in A_2, \ldots, \pi_{s_k} \in A_k\},
\end{aligned}$$

as required.

Therefore it suffices to prove our structure theorem for the projection $(\pi_{S_1})_\# \mu$ instead of for the original measure $\mu$. This sleight of hand now allows us to use the remaining random variables $(\pi_s)_{s \in S_2}$ as 'poll vertices': we shall deduce the required form for the joint law of $(\pi_{S_1})_\# \mu$ by showing that the $(\pi_s)_{s \in S_1}$ are conditionally independent given the behaviour of the $(\pi_s)_{s \in S_2}$.

It turns out that with a little more care it is possible to do without this decomposition of $S$, effectively by letting $S_2 \supset S_3 \supset \cdots$ be a descending sequence of infinite-coinfinite subsets of $S$ with empty intersection, and conditioning only on the tail $\sigma$-algebra $\bigcap_{i \geq 2} \sigma(\pi_{S_i})$: it turns out that this retains 'enough information' about the joint behaviour of the individual $\pi_s$ that we can still prove the necessary conditional independence. However, the analog of this argument fails in the more complex cases of exchangeability to be considered later (we will justify this with an example in Subsection 3.6), and so we shall not pursue it further here.

Thus, writing $Z_0 := K^{S_2}$ and letting $\mu_0 \in \Pr Z_0$ be the projection $(\pi_{S_2})_\# \mu$ of $\mu$ onto the auxiliary coordinates, these will themselves be our structural ingredients.

**Theorem 3.1** (Variant of de Finetti's Theorem). *The variables $\pi_s$, $s \in S_1$, are relatively independent over $\pi_{S_2}$ and are such that for any $f \in L^\infty(K)$ the conditional expectations $\mathsf{E}_\mu[f \circ \pi_s \,|\, \pi_{S_2}]$ are (almost surely) equal for all $s \in S_1$.*

*Proof.* An essentially analogous proof but using the abovementioned tail $\sigma$-algebra can be found under Theorem 11.10 in Kallenberg [45]. Suppose that $F = \{s_1, s_2, \ldots, s_r\} \subset S_1$ and that we are given $f_1, f_2, \ldots, f_r \in L^\infty(K)$. We need to show that

$$\mathsf{E}_\mu\left[\prod_{i \leq r} f_i \circ \pi_{s_i} \,\bigg|\, \pi_{S_2}\right] = \prod_{i \leq r} \mathsf{E}_\mu[f_i \circ \pi_{s_i} \,|\, \pi_{S_2}].$$

By induction on $r$ it clearly suffices to prove that

$$\mathsf{E}_\mu\left[\prod_{i\leq r} f_i \circ \pi_{s_i} \,\Big|\, \pi_{S_2}\right] = [f_i \circ \pi_{s_i} \,|\, \pi_{S_2}] \cdot \mathsf{E}_\mu\left[\prod_{i=2}^r f_i \circ \pi_{s_i} \,\Big|\, \pi_{S_2}\right],$$

and hence that for any $A \in \sigma(\pi_{S_2})$ we have

$$\int_A \prod_{i\leq r} f_i \circ \pi_{s_i} \, \mathrm{d}\mu = \int_A \mathsf{E}_\mu[f_i \circ \pi_{s_i} \,|\, \pi_{S_2}] \cdot \left(\prod_{i=2}^r f_i \circ \pi_{s_i}\right) \mathrm{d}\mu.$$

By the approximability in $\mu$ of any Borel subset of $K^{S_2}$ by a finite-dimensional cylinder set, we may assume that $A$ depends only on coordinates in some finite $E \subset S_2$. Therefore, choosing coordinates $s'_2, s'_3, \ldots, s'_r \in S_2 \setminus E$, by exchangeability we know that

$$\int_A \prod_{i\leq r} f_i \circ \pi_{s_i} \, \mathrm{d}\mu = \int_A f_1 \circ \pi_{s_1} \cdot \left(\prod_{i=2}^r f_i \circ \pi_{s'_i}\right) \mathrm{d}\mu,$$

and now both $A$ and the functions $f_i \circ \pi_{s'_i}$, $i = 2, 3, \ldots, r$, depend only on the coordinates $\pi_s$ for $s \in S_2$, so certainly we may replace $f_1 \circ \pi_{s_1}$ with $\mathsf{E}_\mu[f_1 \circ \pi_{s_1} \,|\, \pi_{S_2}]$ in the right-hand-side above:

$$\int_A \prod_{i\leq r} f_i \circ \pi_{s_i} \, \mathrm{d}\mu = \int_A \mathsf{E}_\mu[f_1 \circ \pi_{s_1} \,|\, \pi_{S_2}] \cdot \left(\prod_{i=2}^r f_i \circ \pi_{s'_i}\right) \mathrm{d}\mu.$$

This is almost what we want: it remains only to move the vertices $s'_i$ back to $s_i$ to leave the desired equality. However, this requires just a little care: while we certainly still have exchangeability to appeal to, we must make sure that the conditional expectation $\mathsf{E}_\mu[f_1 \circ \pi_{s_1} \,|\, \pi_{S_2}]$ is actually independent of the coordinates $\pi_{s'_2}, \pi_{s'_3}, \ldots, \pi_{s'_r}$, so that upon swapping them back with $\pi_{s_2}, \pi_{s_3}, \ldots, \pi_{s_r}$ this conditional expectation is unchanged.

This is so because for any infinite $S_3 \subset S_2 \setminus (E \cup \{s'_2, s'_3, \ldots, s'_r\})$, we know that $(\pi_{S_1 \cup S_3})_\#\mu$ is just another copy of $\mu$, and hence, in particular, that

$$\left\|\mathsf{E}_\mu[f_1 \circ \pi_{s_1} \,|\, \pi_{S_3}]\right\|_2^2 = \left\|\mathsf{E}_\mu[f_1 \circ \pi_{s_1} \,|\, \pi_{S_2}]\right\|_2^2.$$

Since, on the other hand,

$$\mathsf{E}_\mu[f_1 \circ \pi_{s_1} \,|\, \pi_{S_3}] = \mathsf{E}_\mu\big[\mathsf{E}_\mu[f_1 \circ \pi_{s_1} \,|\, \pi_{S_2}] \,\big|\, \pi_{S_3}\big],$$

this equality of $L^2$-norms can hold only if the two conditional expectations are actually equal:

$$\mathsf{E}_\mu[f_1 \circ \pi_{s_1} \,|\, \pi_{S_3}] = \mathsf{E}_\mu[f_1 \circ \pi_{s_1} \,|\, \pi_{S_2}].$$

We will sometimes refer back to this as the 'tail property' of these conditional expectations.

It follows, in particular, that $\mathsf{E}_\mu[f_1 \circ \pi_{s_1} \,|\, \pi_{S_2}]$ does not depend on $\pi_{s'_2}, \pi_{s'_3}, \ldots, \pi_{s'_r}$; so it is invariant under the coordinate-permutation that swaps $(s_2, s_3, \ldots, s_r)$ with $(s'_2, s'_3, \ldots, s'_r)$, and we are left with the desired equality. □

**Remark.** It is clearly essential for this proof that $S$ be infinite, and in fact it is not hard to construct examples of measures enjoying the analogous symmetry with finite $S$ for which the conclusion fails. ◁

We can now deduce Theorem 2.5 immediately.

*Proof of Theorem 2.5.* Conditioned on $\pi_{S_2} : K^S \to Z_0$, the measure $(\pi_{S_1})_\#\mu$ disintegrates as an integral of fibre measures; and by the above these are all product measures, giving the desired form

$$(\pi_{S_1})_\#\mu = \int_{Z_0} P_1^{\otimes S_1}(z, \,\cdot\,)\, \mu_0(\mathrm{d}z),$$

where for $P_1 : Z_0 \rightsquigarrow K$ we simply take the conditional law of any individual one of the random variables $\pi_s$ for $s \in S_1$ over $\pi_{S_2}$. By the equivalence between $\mu$ and $(\pi_{S_1})_\#\mu$ described previously, this completes the proof. □

We will later return repeatedly to the above trick of considering the structure of the projected measure $(\pi_{T_1})_\#\mu$ instead of that of $\mu$ itself for some suitable infinite subset $T_1$ of the index set $T$, in order to set at our disposal an additional infinite pool of random variables $(\pi_t)_{t \in T \setminus T_1}$ corresponding to 'reference vertices'. In the coloured hypergraph setting with $T = \binom{S}{k}$ we will take $T_1 := \binom{S_1}{k}$ for some infinite-coinfinite $S_1 \subset S$, and will then refer to $(\pi_{\binom{S_1}{k}})_\#\mu$ as the **induced random hypergraph on** $S_1$, by analogy with the usual combinatorial term for an induced subhypergraph. This idea of selecting some infinite subset of vertices to 'poll' when analyzing the part played by other vertices appears explicitly in Tao's approach to an infinitary hypergraph removal lemma ([62]).

### 3.2. Exchangeable random graph colourings

We will now consider the slightly more complicated setting of exchangeability for random graph colourings. Naturally, it will be subsumed by the treatment of hypergraphs of arbitrary rank in the next section.

To prove Theorem 2.9 in the special case $k = 2$ we must extract from a given exchangeable random $K$-colouring of $\binom{S}{2}$ an exchangeable random $Z_1$-colouring of $S$ for some suitable auxiliary colour space $Z_1$. As in the case of de Finetti's Theorem, this $Z_1$ will emerge explicitly from a decomposition of $K^{\binom{S}{2}}$ corresponding to a partition of $S$ into two infinite subsets, $S = S_1 \cup S_2$. Our first step, then, is to write

$$K^{\binom{S}{2}} = K^{\binom{S_1}{2}} \times K^{S_1 \times S_2} \times K^{\binom{S_2}{2}} = K^{\binom{S_1}{2}} \times (K^{S_2})^{S_1} \times K^{\binom{S_2}{2}},$$

where we identify $K^{\{e \in \binom{S}{2}: \, |e \cap S_1| = |e \cap S_2| = 1\}}$ with $K^{S_1 \times S_2}$ and then with $(K^{S_2})^{S_1}$ in the natural way.

Now let $Z_1 := K^{S_2} \times K^{\binom{S_2}{2}}$, and for $s \in S_1$ write $\pi_s^{Z_1}$ for the composition of the projection $\pi_{\binom{\{s\} \times S_2}{2}} : K^{\binom{S}{2}} \to K^{\{e \in \binom{\{s\} \cup S_2}{2}: \, e \ni s\}} \times K^{\binom{S_2}{2}}$ with the natural identification of the first of these factors with $Z_1$. It is clear that $\pi_s^{Z_1}$ is invariant under $\tau^g$ for any $g \in \mathrm{Sym}_0(S)$ that fixes $\{s\} \cup S_2$, that if $s \in S_1$ and $g \in \mathrm{Sym}_0(S_1)$ then $\pi_s^{Z_1} \circ \tau^g = \pi_{g(s)}^{Z_1}$, and that the $\sigma$-subalgebra $\bigvee_{s \in S_1} \sigma(\pi_s^{Z_1})$ is $\mathrm{Sym}_0(S_1)$-invariant. In addition, for any $f \in L^\infty(K)$ we have

$$\mathsf{E}_\mu[f \circ \pi_{e_1} \mid \pi_s^{Z_1}] = \mathsf{E}_\mu[f \circ \pi_{e_2} \circ \tau^g \mid \pi_s^{Z_1}]$$

whenever $e_1, e_2 \in \binom{S_1}{2}$ share the vertex $s$ and $g \in \mathrm{Sym}_0(S_1)$ fixes $s$ and sends $e_2$ onto $e_1$: this follows from the symmetry and the fact that $\pi_s^{Z_1}$ is fixed by any such permutation $g$.

**Lemma 3.2.** *The random variables $(\pi_s^{Z_1})_{s \in S_1}$ are exchangeable.*

*Proof.* This is a direct verification: if $s_1, s_2, \ldots, s_k \in S_1$, $g \in \mathrm{Sym}_0(S_1)$ and $A_i \in \sigma(\pi_{s_i}^{Z_1})$ for each $i \leq k$, then

$$\mu\left(\bigcap_{i \leq k} A_i\right) = \mu\left(\bigcap_{i \leq k} \tau^g(A_i)\right)$$

(by the permutation invariance of $\mu$); but if $A_i = (\pi_{s_i}^{Z_1})^{-1}(B_i)$ with $B_i \in \Sigma_{Z_1}$ for each $i \leq k$ then the sets involved in the two sides of this equation are precisely $\bigcap_{i \leq k} \{\pi_{s_i}^{Z_1} \in B_i\}$ and $\bigcap_{i \leq k} \{\pi_{g(s_i)}^{Z_1} \in B_i\}$ respectively. $\square$

Now we turn to the main estimates that will underly our structure theorem. For convenience, let us write $\mathrm{T}(s)$ for the $\sigma$-subalgebra $(\pi_s^{Z_1})^{-1}(\Sigma_{Z_1}) \subseteq \Sigma_K^{\otimes \binom{S}{2}}$, and set $\mathrm{T} := \mathrm{T}^{S_2} := \bigvee_{s \in S_1} \mathrm{T}(s)$.

Note that we could repeat the whole of the preceding discussion but defining the maps $\pi_s^{Z_1, S_3}$ using only $s$ and vertices inside an infinite subset $S_3 \subseteq S_2$, rather than all of $S_2$; we will need to have this modified construction at our disposal also, and so write $\mathrm{T}^{S_3}$ for the $\sigma$-subalgebra obtained just as was $\mathrm{T}^{S_2}$ but from these more restricted factor maps $\pi_s^{Z_1, S_3}$ (using only the smaller infinite subset $S_3 \subset S_2$ as the set of auxiliary vertices).

**Lemma 3.3** (Tail property of conditional expectation). *Whenever $S_3 \subseteq S_2$ is a further infinite subset, $e \in \binom{S_1}{2}$ and $f \in L^\infty(K)$, the conditional expectation $\mathsf{E}_\mu[f \circ \pi_e \mid \mathrm{T}]$ is actually measurable with respect to $\mathrm{T}^{S_3}$.*

**Remark.** This assertion is analogous to the 'tail property' of the conditional expectations that we used in the last step of the proof of de Finetti's Theorem 3.1. However, whereas in that case this property is actually enough to justify the use of the tail $\sigma$-algebra (as we remarked after the statement of the theorem), here this next step fails: the point is that we would now need to take an intersection of the form $\bigcap_{i \geq 2} \mathrm{T}^{S_i} \vee \mathrm{T}^{S_i}$ for some decreasing sequence $S_2 \supseteq S_3 \supseteq \ldots$ of infinite sets with empty intersection, but we cannot in general then equate this with $\left(\bigcap_{i \geq 2} \mathrm{T}^{S_i}\right) \vee \left(\bigcap_{i \geq 2} \mathrm{T}^{S_i}\right)$ (even up to $\mu$-negligible sets), as we would need to do in order to work with a separate tail $\sigma$-algebra for each vertex $s \in S$. In Subsection 3.6 we will see an example that witnesses this. ◁

*Proof.* This follows from exchangeability and a simple 'energy squeeze' argument. Since $T^{S_3} \subseteq T^{S_2}$, by the iterability of conditional expectation we certainly have

$$\mathsf{E}_\mu[f \circ \pi_e \,|\, T^{S_3}] = \mathsf{E}_\mu\big[\mathsf{E}_\mu[f \circ \pi_e \,|\, T^{S_2}] \,\big|\, T^{S_3}\big],$$

and so, in particular, $\mathsf{E}_\mu[f \circ \pi_e \,|\, T^{S_3}]$ is the image of $\mathsf{E}_\mu[f \circ \pi_e \,|\, T^{S_2}]$ under an orthogonal projection in $L^2(\mu)$. It will therefore follow that they are equal if they have the same norms in $L^2(\mu)$. However, we know that if $e = \{s, t\}$ then

$$\big\|\mathsf{E}_\mu[f \circ \pi_e \,|\, T^{S_2}]\big\|_2^2 = \sup_{F \subset S_2 \text{ finite}} \big\|\mathsf{E}_\mu[f \circ \pi_e \,|\, \sigma((\pi_a)_{a \in \binom{\{s\} \cup F}{2}}) \vee \sigma((\pi_a)_{a \in \binom{\{t\} \cup F}{2}})]\big\|_2^2,$$

and in the expressions within this supremum the value of

$$\big\|\mathsf{E}_\mu[f \circ \pi_e \,|\, \sigma((\pi_a)_{a \in \binom{\{s\} \cup F}{2}}) \vee \sigma((\pi_a)_{a \in \binom{\{t\} \cup F}{2}})]\big\|_2^2$$

depends only on $|F|$, since by exchangeability we may choose a permutation of $S$ that fixes $S_1$ and swaps $F$ with any other subset of $S_2$ of the same size, and this will act as an isometry on $L^2(\mu)$. Since $S_3$ is still infinite, and so contains arbitrarily large finite subsets, it follows that also

$$\begin{aligned}
\big\|\mathsf{E}_\mu[f \circ \pi_e \,|\, T^{S_2}]\big\|_2^2 &\geq \big\|\mathsf{E}_\mu[f \circ \pi_e \,|\, T^{S_3}]\big\|_2^2 \\
&\geq \sup_{F \subset S_2 \text{ finite}} \big\|\mathsf{E}_\mu[f \circ \pi_e \,|\, \sigma((\pi_a)_{a \in \binom{\{s\} \cup F}{2}}) \vee \sigma((\pi_a)_{a \in \binom{\{t\} \cup F}{2}})]\big\|_2^2 \\
&= \big\|\mathsf{E}_\mu[f \circ \pi_e \,|\, T^{S_2}]\big\|_2^2,
\end{aligned}$$

and so the desired energies must be equal and the proof is complete. $\square$

**Proposition 3.4** (Relative independence)**.** *The random variables $\pi_e$ for $e \in \binom{S_1}{2}$ are relatively independent over $T$: for any finite $F \subset \binom{S_1}{2}$ and bounded measurable functions $f_e : K \to \mathbb{R}$ for $e \in F$ we have*

$$\mathsf{E}_\mu\bigg[\prod_{e \in F} f_e \circ \pi_e \,\bigg|\, T\bigg] = \prod_{e \in F} \mathsf{E}_\mu[f_e \circ \pi_e \,|\, T].$$

*Proof.* Enumerate $F$ as $\{e_1, e_2, \ldots, e_r\}$. By induction on $r$ it suffices to show that

$$\mathsf{E}_\mu\bigg[\prod_{p=1}^r f_{e_p} \circ \pi_{e_p} \,\bigg|\, T\bigg] = \mathsf{E}_\mu[f_{e_1} \circ \pi_{e_1} \,|\, T] \cdot \mathsf{E}_\mu\bigg[\prod_{p=2}^r f_{e_p} \circ \pi_{e_p} \,\bigg|\, T\bigg].$$

This identity asserts that whenever $A \in T$ we have

$$\int_A \prod_{e \in F} f_e \circ \pi_e \,\mathrm{d}\mu = \int_A \mathsf{E}_\mu[f_{e_1} \circ \pi_{e_1} \,|\, T] \cdot \mathsf{E}_\mu\bigg[\prod_{p=2}^r f_{e_p} \circ \pi_{e_p} \,\bigg|\, T\bigg] \,\mathrm{d}\mu,$$

and so, since both $A$ and $\mathsf{E}_\mu[f_{e_1} \circ \pi_{e_1} \mid \mathrm{T}]$ are T-measurable, this will follow if we only prove that

$$\int_A \prod_{e \in F} f_e \circ \pi_e \, \mathrm{d}\mu = \int_A \mathsf{E}_\mu[f_{e_1} \circ \pi_{e_1} \mid \mathrm{T}] \cdot \prod_{p=2}^r f_{e_p} \circ \pi_{e_p} \, \mathrm{d}\mu.$$

By continuity of expectation, we may restrict our attention to the case of $A$ a finite-dimensional cylinder, say $A = A_1 \cap A_2 \cap \cdots \cap A_m$ for $A_j \in \sigma\big((\pi_e)_{e \in \binom{\{y_j\} \cup R_j}{2}}\big)$ for some $y_j \in S_1$ and finite $R_j \subset S_2$, $1 \le j \le m$.

Let $e_1 = \{s_1, s_2\}$ and enumerate all the remaining vertices in $(\{y_1, y_2, \ldots, y_m\} \cup e_1 \cup e_2 \cup \cdots \cup e_r) \setminus e_i$ as $z_1, z_2, \ldots, z_q$. Select distinct vertices $w_1, w_2, \ldots, w_q \in S_2 \setminus (R_1 \cup R_2 \cup \cdots \cup R_m)$ and let $g \in \mathrm{Sym}_0(S)$ be the permutation that swaps $z_i$ and $w_i$ for all $i \le q$.

By the exchangeability we know that

$$\int_A \prod_{e \in F} f_e \circ \pi_e \, \mathrm{d}\mu = \int_{\tau^g(A)} \prod_{e \in F} f_{g(e)} \circ \pi_{g(e)} \, \mathrm{d}\mu = \int_{\tau^g(A)} f_{e_1} \circ \pi_{e_1} \cdot \prod_{p=2}^r f_{g(e_p)} \circ \pi_{g(e_p)} \, \mathrm{d}\mu,$$

since $g$ does not move $s_1$ or $s_2$. But now both $\tau^g(A)$ and $\prod_{p=2}^r f_{g(e_p)} \circ \pi_{g(e_p)}$ are clearly T-measurable, since through the action of $g$ we have moved all relevant vertices from $S_1$ into $S_2$, except for $s_1$ and $s_2$; in particular, every edge $g(e_j)$ for $j \ne i$ now has at least one end-point in $S_2$. Therefore the left-hand-side above must equal

$$\int_{\tau^g(A)} \mathsf{E}_\mu[f_{e_1} \circ \pi_{e_1} \mid \mathrm{T}] \cdot \prod_{p=2}^r f_{g(e_p)} \circ \pi_{g(e_p)} \, \mathrm{d}\mu.$$

It remains to put the vertices we have just moved back to their original positions. As in the proof of de Finetti's Theorem, we cannot quite naïvely re-apply the permutation $g$, since the $\sigma$-algebra $\mathrm{T} = \mathrm{T}^{S_2}$ is not fixed under $g$. However, it turns out that the conditional expectation $\mathsf{E}_\mu[f_{e_1} \circ \pi_{e_1} \mid \mathrm{T}]$ does enjoy this invariance; for the infinite subset $S_3 := S_2 \setminus (R_1 \cup R_2 \cup \cdots \cup R_m \cup \{w_1, w_2, \ldots, w_q\})$ (and so also the $\sigma$-algebra $\mathrm{T}^{S_3}$) is invariant under $g$, and by Lemma 3.3 we know that we may replace $\mathrm{T}^{S_2}$ with $\mathrm{T}^{S_3}$ in the conditional expectation of interest. Therefore re-applying $g$ after making this replacement shows that the above integral is equal to

$$\int_A \mathsf{E}_\mu[f_{e_1} \circ \pi_{e_1} \mid \mathrm{T}^{S_3}] \cdot \prod_{p=2}^r f_{e_p} \circ \pi_{e_p} \, \mathrm{d}\mu.$$

Now switching back from $\mathrm{T}^{S_3}$ to $\mathrm{T}^{S_2}$ completes the proof. □

**Corollary 3.5** (Structure theorem for random coloured graphs). *For any exchangeable random coloured graph $\mu$ there are ingredients $(Z_0, \mu_0), (Z_1, P_1), (K, P_2)$ which yield $\mu$ upon following the standard recipe.*

*Proof.* As remarked previously, we may easily find an isomorphism between $\mu$ and the infinite induced random coloured subgraph $(\pi_{\binom{S_1}{2}})_\#\mu$, and so we may prove the result for the subgraph and then simply pull the whole structure back under composition with some bijection $\psi : S \to S_1$. However, Lemma 3.2 and the argument of Proposition 3.4 show that our projections $(\pi_s^{Z_1})_{s \in S_1}$ give the desired relative independence properties for the edge colours in $S_1$ and that their joint law $\mu_1$ is an exchangeable law on $Z_1^{S_1}$. Now de Finetti's Theorem 3.1 allows us to disintegrate the exchangeable measure $\mu_1$ further over some other Lebesgue space $(Z_0, \mu_0)$, as described, completing the proof. □

We note here that the form of our Structure Theorem 2.9 is somewhat different from that appearing in the original probabilistic literature on exchangeability. For example, the analog of Theorem 3.5 that follows most directly from a suitable specialization of Theorem 14.11 in Aldous [3] is the following.

**Theorem 3.6.** *For any exchangeable random $K$-coloured graph $\mu$ there is some measurable function $f : [0,1]^4 \to K$ such that given any uniform $[0,1]$-valued random variable $\xi_\emptyset$ and uniform $[0,1]$-valued processes $(\xi_s)_{s \in S}$ and $(\xi_e)_{e \in \binom{S}{2}}$, all independent, then the stochastic process $\big(f(\xi_\emptyset, \xi_s, \xi_t, \xi_{\{s,t\}})\big)_{\{s,t\} \in \binom{S}{2}}$ has joint law $\mu$.*

Thus, this result represents the joint law $\mu$ using the rather different ingredient of a measurable function $[0,1]^{\binom{[2]}{\leq 2}} \to K$. However, it is relatively easy to move between these two versions of the theorems using standard results on replacing probability kernels with dependence on additional independent random variables; we return to these matters in more detail and exhibit the necessary techniques in Subsections 3.7 and 4.2.

### 3.3. Exchangeable random hypergraph colourings

We will now turn to the full version of Theorem 2.9. We begin by introducing some more notation; this will not only be important for the organization of the present subsection, but will recur in the proof of the positive results of [10].

The skeleton of our proof is largely as in the case of graphs; however, the reduction to de Finetti's Theorem for the last step of the proof of Corollary 3.5 now grows into an induction on the rank $k$. Importantly, even if we begin our analysis in the case of a uniform exchangeable random hypergraph colouring, say $\mu \in K^{\binom{S}{k}}$ for some fixed $k$, the inductive step of the argument will generally introduce a *non*-uniform rank-$(k-1)$ hypergraph colouring on a space of the form $Z_0 \times Z_1^S \times \cdots \times Z_{k-1}^{\binom{S}{k-1}}$ for some nontrivial auxiliary standard Borel spaces $Z_i$.

It is simplest to work in this larger category from the beginning, and so let us first make the obvious extensions of Definitions 2.1 and 2.2 to this purpose:

- For $k \geq 1$, by a $k$-**palette** we shall understand a sequence $K = (K_i)_{i=0}^k$ of standard Borel spaces (not necessarily finite sets), and shall refer to $k$

  as the **rank** of the palette as in the finite case of Definition 2.1
- Given a vertex set $S$ and a general $k$-palette $K = (K_i)_{i=0}^k$ we shall define a **$K$-coloured hypergraph on** $S$ to be a sequence $H_i$ of $K_i$-coloured $i$-uniform hypergraphs on $S$.

Let us also adopt in this subsection the notation

$$K^{(S)} := K_0 \times K_1^S \times K_2^{\binom{S}{2}} \times \cdots \times K_k^{\binom{S}{k}}$$

for the space of all $K$-coloured hypergraphs on $S$.

**Definition 3.7** (Exchangeable random coloured hypergraphs). Given a $k$-palette $K$, by an **exchangeable random $K$-coloured hypergraph on** $S$ we understand a Borel probability measure $\mu$ on $K^{(S)}$ that is invariant under the coordinate-permuting action $\tau$ of $\mathrm{Sym}_0(S)$ defined by

$$\tau^g\big((\omega_e)_{e \in \binom{S}{\leq k}}\big) = (\omega_{g(e)})_{e \in \binom{S}{\leq k}}.$$

It is clear that for each $i \leq k$ the projection of $\mu$ onto the factor $K_i^{\binom{S}{i}}$ is an exchangeable random $i$-uniform $K_i$-coloured hypergraph $\mu_i$; the overall measure $\mu$ may simply be regarded as a $\mathrm{Sym}_0(S)$-invariant coupling of these $\mu_i$.

Having set up this slightly more general class of spaces, we shall deduce the Structure Theorem 2.9 from a generalization to the non-uniform case, relying on suitable extensions of Definitions 2.7 and 2.8.

**Definition 3.8.** By a sequence of **ingredients** we shall now understand a sequence $Z_0, Z_1, \ldots, Z_{k-1}$ of standard Borel spaces, a probability measure $\mu_0$ on $Z_0$, a Borel map $\kappa_i : Z_i \to K_i$ for each $i = 0, 1, \ldots, k-1$ and (setting $Z_k := K$) a probability kernel

$$P_i : Z_0 \times Z_1^i \times Z_2^{\binom{i}{2}} \times \cdots \times Z_{i-1}^{\binom{i}{i-1}} \rightsquigarrow Z_i$$

for each $i = 1, 2, \ldots, k$ that is symmetric under the natural coordinate-permuting action of $\mathrm{Sym}([i])$ on the domain.

**Definition 3.9.** Given a sequence of ingredients

$$(Z_0, \kappa_0, \mu_0), (Z_1, \kappa_1, P_1), \ldots, (Z_{k-1}, \kappa_{k-1}, P_{k-1}), (K, P_k),$$

by the **standard recipe** we understand the following procedure for picking a random member of $K^{(S)}$:

- Choose $z_\emptyset \in Z_0$ at random according to the law $\mu_0$, and set $\omega_\emptyset := \kappa_0(z_\emptyset)$;
- Colour each vertex $s \in S$ by some $z_s \in Z_1$ chosen independently according to $P_1(z_\emptyset, \,\cdot\,)$, and set $\omega_s := \kappa_1(z_s)$
- Colour each edge $a = \{s, t\} \in \binom{S}{2}$ by some $z_a \in Z_2$ chosen independently according to $P_2(z_\emptyset, z_s, z_t, \,\cdot\,)$, and set $\omega_a := \kappa_2(z_a)$;

$$\vdots$$

- Colour each $(k-1)$-hyperedge $u \in \binom{S}{k-1}$ by some $z_u \in Z_{k-1}$ chosen independently according to $P_{k-1}(z_\emptyset, (z_s)_{s \in u}, \ldots, (z_v)_{v \in \binom{u}{k-2}}, \cdot)$, and set $\omega_u := \kappa_{k-1}(z_u)$;
- Colour each $k$-hyperedge $e \in \binom{S}{k}$ by some colour $\omega_e \in K$ chosen independently according to $P_k(z_\emptyset, (z_s)_{s \in e}, \ldots, (z_u)_{u \in \binom{e}{k-1}}, \cdot)$.

If $\mu$ is the law of the resulting random $K$-coloured $k$-uniform hypergraph, we will say that the ingredients **yield** $\mu$ upon following the standard recipe.

Thus, in this non-uniform version of the standard recipe, we retain some of the information contained in the choice of $z_u$ for each $u \in \binom{S}{\leq k}$ (specifically, its image $\kappa_{|u|}(z_u)$), whereas in the $k$-uniform case we threw all of this away at the last step save for the top-rank data $(\omega_e)_{e \in \binom{S}{k}}$.

**Theorem 3.10** (Structure theorem for exchangeable random hypergraph colourings). *For any exchangeable random $K$-coloured hypergraph $\mu$ there is some sequence of ingredients $(Z_0, \kappa_0, \mu_0)$, $(Z_1, \kappa_1, P_1), \ldots, (Z_{k-1}, \kappa_{k-1}, P_{k-1})$, $(K, P_k)$ which yields $\mu$ upon following the standard recipe.*

Let us begin our analysis. Suppose that $K = (K_i)_{i=0}^k$ is a $k$-palette, and again partition $S$ into two infinite subsets, say $S_1$ and $S_2$; as in the cases $k = 1$ and $k = 2$ it suffices to prove the structure theorem for the projection $(\pi_{\binom{S_1}{\leq k}})_\# \mu$. Given this partition, we write

$$K_k^{\binom{S}{k}} = \prod_{i=0}^k K_k^{\{e \in \binom{S}{k}: |e \cap S_1| = i\}} = \prod_{i=0}^k K_k^{\binom{S_1}{i} \times \binom{S_2}{k-i}} = \prod_{i=0}^k \left(K_k^{\binom{S_2}{k-i}}\right)^{\binom{S_1}{i}},$$

where we identify $e \in \binom{S}{k}$ having $|e \cap S_1| = i$ with $(e \cap S_1, e \cap S_2) \in \binom{S_1}{i} \times \binom{S_2}{k-i}$ in the obvious way. Note that we are currently going to use this factorization only for the top rank $k$.

Now let $Y_i := K_k^{\binom{S_2}{k-i}}$ for each $i \leq k-1$, and write $\pi_a^{Y_i}$ for the projection from $Y_i^{\binom{S_1}{i}}$ onto the coordinate indexed by $a \in \binom{S_1}{i}$.

Considering the maps $\pi_a$ for $a \in \binom{S_1}{\leq k}$ and also $\pi_a^{Y_i}$ for $a \in \binom{S_1}{\leq k-1}$, we can obtain from a $K$-coloured hypergraph on $S$ a hypergraph on $S_1$ coloured by the palette $\tilde{K} := (K_0 \times Y_0, K_1 \times Y_1, \ldots, K_{k-1} \times Y_{k-1}, K_k)$. We will show that for the $\tilde{K}$-coloured random hypergraph on $S_1$ obtained from $\mu$ by pushforward, say $\tilde{\mu}$, the $K_k$ colours $\pi_e$ of the $k$-hyperedges $e \in \binom{S_1}{k}$ are relatively independent conditioned on the $(K_i \times Y_i)$ colours of all the sub-edges $a \subset e$. Thus, we will have obtained a random rank-$(k-1)$ coloured hypergraph $\tilde{\mu}$ with enlarged palette $\tilde{K}$ such that the original law $\mu$ can be obtained as the image

$$\mu = \left(\kappa_0, \kappa_1, \ldots, \kappa_{k-1}, \bigotimes_{e \in \binom{S}{k}} P_k \circ \pi_{0, e, \binom{e}{2}, \ldots, \binom{e}{k-1}}\right)_\# \tilde{\mu},$$

where $P_k$ may be obtained explicitly as the conditional distribution of any one of the projections $\pi_e$ for $e \in \binom{S_1}{k}$ over $\pi_{\binom{e}{\leq k-1}}^{\tilde{K}}$ and $\kappa_i : K_i \times Y_i \to K_i$ is the projection onto the first coordinate. This $P_k$ will become the top-level kernel appearing in our list of ingredients, and $K_{k-1} \times Y_{k-1}$ will play the rôle of $Z_{k-1}$. We will then be able to complete the proof by induction on $k$, since we may make an inductive appeal to the structure theorem for the rank-$(k-1)$ random hypergraph colouring $(\pi^{\tilde{K}}_{\binom{S_1}{\leq k-1}})_\#\tilde{\mu}$ to obtain the remaining ingredients.

For convenience, let us write T($u$) for the $\sigma$-subalgebra of $\Sigma_K^{\otimes\binom{S}{k}}$ generated by the maps $(\pi_a, \pi_a^{Y_{|a|}})_{a\subset u}$, and set $\mathrm{T} := \bigvee_{u \in \binom{S_1}{k-1}} \mathrm{T}(u)$. Just as in the graph case, we shall also want to use the $\sigma$-subalgebra $\mathrm{T}^{S_3}$ obtained analogously to the above but using only the smaller infinite subset $S_3 \subset S_2$ as the set of auxiliary vertices.

An energy-increment argument exactly analogous to that for Lemma 3.3 gives the following.

**Lemma 3.11** (Tail property of conditional expectation). *Whenever $S_3 \subseteq S_2$ is a further infinite subset, $e \in \binom{S_1}{k}$ and $f \in L^\infty(K)$, the conditional expectation $\mathsf{E}_\mu[f \circ \pi_e \mid \mathrm{T}]$ is actually measurable with respect to $\mathrm{T}^{S_3}$.* □

**Proposition 3.12** (Relative independence). *The random variables $\pi_e$ for $e \in \binom{S_1}{k}$ are conditionally independent over $\mathrm{T}$: for any finite $F \subset \binom{S_1}{k}$ and bounded measurable functions $f_e : K \to \mathbb{R}$ for $e \in F$ we have*

$$\mathsf{E}_\mu\left[\prod_{e \in F} f_e \circ \pi_e \,\bigg|\, \mathrm{T}\right] = \prod_{e \in F} \mathsf{E}_\mu[f_e \circ \pi_e \mid \mathrm{T}].$$

*Proof.* Enumerate $F$ as $\{e_1, e_2, \ldots, e_r\}$. Arguing just as for Proposition 3.4, it suffices to prove that

$$\int_A \prod_{e \in F} f_e \circ \pi_e \,\mathrm{d}\mu = \int_A \mathsf{E}_\mu[f_{e_1} \circ \pi_{e_1} \mid \mathrm{T}] \cdot \prod_{p=2}^r f_{e_p} \circ \pi_{e_p} \,\mathrm{d}\mu$$

for any finite-dimensional cylinder $A = A_1 \cap A_2 \cap \cdots \cap A_m$ for $A_j \in \sigma((\pi_u)_{u \in \binom{a_j \cup R_j}{\leq k}})$ for some $a_j \in \binom{S_1}{k-1}$ and finite subsets $R_j \subset S_2$, $1 \leq j \leq m$.

Let $e_1 = \{s_1, s_2, \ldots, s_k\}$ and enumerate all the remaining vertices in $(a_1 \cup a_2 \cup \cdots \cup a_m \cup e_1 \cup e_2 \cup \cdots \cup e_r) \setminus e_i$ as $z_1, z_2, \ldots, z_q$. Select distinct vertices $w_1, w_2, \ldots, w_q \in S_2 \setminus (R_1 \cup R_2 \cup \cdots \cup R_m)$ and let $g \in \mathrm{Sym}_0(S)$ be the permutation that swaps $z_i$ and $w_i$ for all $i \leq q$.

By the exchangeability we know that

$$\int_A \prod_{e \in F} f_e \circ \pi_e \,\mathrm{d}\mu = \int_{\tau^g(A)} \prod_{e \in F} f_{g(e)} \circ \pi_{g(e)} \,\mathrm{d}\mu = \int_{\tau^g(A)} f_{e_1} \circ \pi_{e_1} \cdot \prod_{p=2}^r f_{g(e_p)} \circ \pi_{g(e_p)} \,\mathrm{d}\mu,$$

since $g$ does not move any vertex of $e_1$. Now both $\tau^g(A)$ and $\prod_{p=2}^r f_{g(e_p)} \circ \pi_{g(e_p)}$ are T-measurable, since through the action of $g$ we have moved all relevant vertices from $S_1$ into $S_2$, except for those in $e_1$; and $|e_1 \cap e_p| \leq k-1$ for any $p \geq 2$. In particular, every edge $g(e_p)$ for $p \geq 2$ now has at least one end-point in $S_2$. Therefore the left-hand-side above must equal

$$\int_{\tau^g(A)} \mathsf{E}_\mu[f_{e_1} \circ \pi_{e_1} \mid \mathrm{T}] \cdot \prod_{p=2}^r f_{g(e_p)} \circ \pi_{g(e_p)} \, \mathrm{d}\mu.$$

Finally, by Lemma 3.11 the conditional expectation $\mathsf{E}_\mu[f_{e_1} \circ \pi_{e_1} \mid \mathrm{T}]$ is actually measurable with respect to $\mathrm{T}^{S_3}$ for $S_3 := S_2 \setminus (R_1 \cup R_2 \cup \cdots \cup R_m \cup \{w_1, w_2, \ldots, w_q\})$, and is therefore $\tau^g$-invariant; thus, re-applying $g$ we may return the vertices $w_j$ to their original locations $z_j$, and so deduce that

$$\int_A \prod_{e \in F} f_e \circ \pi_e \, \mathrm{d}\mu = \int_A \mathsf{E}_\mu[f_{e_1} \circ \pi_{e_1} \mid \mathrm{T}] \cdot \prod_{p=2}^r f_{e_p} \circ \pi_{e_p} \, \mathrm{d}\mu.$$

This completes the proof. □

Letting

$$P_k : (K_0 \times Y_0) \times (K_1 \times Y_1)^k \times \cdots \times (K_{k-1} \times Y_{k-1})^{\binom{k}{k-1}} \rightsquigarrow K_k$$

be the distribution of any arbitrarily chosen $\pi_e$ for $e \in \binom{S_1}{k}$ (by exchangeability the choice of $e$ will not matter), and writing $Z'_i$ for $K_i \times Y_i$ when $i \leq k-1$ and for $K_k$ alone when $i = k$, the above relative independence immediately gives the following.

**Proposition 3.13.** *If $K$ is a $k$-palette and $\mu$ is an exchangeable random $K$-coloured hypergraph then we can introduce new standard Borel spaces $Z'_0$, $Z'_1$, ..., $Z'_k$ forming a $k$-palette $Z'$, Borel maps $\kappa_i : Z'_i \to K_i$ with $\kappa_k = \mathrm{id}_{K_k}$, an exchangeable random $Z'$-coloured rank-$(k-1)$ hypergraph $\mu_{k-1}$ coloured by them, and a $\mathrm{Sym}([k])$-symmetric probability kernel*

$$P_k : Z'_0 \times (Z'_1)^k \times \cdots \times (Z'_{k-1})^{\binom{k}{k-1}} \rightsquigarrow K_k$$

*such that*

$$\mu = \left(\kappa_0, \kappa_1^{\otimes S}, \kappa_2^{\otimes \binom{S}{2}}, \ldots, \kappa_{k-1}^{\otimes \binom{S}{k-1}}, \bigotimes_{e \in \binom{S}{k}} P_k \circ \pi_{0,e,\binom{e}{2},\ldots,\binom{e}{k-1}}\right)_{\#} \mu_{k-1}.$$

□

We can now prove Theorem 2.9.

*Proof of Theorem 2.9.* Starting from $\mu$, we may apply Theorem 3.13 to obtain spaces $Z'_i$ for $0 \leq i \leq k-1$ enlarging the colour spaces $K_i$, together with a probability kernel $P_k$ as in that theorem. This gives the desired structure at the 'top level' of $\mu$. However, having obtained the new exchangeable random $Z'$-coloured hypergraph $\mu'$ with rank at most $k-1$ from which $\mu$ can thus be extracted, we may iterate the conclusion of Theorem 3.13 to obtain greater and greater enlargements $Z'_i \leftarrow Z_i^{(2)} \leftarrow \cdots \leftarrow Z_i^{(k-i)}$, and at the $i^{\text{th}}$ step also an additional probability kernel

$$P_{k-i+1}^{(i)} : Z_0^{(i)} \times (Z_1^{(i)})^{k-i+1} \times \cdots \times (Z_{k-i}^{(i)})^{\binom{k-i+1}{k-i}} \rightsquigarrow Z_{k-i+1}^{(i-1)},$$

so that if once this iteration is complete we take $Z_i := Z_i^{(k-i)}$ for each $i \leq k-1$ and $P_{k-i+1}$ to be $P_{k-i+1}^{(i)}$ composed with the factors $Z_j \to Z_j^{(i)}$ for $j \leq k-i$, then these ingredients recover the full structure that we want. □

As in Subsection 3.2, let us note that the form of our Structure Theorem 2.9 is somewhat different from that following directly from the probabilistic results of Kallenberg [44]. We shall address this issue in more detail in Subsection 3.7, when we discuss the setting of colourings of partite graphs, since it is actually these that correspond most closely to the original probabilistic results; here let us simply state the extension to arbitrary-$k$ case of the older formalism:

**Theorem 3.14.** *For any $k$-uniform exchangeable random $K$-coloured hypergraph $\mu$ there is some measurable function $f : [0,1]^{\binom{[k]}{k}} \to K$ such that, if $(\xi_a)_{a \in \binom{S}{\leq k}}$ is a collection of independent uniform $[0,1]$-valued random variables then the stochastic process $\big(f(\xi_a)_{a \subseteq e}\big)_{e \in \binom{S}{k}}$ has joint law $\mu$.*

### 3.4. Finer topological consequences of the structure theorem

In later subsections we will describe modifications of Theorem 2.9 to a number of related notions of exchangeability. However, before doing so we wish to elaborate some more delicate topological features of the proof of the last section, in the case when the original palette $K$ consists of finite sets.

Our interest is in a sense in which our extraction of suitable structural ingredients for a given $\mu$ was 'continuous' in that original measure $\mu$. This aspect to these theorems, while an immediate consequence of the proofs, seems to have been ignored in earlier works, since its consequences for probability theory are slight. However, they will be crucial for the use to which the structure theorem will be put in [10], and so we recount them here, and will restate them when we turn to directed hypergraph colourings in the next subsection (again, for the sake of applications in [10]).

Let us re-describe the proof of the previous section in slightly different terms. In order to deduce the structure of the exchangeable law $\mu$, we first identify $S$ with an infinite-co-infinite subset $S_1$ of itself through some bijection $\psi : S \to S_1$, and then agree that by exchangeability it suffices to understand the structure of the projection $(\pi_{\binom{S_1}{\leq k}})_\# \mu$. We then make use of certain $\sigma$-algebras generated by the colours $\pi_e$ of the spare edges in $\binom{S}{k} \setminus \binom{S_1}{k}$ to describe that structure, enough to justify the recursion clause of an inductive argument.

If we unravel the induction, we may write instead that we have chosen a partition $S_{2,1} \cup S_{2,2} \cup \cdots S_{2,k}$ of $S \setminus S_1$ into infinite subsets, and now the estimate of Proposition 3.12 shows that

- defining $\pi_u^{Z_{k-1}}$ to be the cluster of colours

$$(\pi_{e'})_{e' \subseteq S_1 \cup S_{2,1} \cup \cdots \cup S_{2,k},\ e' \cap (S_1 \cup S_{2,1} \cup \cdots \cup S_{2,k-1}) \subseteq u}$$

  for each $u \in \binom{S_1 \cup S_{2,1} \cup \cdots \cup S_{2,k-1}}{k-1}$ (taking values in a space $Z_{k-1}$ given by the obvious product of copies of $K_i$'s), the highest-rank colours $\pi_e$ for $e \in \binom{S_1}{k}$ are relatively independent over $(\pi_u^{Z_{k-1}})_{u \in \binom{S_1}{k-1}}$ (we defined $\pi_u^{Z_{k-1}}$ also for $u$ in the larger collection $\binom{S_1 \cup S_{2,1} \cup \cdots \cup S_{2,k-1}}{k-1}$ for auxiliary use in subsequent steps);

- next, forgetting about the $K_k$-colours of any $k$-subsets not contained in $S_1$, and defining $\pi_v^{Z_{k-2}}$ to be the cluster of colours

$$(\pi_{u'}^{Z_{k-1}})_{u' \subseteq S_1 \cup S_{2,1} \cup \cdots \cup S_{2,k-1},\ u' \cap (S_1 \cup S_{2,1} \cup \cdots \cup S_{2,k-2}) \subseteq v}$$

  for each $v \in \binom{S_1 \cup S_{2,1} \cup \cdots \cup S_{2,k-2}}{k-2}$ (taking values in new larger product space $Z_{k-2}$ similarly to above), now the colours $\pi_u^{Z_{k-1}}$ for $u \in \binom{S_1}{k-1}$ are relatively independent over $(\pi_v^{Z_{k-2}})_{v \in \binom{S_1}{k-2}}$;

$$\vdots$$

Simply by keeping track of which clusters of colours of subsets of $S_1 \cup S_{2,1} \cup \cdots \cup S_{2,k}$ end up inside each of the spaces $Z_i$ during the above construction, and which are simply thrown away, we find that we have proved a 'continuous and functorial' version of the structure theorem in case $K$ is finite (in fact, in case $K$ is itself compact, metrizable and totally disconnected). This is more easily stated with help from a little additional notation.

**Definition 3.15** (Hypergraph powers of maps). Suppose that $Y = (Y_i)_{i=0}^k$ and $Z = (Z_i)_{i=0}^k$ are $k$-palettes, that $\Lambda$ is a sequence of Borel maps $\Lambda_i : Y_i \to Z_i$ and that $S$ is some vertex set. Then by the **hypergraph power of $\Lambda$ over $S$** we understand the map $\Lambda^{(S)} : Y^{(S)} \to Z^{(S)}$ defined by

$$\Lambda^{(S)}\big(y_\emptyset, (y_i)_{i \in S}, \ldots, (y_e)_{e \in \binom{S}{k}}\big) := \big(\Lambda_0(y_\emptyset), (\Lambda_1(y_i))_{i \in S}, \ldots, (\Lambda_k(y_e))_{e \in \binom{S}{k}}\big).$$

It is clear that $\Lambda^{(S)}$ is continuous if each $Y_i$ and $Z_i$ has a topology and $\Lambda_i$ is continuous. We will use these hypergraph powers through the following construction.

**Definition 3.16** (Internalizing vertices). If $K$ is a $k$-palette and $\mu$ is a random $K$-coloured hypergraph on $S_1 \cup S_2$ for two disjoint infinite sets $S_1$ and $S_2$, then we can define the $S_2$-**internalization of $K$** to be the $k$-palette $K^{\uplus S_2} := (K_i^{\uplus S_2})_{i=0}^k$ given by

$$K_i^{\uplus S_2} := K_k^{\binom{S_2}{k-i}} \times K_{k-1}^{\binom{S_2}{k-i-1}} \times \cdots \times K_i,$$

and the $S_2$-**internalization of** $\mu$ to be the $K^{\uplus S_2}$-coloured hypergraph $\mu^{\uplus S_2}$ on $S_1$ obtained by the canonical identification of $K^{(S_1 \cup S_2)}$ with

$$
\begin{aligned}
(K^{\uplus S_2})^{(S_1)} &= \bigl(K_k^{\binom{S_2}{k}} \times K_{k-1}^{\binom{S_2}{k-1}} \times \cdots \times K_1^{S_2} \times K_0\bigr) \\
&\quad \times \bigl(K_k^{\binom{S_2}{k-1}} \times K_{k-1}^{\binom{S_2}{k-2}} \times \cdots \times K_2^{S_2} \times K_1\bigr)^{S_1} \times \cdots \times K_k^{\binom{S_1}{k}}
\end{aligned}
$$

that results from partitioning

$$
\binom{S_1 \cup S_2}{\leq k} = \bigcup_{i \leq k} \bigcup_{j=0}^{k-i} \left\{ u \in \binom{S_1 \cup S_2}{i+j} : |u \cap S_1| = i \right\}
$$

and identifying of $\{u \in \binom{S_1 \cup S_2}{i+j} : |u \cap S_1| = i\}$ with $\binom{S_1}{i} \times \binom{S_2}{j}$.

Now the re-write of Theorem 2.9 to which we have been building is the following.

**Theorem 3.17.** *Suppose that $K$ is a finite k-palette, $\mu$ is an exchangeable random $K$-coloured hypergraph on $S_1$, $S_2$ is some additional countably infinite pool of vertices and $\mu'$ is the unique exchangeable extension of $\mu$ to a random $K$-coloured hypergraph on $S_1 \cup S_2$. Then there are a k-palette $Z = (Z_i)_{i=0}^k$ comprising totally disconnected compact metric spaces (which are therefore homeomorphic to closed subsets of the Cantor space) and collections of continuous maps $\kappa_i : Z_i \to K_i$ and $\Lambda_i : K_i^{\uplus S_2} \to Z_i$ such that under $\Lambda_\#^{(S_1)}(\mu')^{\uplus S_2}$ the random variables $\pi_e$ for $|e| = i$ are relatively independent when conditioned on all the random variables $\pi_u$ with $|u| < i$, and $\mu = \kappa_\#^{(S_1)} \Lambda_\#^{(S_1)}(\mu')^{\uplus S_2}$.* □

It is this version of our result, together with its directed-hypergraph counterpart, that will underpin the later sections of [10]. (It will be formulated there using slightly more category-theoretic terms, as these will help with the organization of that paper, but the content is easily seen to be identical.)

### 3.5. Exchangeable random directed hypergraph colourings

In the previous section we considered probability measures on the space $\prod_{i \leq k} K_i^{\binom{S}{i}}$ indexed by the collection of subsets of $S$ of size at most $k$; a natural extension of this setting is to work instead in the space $\prod_{i \leq k} K_i^{\mathrm{Inj}([i], S)}$ indexed by the sets $\mathrm{Inj}([i], S)$ of all ordered $i$-tuples with distinct terms in $S$ for $i \leq k$: that is, of $K$-coloured *directed* hypergraphs on $S$. Let us first re-assign our notation from Subsection 3.3 in this new setting, now writing $K^{(S)}$ for this space $\prod_{i \leq k} K_i^{\mathrm{Inj}([i], S)}$. Then we can make the obvious analog of Definition 3.7.

**Definition 3.18** (Exchangeable random coloured directed hypergraphs)**.** Given the same data as above, by an **random $K$-coloured directed hypergraph on $S$** we understand a Borel probability measure $\mu$ on $K^{(S)}$ that is invariant under the coordinate-permuting action $\tau$ of $\mathrm{Sym}_0(S)$ defined by

$$
\tau^g \bigl((\omega_\phi)_{\phi \in \bigcup_{i \leq k} \mathrm{Inj}([i], S)}\bigr) = (\omega_{g \circ \phi})_{\phi \in \bigcup_{i \leq k} \mathrm{Inj}([i], S)}.
$$

In order to handle this new setup efficiently, we turn to the partitions $\mathrm{Inj}([i], S) := \bigcup_{e \in \binom{S}{i}} \mathrm{Inj}([i], e)$. For each $i \leq k$ we are considering a colouring from some fixed standard Borel spaces $K_i$ of directed hyperedges $\phi \in \mathrm{Inj}([i], S)$. This can be viewed as assigning to every *undirected* $i$-hyperedge $e$ a colour from the larger space $K_i^{\mathrm{Inj}([i],e)}$, with one factor of $K_i$ appearing for each of the $i!$ possible enumerations of $e$, although there is no canonical way to choose a bijection between the enumerations of one $i$-hyperedge $e$ and of another $i$-hyperedge $e'$.

We will find that our structure theorem for random directed hypergraph colourings requires us to work with these clusters of colours indexed by the corresponding undirected edges. This will lead to an important difference from Theorem 2.9. In the undirected setting, at the $i^{\mathrm{th}}$ step of our recipe we choose the $i$-hyperedge colours independently at random from the distributions given by the kernels $P_i$ depending on the previously-chosen lower-rank hyperedge colours. For directed hypergraphs, however, the choices of colours of directed $i$-hyperedges $\phi$ and $\psi$ with *different images* are still independent, but for $\phi$ and $\psi$ two different injections from $[i]$ into the same undirected $i$-hyperedge $a \in \binom{S}{i}$ they are usually not. Rather, their dependence is coded into the probability kernel $P_i$, and is quite arbitrary subject only to the $\mathrm{Sym}([i])$-covariance condition on $P_i$.

We cannot hope for these last to be independent because the distinct random variables indexed by different enumerations of a single $i$-set $a$ are always moved around together under vertex-permutations, and so, aside from moving the vertices around inside $a$, no averaging over different images and subsequent estimates on covariances can be used to pry apart their joint distribution.

Nevertheless, this slight additional complexity having been taken into account, the basic methods of the preceding sections give an analogous result to Theorem 2.9. We shall give the details of this statement here, but leave the (essentially unchanged) proof to the reader.

**Definition 3.19** (Directed ingredients)**.** By a sequence of **ingredients** we understand a sequence $Z_0, Z_1, \ldots, Z_{k-1}$ of standard Borel spaces, a probability measure $\mu_0$ on $Z_0$ and (setting $Z_k := K_k$) Borel maps $\kappa_i : Z_i \to K_i$ for each $i = 0, 1, \ldots, k-1$ and probability kernels

$$P_i : Z_0 \times Z_1^{\mathrm{Inj}([1],[i])} \times Z_2^{\mathrm{Inj}([2],[i])} \times \cdots \times Z_{i-1}^{\mathrm{Inj}([i-1],[i])} \leadsto Z_i^{\mathrm{Sym}([i])}$$

for each $i = 1, 2, \ldots, k$ that are covariant under the natural coordinate-permuting actions of $\mathrm{Sym}([i])$ on the domain and on the target.

**Definition 3.20** (Directed recipe)**.** Given a sequence of ingredients

$$(Z_0, \kappa_0, \mu_0), (Z_1, \kappa_1, P_1), \ldots, (Z_{k-1}, \kappa_{k-1}, P_{k-1}), (K_k, P_k),$$

by the **standard recipe** we understand the following procedure for picking a random member of $(\omega_\phi)_{\phi \in \bigcup_{i \leq k} \mathrm{Inj}([i],S)} \in K^{(S)}$:

- Choose $z_\emptyset \in Z_0$ according to $\mu_0$, and set $\omega_\emptyset := \kappa_0(z_\emptyset)$;
- Colour each vertex $s \in S$ by $z_s \in Z_1$ independently according to $P_1(z_\emptyset, \cdot)$, and set $\omega_s := \kappa_1(z_s)$;

- Colour each $a = \{s,t\} \in \binom{S}{2}$ by some pair $(z_\phi)_{\phi \in \text{Inj}([2],a)} \in Z_2^{\text{Inj}([2],a)}$ chosen independently according to $P_2(z_\emptyset, z_s, z_t, \cdot)$, and set $\omega_\phi := \kappa_2(z_\phi)$ for each $\phi \in \text{Inj}([2], S)$;

$$\vdots$$

- Colour each $u \in \binom{S}{k-1}$ by some $(k-1)!$-tuple $(z_\phi)_{\phi \in \text{Inj}([k-1],u)}$ in $Z_{k-1}^{\text{Inj}([k-1],u)}$ chosen independently according to $P_{k-1}(z_\emptyset, (z_s)_{s \in u}, \ldots, (z_\psi)_{\psi \in \text{Inj}([k-2],u)}, \cdot)$, and set $\omega_\phi := \kappa_{k-1}(z_\phi)$ for each $\phi \in \text{Inj}([k-1], S)$;
- Colour each $e \in \binom{S}{k}$ by some $k!$-tuple $(\omega_\phi)_{\phi \in \text{Inj}([k],e)} \in K_k^{\text{Inj}([k],e)}$ independently according to $P_k(z_\emptyset, (z_s)_{s \in e}, \ldots, (z_\psi)_{\psi \in \text{Inj}([k-1],e)}, \cdot)$.

If $\mu$ is the law of the resulting random $K$-coloured $k$-uniform hypergraph, we will say that the ingredients **yield** $\mu$ upon following the standard recipe.

**Theorem 3.21** (Structure theorem for exchangeable random directed hypergraph colourings). *For any exchangeable random $K$-coloured directed hypergraph $\mu$ there is some sequence of ingredients $(Z_0, \kappa_0, \mu_0)$, $(Z_1, \kappa_1, P_1)$, ..., $(Z_{k-1}, \kappa_{k-1}, P_{k-1})$, $(K_k, P_k)$ which yields $\mu$ upon following the standard recipe.*

As in the case of undirected graphs, our application of this result in [10] will require the following 'continuous and functorial' version, which follows from the above argument just as did Theorem 3.17 from the argument of Subsection 3.3. We first make an analogous directed version of Definition 3.16; we will use the same notation, since no confusion should arise and this definition collapses to its predecessor if we identify undirected graphs as a special class of directed graphs. We will also use hypergraph powers of maps as in Definition 3.15, but this definition is essentially unchanged.

**Definition 3.22** (Internalizing vertices). If $K$ is a $k$-palette and $\mu$ is a random $K$-coloured directed hypergraph on $S_1 \cup S_2$ for two disjoint infinite sets $S_1$ and $S_2$, then we define the $S_2$**-internalization of** $K$ to be the $k$-palette $K^{\uplus S_2} := (K_i^{\uplus S_2})_{i=0}^k$ given by

$$K_i^{\uplus S_2} := K_k^{\binom{[k]}{[i]} \times \text{Inj}([k-i], S_2)} \times K_{k-1}^{\binom{[k-1]}{[i]} \times \text{Inj}([k-i-1], S_2)} \times \cdots \times K_i,$$

and the $S_2$**-internalization of** $\mu$ to be the $K^{\uplus S_2}$-coloured directed hypergraph $\mu^{\uplus S_2}$ on $S_1$ obtained by the canonical identification of $K^{(S_1 \cup S_2)}$ with

$$(K^{\uplus S_2})^{(S_1)} = \prod_{i=0}^k \left( K_k^{\binom{[k]}{[i]} \times \text{Inj}([k-i], S_2)} \times K_{k-1}^{\binom{[k-1]}{[i]} \times \text{Inj}([k-i-1], S_2)} \times \cdots \times K_i \right)^{\text{Inj}([i], S_1)}$$

that results from partitioning

$$\bigcup_{i \leq k} \text{Inj}([i], S_1 \cup S_2) = \bigcup_{i \leq k} \bigcup_{j=0}^{k-i} \{\phi \in \text{Inj}([i+j], S_1 \cup S_2) : |\phi([i+j]) \cap S_1| = i\}$$

and identifying $\{\phi \in \text{Inj}([i+j], S_1 \cup S_2) : |\phi([i+j]) \cap S_1| = i\}$ with $\text{Inj}([i], S_1) \times \binom{[i+j]}{[i]} \times \text{Inj}([j], S_2)$ in the natural way.

**Theorem 3.23.** *Suppose that $K$ is a finite $k$-palette, $\mu$ is an exchangeable random $K$-coloured directed hypergraph on $S_1$, $S_2$ is some additional countably infinite pool of vertices and $\mu'$ is the unique exchangeable extension of $\mu$ to a random $K$-coloured directed hypergraph on $S_1 \cup S_2$. Then there are a $k$-palette $Z = (Z_i)_{i=0}^k$ comprising totally disconnected compact metric spaces and collections of continuous maps $\kappa_i : Z_i \to K_i$ and $\Lambda_i : K_i^{\uplus S_2} \to Z_i^{\mathrm{Sym}([i])}$ such that under $\Lambda_\#^{(S_1)}(\mu')^{\uplus S_2}$ the clusters of random variables $(\pi_\phi)_{\phi \in \mathrm{Inj}([i],u)}$ for distinct $u \in \binom{S_1}{i}$ are relatively independent when conditioned on all the random variables $\pi_\psi$ with $\psi \in \mathrm{Inj}(<i, S_1)$, and $\mu = \kappa_\#^{(S_1)} \Lambda_\#^{(S_1)} (\mu')^{\uplus S_2}$.* □

Before leaving this subsection, let us mention the following possible common generalization of the settings of undirected and directed exchangeable random coloured hypergraphs. Let us consider again only the $k$-uniform case, as this illustrates all the new ideas. We now posit a whole family of colour spaces $(K_e)_{e \in \binom{S}{k}}$ indexed by the hyperedges, together with an action $g \mapsto \tau^g$ of $\mathrm{Sym}_0$ on the space $\prod_{e \in \binom{S}{k}} K_e$ that involves both the composition with the usual action of $g$ on the base space $\binom{S}{k}$ through permuting $S$, and also some nontrivial 'local' transformation by $g$ of the value of $\pi_e$ for an individual hyperedge $e$: precisely, $\tau : \mathrm{Sym}_0(S) \curvearrowright \prod_{e \in \binom{S}{k}} K_e$ is of the form

$$\tau^g\big((\omega_e)_{e \in \binom{S}{k}}\big) = \big(\sigma_e^g(\omega_{g(e)})\big)_{e \in \binom{S}{k}},$$

where for each $g \in \mathrm{Sym}_0(S)$ and $e \in \binom{S}{k}$ the map $\sigma_e^g$ is a Borel isomorphism $K_{g(e)} \to K_e$ which depends only on the restriction $g|_e$ and is subject to the 'cocycle condition' that $\sigma_e^g \circ \sigma_{g(e)}^{g'} = \sigma_e^{g' \circ g}$ for all $g$, $g'$ (so that we obtain a well-defined group action overall).

Thus the map $\sigma_e^g$ provides a local transformation from the individual fibre $K_{g(e)}$ to $K_e$, but different $g$ with the same image $g(e)$ may give different such transformations. In particular, restricting to those $g$ with $g(e) = e$, we obtain a possibly non-trivial Borel action of $\mathrm{Sym}(e)$ on $K_e$ for each $e$. (We warn the reader that it can*not* in general be reconstructed uniquely from the knowledge of only the action of $\mathrm{Sym}(e)$ on $K_e$ given by those $\sigma_e^g$ with $g(e) = e$.) Given an action $\tau$ of this form involving nontrivial local transformations, we can study those probability measures $\mu$ on $\prod_{e \in \binom{S}{k}} K_e$ that are $\tau$-invariant.

It is straightforward to check that the basic method we have outlined above can be adapted to give a structure theorem for such invariant measures $\mu$, simply upon making all of our ingredients suitably covariant under the local transformations $\sigma_e^g$. From this we can recover Theorem 2.9 by taking $K_e := K$ and each $\sigma_e^g$ to be the identity, and Theorem 3.21 by taking $K_e := K^{\mathrm{Inj}([k],e)}$ and $\sigma_e^g$ to be the natural coordinate-permuting map

$$(\omega_\phi)_{\phi \in \mathrm{Inj}([k],e)} \mapsto (\omega_{\phi \circ g})_{\phi \in \mathrm{Inj}([k],e)}$$

whenever $g(e) = e$. In this second setting we see the clusters of colours indexed by all directed edges with a given underlying undirected emerge again, just as in the earlier arguments of the present subsection.

### 3.6. Two counterexamples

*The need for several type spaces*

We will now give the first of the counterexamples promised in the last part of Subsection 2.2. We show that for simple exchangeable random hypergraph colourings the Structure Theorem 2.9 really does need all of the intermediate type spaces $Z_i$: it cannot be simplified in general to a construction requiring only sampling from some space $Z_1$ at rank 1.

We will give an example of an exchangeable random 3-uniform hypergraph colouring that cannot be constructed as a mixture of vertex-sampling random 3-uniform hypergraph colourings; it is easy to construct similar higher-rank examples. First choose for each $e \in \binom{S}{2}$ independently and uniformly a random colour from the set $\{\text{red}, \text{blue}\}$. Now for each $u \in \binom{S}{3}$ let $\omega_u = 1$ if at least one of the 2-subsets of $u$ is coloured red; otherwise set $\omega_u = 0$. Our example will be the law $\mu$ of the resulting random 3-uniform $\{0,1\}$-coloured hypergraph $(\omega_u)_{u \in \binom{S}{3}}$.

An easy check shows that this $\mu$ is $\mathrm{Sym}_0(S)$-ergodic, so cannot be a nontrivial mixture of other exchangeable random hypergraphs. We will show that it also cannot be obtained by sampling vertex colours at random from some probability space $(V, \nu)$ and then applying some symmetric probability kernel $P : V^3 \rightsquigarrow \{0,1\}$ (that is, we cannot jump over the intermediate step of randomly colouring 2-edges in our recipe). Indeed, suppose that $\mu'$ is any exchangeable random hypergraph that can be constructed in this way, and consider, as in the main argument, a partition $S = S_1 \cup S_2$ of $S$ into two disjoint infinite subsets. For each $s \in S_1$ we introduce the $\sigma$-algebra $\mathrm{T}'(s) := \sigma((\pi_u)_{u \in \binom{\{s\} \cup S_2}{3}})$, and can now compute explicitly from the assumed vertex-sampling-only structure of $\mu'$ that under $\mu'$ the 3-hyperedge random variables $\pi_u$ for $u \in \binom{S}{3}$ must be relatively independent over the $\sigma$-algebra $\mathrm{T}' := \bigvee_{s \in S_1} \mathrm{T}'(s)$. On the other hand, a quick calculation shows that if $u_1, u_2 \in \binom{S}{3}$ share exactly one 2-edge $e$, then under the law $\mu$ the joint behaviour of $\pi_{u_1}$ and $\pi_{u_2}$ is already independent from the whole of $\mathrm{T}'$, but

$$\mu\{\pi_{u_1} = \pi_{u_2} = 1\} = \frac{25}{32} \neq \left(\frac{7}{8}\right)^2 = \mu\{\pi_{u_1} = 1\}\mu\{\pi_{u_2} = 1\}.$$

This contradicts the abovementioned relative independence over $\mathrm{T}'$ if $\mu = \mu'$, and so completes our argument.

It is worth mentioning the similarity between the necessity of introducing a whole list of auxiliary spaces $Z_i$ and a general phenomenon surrounding versions of the regularity lemma for higher-rank hypergraphs: in the case of the latter, in order to yield effective counting lemmas they all require the introduction of a whole hierarchy of vertex-partitions, each dedicated to the control of a particular size of subsets of the vertex set $V$ up to size $k$. We will not examine this similarity in more detail here, but refer the reader to Gowers' discussion in [38, 36]. ◁

*The need for unrecoverable quasifactors*

Intuitively, one might wish to do a little more than the proof we have given of Theorem 2.9: to wit, obtain the spaces $Z_i^{\otimes \binom{S}{i}}$ for our list of ingredients as factors of the original hypergraph colouring probability space $(K^{(S)}, \mu)$ (we return to the convention of Subsection 3.3 for the notation '$K^{(S)}$'), without the trick of passing to an induced random coloured subhypergraph on $K^{(S_1)}$. In case $k = 2$ this would amount to showing that we can find an $S$-indexed family of measurable maps $\phi_s : K^{(S)} \to Z_1$ having the following two properties:

- covariance: $\phi_s \circ \tau^g = \phi_{g(s)}$ for every $g \in \mathrm{Sym}_0$;
- the edge random variables $\pi_e$ are relatively independent over the combined map $(\phi_s)_{s \in S}$.

We will now show that this stronger conclusion can fail (and similarly in case $k \geq 2$), and so some trick like the passage to induced subgraphs followed by structural pull-back is really necessary. We will give an example with $k = 2$ in which, intuitively, while no one vertex of $S_2$ is important for the construction of the auxiliary spaces $Z_1$, we cannot do away with *all* vertices of $S_2$ without losing relevant information.

We first derive some consequences for a family $(\phi_s)_{s \in S}$ of measurable maps of the first of the above conditions in case $\mu$ is an exchangeable random $K$-coloured graph (for a single colour-space $K$, not a more elaborate 2-palette); we will then show for our example $\mu$ that we cannot find a family that also satisfies the second condition.

For such $\mu$, the $\sigma$-subalgebra $\phi_s^{-1}(\Sigma_{Z_1}) \subseteq \Sigma_{K_2}^{\otimes \binom{S}{2}}$ must be contained in the $\mu$-completion of the $\sigma$-subalgebra $\Xi_s \subseteq \Sigma$ of all those $A \in \Sigma$ that are invariant under $\mathrm{Stab}_{\mathrm{Sym}_0}(s) := \{\tau^g : g(s) = s\}$: indeed, by the covariance of these maps $\phi_s$, if $B \in \Sigma_{Z_1}$ and $g(s) = s$ then

$$(\tau^g)^{-1}(\phi_s^{-1}(B)) = \phi_{g(s)}^{-1}(B) = \phi_s^{-1}(B).$$

We will use a certain approximability result for this $\sigma$-algebra $\Xi_s$. Since $\mu$ is a Radon measure, for any $A \in \Xi_s$ and any $\varepsilon > 0$ we may find some finite $I \subset S \setminus \{s\}$ and some finite-dimensional subset $B \subseteq \{0,1\}^{\binom{\{s\} \cup I}{2}}$ such that $\mu(A \triangle B) < \varepsilon$, and so $\|1_A - 1_B\|_{L^1(\mu)} < \varepsilon$. However, since $A$ is fixed by $\mathrm{Stab}_{\mathrm{Sym}_0}(s)$, this gives also $\|1_A - 1_B \circ \tau^g\|_{L^1(\mu)}$ for every $g \in \mathrm{Stab}_{\mathrm{Sym}_0}(s)$, and hence

$$\left\| 1_A - \frac{1}{|F|} \sum_{g \in F} 1_{\tau^g(B)} \right\|_{L^1(\mu)} < \varepsilon$$

for all finite $F \subset \mathrm{Stab}_{\mathrm{Sym}_0}(s)$. Since this group (being isomorphic to $\mathrm{Sym}_0$, for example) is amenable, by the ergodic theorem for amenable groups we may take a pointwise and $L^1$-limit of averages over the members of some Følner sequence for $\mathrm{Stab}_{\mathrm{Sym}_0}(s)$, and so obtain from $1_B$ some function $f$ that is $\Xi_s$-measurable and such that $\|1_A - f\|_{L^1(\mu)} < \varepsilon$. Hence, since $\varepsilon$ was arbitrary, we have proved that any $A \in \Xi_s$ lies in the $\mu$-completion of the $\sigma$-subalgebra $\mathrm{T}_\infty(s) \subseteq \Sigma$ generated by functions, such as $f$, obtained as $\mathrm{Stab}_{\mathrm{Sym}_0}(s)$-ergodic averages of finite-dimensional functions. Moreover, since $1_B$ may be written as a finite sum of indicator functions $1_{\{\omega|_{\binom{\{s\} \cup I}{2}} = \eta\}}$ for various $\eta \in \{0,1\}^{\binom{\{s\} \cup I}{2}}$, we may further reduce to considering $B$ of this particularly simple form.

We are now ready for our example. Consider the exchangeable random ($\{0,1\}$-coloured) hypergraph obtained from ingredients $Z_0 := \{*\}$, $Z_1 := \{\mathrm{red}, \mathrm{blue}\}$, and with conditional edge-colour distributions $P_2(z_1, z_2, \cdot)$ given by

$$P_2(z_1, z_2, \{1\}) = 1 - \delta_{z_1, z_2}:$$

that is, we sample a family $(z_s)_{s \in S}$ of colours, red or blue, independently and fairly at random, and then insert an edge between $s, t \in S$ if and only if $s$ and $t$ are different colours. This is simply a random complete bipartite graph.

We will show that for this $\mu$ the $\sigma$-subalgebras $\Xi_s$ of $\Sigma$ are 'too small' to allow the required conditional independence of the edge random variables $\pi_e$.

Consider a finite-dimensional cylinder set $B \subseteq \{0,1\}^{\binom{\{s\} \cup I}{2}}$ corresponding to some $\eta \in \{0,1\}^{\binom{\{s\} \cup I}{2}}$ as above. Take the Følner sequence for $\mathrm{Stab}_{\mathrm{Sym}_0}(s)$ comprising the family $F_n := \mathrm{Sym}(J_n)$ for some increasing sequence $J_1 \subseteq J_2 \subseteq J_3 \subseteq \ldots$ of finite sets with union $S \setminus \{s\}$. Once $n$ is large enough we must have $I \subseteq J_n$, and then for $\omega \in X$ the function

$$\frac{1}{|\mathrm{Sym}(J_n)|} \sum_{g \in \mathrm{Sym}(J_n)} 1_{\tau^g(B)}(\omega)$$

simply counts the number of isomorphic copies of the graph $\eta$ that occur in $\omega|_{\binom{J_n}{2}}$ if we insist that $s$ does not move but choose an arbitrary size-$|I|$ subset of $J_n$ for the remaining vertices.

For $n$ very large, this statistic cannot distinguish between the two possible vertex colours, red or blue, since there is an automorphism of $(Z_1, P_2)$ that swaps them: if we construct the probability measure on the whole of $\{0,1\}^{\binom{S}{2}} \times Z_1^S$ that arises from our standard recipe, then conditional on either of the events '$s$ is coloured red' and '$s$ is coloured blue', our above ergodic averages converge to the same ($\mu$-almost surely constant) value. Heuristically, for any vertex, an observer situated at that vertex and able to observe its connectivity to all the other vertices and the connectivity among those other vertices can accurately discern the two infinite clusters of the bipartition of the graph, but has no statistical way to tell which is red and which blue, and so cannot tell the colour of any given vertex, including their own. More formally, since the above average converges to a constant for any such $\eta$, we conclude that any ergodic-average function such as $f$ above is $\mu$-almost surely constant, and so the measure-algebra of $(K^{(S)}, \Xi_s, \mu|_{\Xi_s})$ is trivial for any $s \in S$. It follows that for $e = \{s, t\}$ the conditional distribution of $\pi_e$ relative to $\Xi_s \vee \Xi_t$ must be $\mu$-almost surely constant; since, on the other hand, the edge random variables $\pi_e$ for different $e$ are not independent under $\mu$, it follows that our putative factor maps $\phi_s$ cannot have the second of the properties listed above. ◁

### 3.7. Partite hypergraphs and Gowers norms

Let $S_1, S_2, \ldots, S_d$ be a finite list of disjoint countably infinite vertex set and $k \leq d$. We will now briefly treat random colourings of the complete $k$-uniform hypergraph on the vertex classes $S_1, S_2, \ldots, S_d$ (restricting here to the case of a fixed rank $k$ for simplicity). For each $e \in \binom{[d]}{k}$ we refer to the points of the product space $\prod_{i \in e} S_i$ as the **hyperedges above** $e$. Let us introduce the notation $\binom{S_1, S_2, S_3, \ldots, S_d}{k}$ for the set $\bigcup_{e \in \binom{[d]}{k}} \prod_{i \in e} S_i$ of all $k$-hyperedges of the complete such partite hypergraph. Suppose that $K_e$ is some standard Borel space for each $e \in \binom{[d]}{k}$, and consider a $(K_e)_{e \in \binom{[d]}{k}}$-coloured exchangeable $d$-partite $k$-uniform hypergraph $\mu$ on these vertex sets: that is, a probability measure on the product space $\prod_{e \in \binom{[d]}{k}} K_e^{\prod_{i \in e} S_i}$ that is invariant under the separate coordinate-permutation action of each $\mathrm{Sym}_0(S_i)$ on $S_i$.

In this case one can prove another structure theorem along just the same lines as those reviewed above. It tells us that we may always build up an exchangeable random $d$-partite $k$-uniform $K$-colouring $\mu$ according to the following recipe:

**Ingredients**: standard Borel spaces $Z_b$ indexed by all $b \in \binom{[d]}{\leq k}$ with $Z_a = K_a$ when $|a| = k$; a probability measure $\mu_0$ on $Z_\emptyset$; and probability kernels $P_c : \prod_{b \in \binom{c}{|c|-1}} Z_b \rightsquigarrow Z_c$ for all nonempty $c \in \binom{[d]}{k}$.

**Standard recipe**:

- Choose $z_\emptyset \in Z_\emptyset$ at random from $Z_\emptyset$ according to $\mu_0$;
- Choose $z_s \in Z_i$ for each $i \in [d]$ and $s \in S_i$ independently with law $P_i(z_\emptyset, \cdot)$;
- Choose $z_u \in Z_b$ for each $b = \{i,j\} \in \binom{[d]}{2}$ and $u = (s_1, s_2) \in S_i \times S_j$ independently with law $P_b(z_\emptyset, z_{s_1}, z_{s_2}, \cdot)$;

$$\vdots$$

- Choose $\omega_w \in K_e$ for each $e \in \binom{[d]}{k}$ and $w \in \prod_{i \in e} S_i$ independently with law $P_e((z_{w|_b})_{b \in \binom{e}{\leq k-1}}, \cdot)$.

We have assumed no symmetry at all between distinct $S_i$ in the above description, but it is not hard to see that such additional symmetry would be reflected in various further finitary symmetry constraints on the kernels $P_a$. Our main result is now the following.

**Theorem 3.24** (Structure theorem for exchangeable random partite hypergraph colourings). *For any exchangeable random $K$-coloured $k$-uniform hypergraph there is some collection of ingredients $(Z_b)_{b \in \binom{[d]}{\leq k}}$, $\mu_0$ and $(P_c)_{c \in \binom{[d]}{k}}$ which yields $\mu$ upon following the standard recipe.*

The proof that the above recipe can describe all exchangeable partite hypergraph colourings is essentially analogous to that for non-partite hypergraphs given previously, and we shall not spell out the details. In fact, the foundational work by probabilists on exchangeability corresponds most closely to this setting when $d = k$: in particular, exchangeable partite graphs with $d = k = 2$ first appeared in the probabilistic literature under the titles 'separately-exchangeable arrays' or 'row-and-column exchangeable arrays'. See, for example, the discussions in Section 14 of Aldous [3] and the introduction to Kallenberg [44].

These older theorems are not stated in quite the same formalism as that we have adopted in this paper, so let us briefly compare them. We shall take Aldous' treatment of the case $d = k = 2$ in Theorem 14.11 of [3] as representative of the traditional probabilistic approach; the presentation used by Kallenberg in [44] is very similar. Assume that $k = d$ and $S_1 = S_2 = \cdots = S_d$ and write $K$ for $K_e$. Observe that if $(\xi_{b,w})_{b \subseteq [d],\, w \in S^b}$ is any collection of independent uniform $[0,1]$-valued random variables and $f : [0,1]^{\binom{[d]}{\leq d}} \to K$ is any fixed measurable function, then the stochastic process $(\pi_w)_{w \in S^d} := \big(f((\xi_{b,(w_i)_{i \in b}})_{b \subseteq e})\big)_{w \in S^d}$ has law an exchangeable $K$-colouring of $S^b$; in this case we write that $f$ **represents** the law of the process $\pi$. The classical form of the structure theorem is now the following:

**Theorem 3.25.** *Any exchangeable random $K$-coloured $d$-uniform $d$-partite hypergraph $\mu$ is represented by some $f$.*

The difference between this result and Theorem 3.24 above can be seen as one of the choice of representing ingredients; however, it is not too difficult to move between them.

On the one hand, if $\mu$ is representable, then given $f$ we may set $P_e$ to be its conditional distribution on the factor $[0,1]^{\binom{[d]}{\leq d-1}} \leftarrow [0,1]^{\binom{[d]}{d}}$, and now set $Z_b := [0,1]$ when $|b| \leq k - 1$ and $\mu_0$ and $P_c(z, \cdot)$ to be Lebesgue measure for all $z \in [0,1]^{\binom{c}{c-1}}$ whenever $|c| \leq k - 1$; it now follows routinely that these ingredients yield precisely the joint distribution of the process $(\pi_w)_{w \in S^d}$ upon following the standard recipe, and so we have deduced Theorem 3.24 for this instance of a representable process.

In order to prove the reverse implication, it is necessary to process a collection of ingredients $Z_b$, $\mu_0$ and $P_c$ by replacing each space $Z_b$ with a copy of the unit interval $[0,1]$ and then 'outsourcing' the noise in each probability kernel $P_c$ to an independent random variable drawn uniformly from that interval. This is also a very standard argument, but we shall defer presenting this method until Subsection 4.2, where we shall use it to deduce a recent representation theorem of Elek and Szegedy [24] (Theorem 4.1) from our version of the Structure Theorem in the non-partite case, having proved our ability to perform this 'noise outsourcing' in a suitable form in Lemma 4.2.

Before we leave this subsection, it is worth noting that it is $(k + 1)$-partite $k$-uniform hypergraphs, not just ordinary $k$-hypergraphs, that emerge naturally in the recent hypergraph regularity approaches to Szemerédi's Theorem (see, for example, Nagle, Rödl and Schacht [53], Gowers [36] and Chapter 10 of Tao and Vu [64]).

While the parallels between the structure theory for partite hypergraphs and finitary hypergraph regularity lemmas are much as in the non-partite case, in this partite picture it is also quite natural to introduce an infinitary notion of the 'uniformity seminorms' that were first used by Gowers ([38, 36]) on the way to his version of hypergraph regularity.

**Definition 3.26** (Gowers norms). If $\mu$ is an exchangeable law on $\prod_{e \in \binom{[d]}{k}} K_e^{\prod_{i \in e} S_i}$ and $f \in L^\infty(K_e)$ for some $e \in \binom{[d]}{k}$ then we define the **Gowers uniformity seminorm of $f$ above $e$ under $\mu$** to be

$$\|f\|_{\mathrm{U}^e(\mu)} := \left( \int \prod_{\eta \in \{1,2\}^e} f \circ \pi_{w_\eta} \, \mathrm{d}\mu \right)^{2^{-k}},$$

where we have chosen arbitrarily a collection of pairs $s_1^i \ne s_2^i \in S_i$ for $i \in e$ and then written $w_\eta := (s_{\eta_i}^i)_{i \in e}$ for $\eta \in \{1,2\}^e$. It is clear that by exchangeability, any choice of distinct $s_1^i, s_2^i \in S_i$ for each $i \in e$ will give the same value.

By mimicking the analysis in the finitary setting, it is easy to show that $\|\cdot\|_{\mathrm{U}^e(\mu)}$ is a seminorm. Given some complicated exchangeable partite $(K_e)_{e \in \binom{[d]}{k}}$-colouring $\mu$, it is clear that generically all functions $f \in L^\infty(K_e)$ have $\|f\|_{\mathrm{U}^e(\mu)} > 0$; however, after implementing the structure theorem we obtain an extended partite non-uniform exchangeable random hypergraph colouring $\tilde\mu$ of $\bigcup_{j \le k} \binom{S_1, S_2, \ldots, S_k}{j}$ by some auxiliary palette $(Z_a)_{a \subseteq [k]}$ with $Z_e = K_e$ when $|e| = k$, and it can be shown that on this enlarged product space the functions $f \in L^\infty(Z_a)$ for which $\|f\|_{\mathrm{U}^a(\tilde\mu)} = 0$ are precisely those such that $f \circ \pi_w$ is $\tilde\mu$-almost measurable with respect to $(\pi_{w|_b}^{Z_b})_{b \in \binom{a}{\le |a|-1}}$, for any choice of $w \in \prod_{i \in a} S_i$. This closely parallels Lemma 4.3 in Host and Kra's use in [43] of an infinitary analog of the related arithmetic Gowers norms (see [37, 64]), and the proof is exactly similar.

### 3.8. Models of simple theories

Our structural results applied to hypergraphs, directed hypergraphs and towergraphs can be embedded into a somewhat more general setting that has already emerged to a similar purpose in recent work of Razborov [54] (to which we shall return in Subsection 4.3 below). We shall assume various definitions from model theory; see, for example, Chapter 1 of Kopperman [49]. Let $\mathcal{T}$ be a universal first-order theory with equality in a language $\mathcal{L}$ that contains only predicate symbols. Let us suppose first that these symbols have arity at most some finite $k \ge 1$, and (for convenience) assume that $\mathcal{T}$ has only a countable set $\mathcal{S}$ of such predicate symbols; assume further that $\mathcal{T}$ has infinite models. For each $i \le k$ let $\mathcal{S}_i \subseteq \mathcal{S}$ contain those symbols of arity $i$, and let $K_i$ be the space $\{0,1\}^{\mathcal{S}_i}$ with its product topology and Borel $\sigma$-algebra; points of $K_i$ are to be regarded as truth-assignments to the predicates of $\mathcal{S}_i$.

We should stress that we have slipped into the rather abstract lexicon of model theory for its convenience; for theories $\mathcal{T}$ as above, our guiding intuitions will remain those of measures with certain symmetries on a Cantor space.

If the theory $\mathcal{T}$ is free then its models with underlying vertex set $\mathbb{N}$ are precisely the maximal-rank-$k$ directed hypergraph colourings over $\mathbb{N}$ coloured by $K_0, K_1, \ldots, K_k$, except that now we must also allow 'loops': a tuple $(x_1, x_2, \ldots, x_i) \in \mathbb{N}^i$ in which some coordinate appears more than once can also be an argument for an arity-$i$ symbol. Thus in the free case our space of models over the vertex set $\mathbb{N}$ can be identified as $X := \prod_{i \leq k} K^{\mathbb{N}^i}$.

If $\mathcal{T}$ is not free (but still does admit an infinite model), we must correspondingly restrict to the subset $X_\mathcal{T}$ of $X := \prod_{i \leq k} K^{\mathbb{N}^i}$ containing those points that are models of $\mathcal{T}$. This is a closed, hence compact, subset of $X$, since any individual interpretation of a sentence in $\mathcal{T}$ over some particular finite set of vertices in $\mathbb{N}$ simply carves out some clopen subset of $X$ depending only on those vertices as coordinates; the existence of an infinite model is equivalent to the non-emptiness of the intersection of these clopen subsets. The resulting closed subset $X_\mathcal{T}$ is invariant under coordinate-permutation.

Given any such theory $\mathcal{T}$, we can consider the compact convex set $Q_\mathcal{T}$ of Radon probability measures on $X_\mathcal{T}$ invariant under the obvious $\mathrm{Sym}_0(\mathbb{N})$-action, with its vague topology, and ask for a description of the structure of these measures. In fact, such measures have a long history in the model theoretic literature: see, in particular, the papers of Gaifman [34] and Krauss [51], and also the discussion of these actions of $\mathrm{Sym}_0(\mathbb{N})$ as the 'logic actions' (although without the introduction of invariant measures) in Section 2.5 of Becker and Kechris [11]. Of course, we may identify these invariant measures as the exchangeable random hypergraph $K$-colourings that are supported on the closed subset $X_\mathcal{T}$ of $X$, and so we do at least know that they can be described by the standard recipe, but now the additional constraints imposed by $\mathcal{T}$ translate (at least in principle) into additional 'fine-tuning' conditions on the ingredients. Various more precise questions may now be posed about these. For example:

**Question 3.27.** *Given a theory $\mathcal{T}$ having only function symbols of rank* 2, *when is it the case that any (say, ergodic) $\mu$ with support in $X_\mathcal{T}$ can be represented using ingredients $P_2 : Z_1^2 \rightsquigarrow K^{\mathrm{Sym}([2])}$ that can themselves be taken to be deterministic maps, and hence correspond to measurable models of $\mathcal{T}$ with vertex set equal to some fixed copy of the spaces $Z_1$? What happens if the rank is* 3 *or greater?*

It follows from results of Fremlin [28] that equivalence relations (which fit into the above picture with rank 2) do behave in this way, and a positive answer to Fremlin's Problem FY ([27]) would show the same for partial orders. On the other hand, in the free case of graphs and hypergraphs the need for non-deterministic probability kernels $P_j$ in most cases is clear, and so these cannot satisfy the above condition.

A different class of questions pertaining to essentially the same formalism as above are those around the testing of hereditary properties of coloured hypergraph, to which we will return in Subsection 4.5, and in more detail in [10].

### *3.9. A weakened hypothesis: spreadability*

We will finish our review of the classical probabilistic theory by considering another direction in which many of the above forms of exchangeability can be weakened. We now suppose that $T$ is a (necessarily infinite) index set and that $\Gamma$ is a semigroup of self-injections of $T$ (crucially, which may not be invertible), and in this context write that the law $\mu$ of the canonical process $(\pi_t)_{t \in T}$ is **$\Gamma$-spreadable** if $(\pi_{g(t)})_{t \in T}$ still has joint law $\mu$ for any $g \in \Gamma$. As in the exchangeable case, it turns out that if $\Gamma$ is a sufficiently rich class of self-injections then these spreadable laws $\mu$ must still take quite a precise form (and, indeed, very often the resulting structure theorem for a spreadability context subsumes some result for a related exchangeability context on the same index set).

Our leading examples of spreadability, as of exchangeability, correspond to spaces of hypergraph colourings over some countably infinite vertex set $S$, but now with the additional data of a fixed total order $<$ on $S$ and the requirement that our law $\mu$ be invariant under the semigroup $\Gamma$ of *order-preserving* self-injections of $S$. In the special case $k = 1$, the spreadable generalization of de Finetti's Theorem was proved by Ryll-Nardzewski in [58]; the results for higher ranks $k$ were then settled by Kallenberg in [44]. We shall discuss the methods needed for these structural results only cursorily here, referring the reader to this last paper for a complete account. We note that spreadability is referred to as 'spreading-invariance' (in the special case $k = 1$) in Kingman [48] and as 'contractibility' in Kallenberg's more recent book [46].

In the context of such hypergraph spreadable laws, the use of auxiliary vertices often requires the a priori observation that we have some freedom to choose the countable total order $(S, <)$ (up to order-isomorphism), since our law on, say, $K^{\binom{S}{k}}$ will always be determined by its finite-dimensional marginals, and so each spreadable law for any one countably infinite choice $(S, <)$ determines uniquely such a law for any other such choice. Often the proof of the relevant structure theorem proceeds much more smoothly provided we begin with a sufficiently rich total order $(S, <)$, in order to provide enough 'reference vertices' for our subsequent argument.

So we still use the trick of splitting $T$ into a subfamily $T_1$ and an auxiliary collection of 'reference indices', but this now takes typically a much more elaborate form. For example, for Lemma 4.4 of [44], providing a structure theorem for spreadable $\binom{\mathbb{N}}{\leq k}$-indexed processes, Kallenberg must embed $\mathbb{N}$ as a subset of the much more complicated total order on $\mathbb{Q}$, and prove (by a straightforward appeal to the Daniell-Kolmogorov extension theorem) that any spreadable $\binom{\mathbb{N}}{\leq k}$-indexed process can be correspondingly embedded in a spreadable $\binom{\mathbb{Q}}{\leq k}$-indexed process.

The above analysis can also be pushed a little further, by introducing canonical processes $(\pi_\alpha)_{\alpha \in \binom{\mathbb{N}}{<\infty}}$ indexed by *all* finite subsets of $\mathbb{N}$ and with values in some infinite palette $(K_0, K_1, \ldots)$, and ask for the consequences of exchangeability for their structure. This study, too, can be carried out for the weakened notion of spreadability, and an analogous structure theorem is still avail-able: in some sense this is the end-point of the classical examination of exchangeability and spreadability, appearing in Kallenberg [44]. He proves that there always exists an infinite auxiliary palette $(Z_0, Z_1, \ldots)$, probability kernels $P_d : \prod_{\alpha \subset [d]} Z_{|\alpha|} \rightsquigarrow Z_d$ and measurable functions $\kappa_d : Z_d \to K_d$ such the law of the canonical process $(\pi_\alpha)_{\alpha \in \binom{\mathbb{N}}{<\infty}}$ agrees with that of $(\kappa_{|\alpha|}(\xi_\alpha))_{\alpha \in \binom{\mathbb{N}}{<\infty}}$, where $(\xi_\alpha)_{\alpha \in \binom{\mathbb{N}}{<\infty}}$ is the set-indexed process with values in the palette $(Z_0, Z_1, \ldots)$ and law built from the above kernels following the standard recipe. Aside from some finitary symmetry assumptions, the kernels $P_d$ and maps $k_d$ are arbitrary, and this is a more-or-complete explicit structure theorem. It can be proved (essentially) by first establishing a version of the structure theorem for (possibly non-uniform) hypergraphs of finite rank $k$ and then treating a general spreadable set-indexed process as a union of subprocesses corresponding to increasing $k$, so that the final structure emerges from an iterated enlargement of the auxiliary spaces $Z_d$ to handle the sub-processes of $(\pi_\alpha)_{\alpha \in \binom{\mathbb{N}}{<\infty}}$ indexed by larger and large finite subsets.

Let us not leave spreadability without also mentioning the analysis of Fremlin and Talagrand [29], which appeared independently of those mentioned above. Their initial interest is in a class of extremal questions for quite general random graphs on $\mathbb{N}$ that were first raised by Erdős and Hajnal in [25] (and comparable with the questions we shall raise in Subsection 4.6 below). They first reduce to the case of spreadable random graphs (under the name 'deletion-invariant' random graphs), and then resolve them by first proving a version of the structure theorem and then performing a variational analysis on the ingredients that go into it. Their extraction of the structure theorem is rather different from that used, for example, by Kallenberg: having descended to a spreadable random graph they consider the behaviour of all its individual samples along all non-principal ultrafilters on $\mathbb{N}$. These together examine the samples of the graph in much 'finer detail' than our analysis of conditional expectations, and circumvents the need to embed $\mathbb{N}$ into a much larger totally ordered set to provide an infinite pool of auxiliary vertices. This has the advantage that the behaviour of the graph along any one such ultrafilter must be very simple — everything that can converge does converge — but then the Čech-Stone remainder $\beta^*\mathbb{N}$ of $\mathbb{N}$ emerges to play the part of our auxiliary space $Z_1$, equipped with some Radon measure $\nu_1$ defined from the law of the original random graph, and some measurable function $h : \beta^*\mathbb{N} \times \beta^*\mathbb{N} \to [0, 1]$ which (in our terminology) simply corresponds to $P_2(\,\cdot\,, \{1\})$. Their basic construction is outlined in Subsections 2A and 2B of [29], and their version of the structure theorem in Section 5. Let us note that their statement of the structure theorem in terms of this auxiliary function $h$ is arguably closer to the older probabilistic formalism of 'representing' an exchangeable law, as we have discussed for Theorems 3.6, 3.14 and 3.25, except that Fremlin and Talagrand use $\beta^*\mathbb{N}$ in place of $[0, 1]$, and also allow for the construction of different possible probability measures on $\beta^*\mathbb{N}$.

## 4. Relations to finitary combinatorics

We now turn to the parallels between exchangeable random colourings and the 'statistics' of their counterpart finitary structures, building on the general discussion of Subsection 2.3. In the first two subsection we compare exchangeable laws with other notions of 'limit object' for sequences of graphs or hypergraphs, and then give infinitary reformulations of various questions in property testing (Subsection 4.5) and extremal combinatorics (Subsection 4.6). We will finish by comparing the theory outlined above with related ideas in ergodic theory (Subsection 4.7).

### *4.1. The extraction of limit objects*

There are various categories of interest to analysts and geometers in which asymptotic information about a sequence of objects can be parceled into a suitable notion of 'limit object'. This limit object may also lie in the original category, or may require an enlargement through a sort of 'compactification' of that category. If the objects all exist within some larger ambient space with a given topological structure, then this ambient space may itself serve as the larger category; however, if no such superstructure is apparent then a more intrinsic definition of the limit may be necessary. Sometimes both options are available.

Let us first discuss this program in some generality. As illustration we will appeal to two fairly classical examples: the category of metric spaces with Lipshitz maps and that of Banach spaces with continuous linear operators.

Suppose we are working in one of these categories. One option may be to find some large ambient space into which each of a sequence of objects can be mapped, and then consider convergence in terms of some natural topology on the points of that space. Sometimes this choice of ambient space may not be canonical, and sometimes it may not be possible at all; however, once it has been made then we hope that convergence in that ambient space will be closely related to convergence of the parameters or other invariants of the objects that were of interest. Of course, we should not expect the resulting limit object in the ambient space to retain all of the attributes of a typical term of the sequence. It need only distill those properties that are relevant to the convergence, and may suppress the others; thus, having identified the possible limit objects, one should exercise caution before pursuing an overly-elaborate analogy with the original objects. Examples of this construction include Gromov-Hausdorff convergence of metric spaces (see, for example, Chapter 3 of [39]) and Banach-Mazur convergence of Banach spaces of a given dimension (see, for example, Chapter II.E of Wojtaszczyk [65]).

Alternatively we may have at our disposal a rather more forceful, universal construction: that of ultralimits taken along some fixed nonprincipal ultrafilter $\mathcal{U}$. Where ultralimits can be constructed at all, they typically do not depend on the 'convergence' of the sequence of objects in any sense; rather, they merely 'ignore' so much of the sequence that what they retain is forced to converge.

Ultralimits appear in the same categories as did our previous two examples: in metric spaces (Section 3.29 of [39]), where they are used, for example, in the construction of asymptotic cones for various classes of metric space (Section $3.29\frac{1}{2}$); and ultraproducts of Banach spaces have gradually become a standard tool in Banach space local theory (see, for example, Appendix F of Benyamini and Lindenstrauss [14]).

Often both of the above approaches to extracting limit objects are possible. They differ in their merits. The first, more hands-on approach requires a careful check on whether the sequence of objects converges at all (unless one finds that there is some compactness in the ambient object to deduce the existence of some suitably convergent subsequence, as happens in the case of exchangeable random hypergraphs). On the other hand when ultralimits make sense at all their definition is typically quite easy; however, the limit objects that result may be quite unnatural (ultralimits of separable Banach spaces are typically hugely non-separable, for example), and may require further manipulation before they really represent the asymptotic data that was sought. Sometimes, this latter manipulation then leads naturally back to the same limit object as could have been extracted from the first approach; however, in some settings one does not need such precise information about the limit object, and in this case the ultralimit construction may often be much faster.

The situation for dense hypergraphs and exchangeable random hypergraphs seems to be similar. As discussed previously, our construction of sampling random hypergraphs effectively gives an embedding of finite hypergraphs into the single ambient space of exchangeable probability measures on $\{0,1\}^{\binom{S}{k}}$ for some arbitrary countably infinite $S$, and we are now able to take limits quite naturally in the vague topology of measures and identify this with convergence of subhypergraph densities for the original finite hypergraphs. Thus, exchangeable random hypergraphs can serve as limit objects in the first of our two senses above, and with the structure theorem and a little more work it is possible to describe asymptotic features of their leading-order statistics in terms of them.

Indeed, the development of this analogy has an older precedent within combinatorics: that established by the 'objective method' in the study of certain enumeration or optimization problems over large finite networks, surveyed, for example, in Aldous and Steele [5]. The basic finitary objects here are thought of as large networks but with (usually) small individual vertex degrees, and often endowed with a distinguished vertex called the 'root'. Often one is interested in the typical connectivity or other statistics of fixed-radius neighbourhoods around a randomly-chosen root, or of flows through the network, and in this case the relevant numbers can often be extracted into a suitable 'infinite random network' amenable to separate study. We direct the reader also to Benjamini and Schramm [12], Aldous and Lyons [4] and Elek [22, 23] for a view of various instances of this approach, including connections with the theory of graphings of Borel equivalence relations that would take us too far afield here. This approach has recently also entered the arena of property testing in a suitable class of rooted graphs subject to a fixed maximal degree bound (we will treat prop-erty testing for dense graphs and hypergraphs in Subsection 4.5 below) in works of Schramm [59] and Benjamini, Shapira and Schramm [13].

### 4.2. Ultralimits and the work of Elek and Szegedy

While the use of exchangeable laws as limit objects for finite hypergraph colourings amounts to an embedding of finite hypergraphs into a suitable compact metric space, then the approach via ultralimits is also possible. The (easy) extraction and (harder) structural description of these ultralimits has been carried out by Elek and Szegedy in [24], where they too are able to recover certain 'leading order' combinatorial results about finite hypergraphs from general properties of the limit object, including a version of the hypergraph removal lemma.

In fact, Elek and Szegedy's use of ultralimits has a predecessor in Hoover's approach to the basic representation theorems for exchangeable arrays of random variables, as outlined in [42]. However, Hoover uses the ultralimit construction to extract the ingredients for a given exchangeable array, rather than as a possible route to defining limit objects for sequences of finite graphs or hypergraphs.

The analysis of Elek and Szegedy includes a manipulation of the (initially huge) non-separable limit-object to give a separable version capturing essentially the same information; however, the separable limit object that results is not immediately equivalent to our structural ingredients for an exchangeable random hypergraph. Reformulated as a structural result for exchangeable random hypergraphs, their theorem reads as follows.

**Theorem 4.1** (Elek-Szegedy Theorem for exchangeable random hypergraphs). *Suppose that $\mu$ is an exchangeable random $\{0,1\}$-coloured $k$-uniform hypergraph. Then there is some measurable $W : [0,1]^{2^k-1} \to [0,1]$ such that, for any finite-dimensional cylinder set $A$ in $\{0,1\}^{\binom{S}{k}}$ given by*

$$A = \{\omega \in \{0,1\}^{\binom{S}{k}} : \omega|_{\binom{V}{k}} \geq 1_F\}$$

*for some finite $k$-uniform hypergraph $(V, F)$ with vertex set $V \subset S$, we have*

$$\mu(A) = \int_{[0,1]^{\binom{V}{<k}}} \prod_{e \in F} W(x|_{\binom{e}{<k}}) \, \mu_{\mathrm{L}}(\mathrm{d}x),$$

*where $\mu_{\mathrm{L}}$ is Lebesgue measure on $[0,1]^{\binom{V}{<k}}$.*

(We will not consider here the extension of this result to the case of $K$-coloured hypergraphs, but it does not seem to offer serious difficulties.)

We will re-derive this from our structure theorem below, and refer to [24] directly for the more combinatorial consequences of their ultralimit construction. Interestingly, the cubes with Lebesgue measure $[0,1]^{2^k-1}$ also arise explicitly in Aldous' and Kallenberg's chosen versions of the representation theorem for exchangeable (and spreadable) arrays or set-indexed processes ([3, 44]), as the target spaces of collections of independent uniform $[0,1]$-valued random variables.

As promised during our discussion of their versions in Subsection 3.7 above, our deduction of Theorem 4.1 from the structure theorem can be modified directly to deduce the equivalence of Corollary 3.5 and Theorem 3.6, Theorems 2.9 and 3.14 or Theorems 3.24 and 3.25.

We need to remove the complicated structural ingredients $(Z_0, \mu_0)$, $(Z_1, P_1)$, ..., $(Z_k, P_k)$ in favour of a single measurable function $[0,1]^{2^k-2} \to [0,1]$ and Lebesgue measure, at the expense of introducing a measurable function $W$ as above. This reduction loses the description of all the intermediate systems built from the spaces $Z_i^{\binom{S}{i}}$, and so, while quicker for deriving various combinatorial consequences, it seems to suppress some of the probabilistic features of the analysis.

Our recovery of the above version of the Elek-Szegedy result will be based on the following standard fact of probability; see, for example, Theorem 5.10 in Kallenberg [45]. We recall a proof for completeness.

**Lemma 4.2** (Noise outsourcing lemma). *Suppose that $X$ and $Y$ are standard Borel spaces, that $\mu$ is a probability measure on $X$ and that $P : X \rightsquigarrow Y$ is a probability kernel. Let $\mu_L$ be Lebesgue measure on $[0,1]$. Then there is some Borel measurable map $f : X \times [0,1] \to Y$ such that, endowing $X \times [0,1]$ with the product measure $\mu \otimes \mu_L$, the kernel $P(x, \cdot)$ is a version of the conditional distribution of $f(x, \cdot)$ given the first coordinate $x$: that is, for any $A \in \Sigma_X$ and $B \in \Sigma_Y$ we have*

$$P(x, B) = \mathsf{E}_{\mu \otimes \mu_L}[1_{\{f(x,t) \in B\}} \,|\, x]$$

*up to $\mu$-almost everywhere equivalence.*

**Remark.** If we think of the kernel $P(x, \cdot)$ as specifying a random choice of a point $y \in Y$ with law depending on a point $x \in X$, then this lemma tells us that the randomness in this choice can be correctly represented by first choosing independently and uniformly a value for the 'noise parameter' $t \in [0,1]$, and then choosing $y$ according to some deterministic $Y$-valued map on pairs $(x,t)$. ◁

*Proof.* Starting from a probability kernel from $X$ to $Y$, we will need to build a deterministic function with values in $Y$: for this it will be crucial that $Y$ is a standard Borel space, in order that we can specify individual points of it in some measurable way in our definition of $f$. Given this, we can proceed in a number of ways. For example, we can first identify $Y$ with a Borel subset of $[0,1]$, and can now specify a suitable $f$ such that for each $x$ the map $t \mapsto f(x,t)$ is non-decreasing by setting

$$f(x,t) := \sup\{q \in \mathbb{Q} \cap [0,1] : \ P(x, [0,q]) \leq t\}.$$

It is routine to check that this countable supremum is measurable (appealing to the measurability in $x$ of the kernel $P$) and that it has the desired distribution (and so must, in particular, take values in $Y$ almost surely). This construction is essentially a pointwise-in-$x$ implementation of the Skorokhod embedding. □

Given this, we can quickly reconstruct the Elek-Szegedy representation of Theorem 4.1:

*Proof of Theorem 4.1.* We will replace the spaces $Z_i$ by powers of the unit interval level-by-level from below. At stage $i \leq k-1$ we lose some of the information contained in a point of $Z_i$, but work instead with the value of some outsourced noise parameter in $[0,1]$, as given by Lemma 4.2; this will allow us to retain up to stage $i$ the distribution $\mu_i$ on $Z_i^{\binom{S}{i}}$ while replacing all of the lower-rank distributions with simpler outsourced noise. Thus we obtain a sequence of intermediate towers of quasifactors of $\mu$ (after initially also identifying $(Z_0, \mu_0)$ with $([0,1], \mu_L)$), represented by the rows in the following array:

$$\begin{array}{lllllll}
[0,1] \stackrel{P_1}{\rightsquigarrow} Z_1 & Z_{<2} \stackrel{P_2}{\rightsquigarrow} Z_2 & Z_{<3} \stackrel{P_2}{\rightsquigarrow} Z_3 & \ldots & Z_{<k} \stackrel{P_k}{\rightsquigarrow} \{0,1\} \\
[0,1]^{1+1} \stackrel{f_1}{\longrightarrow} Z_1 & Z_{<2} \stackrel{P_2}{\rightsquigarrow} Z_2 & Z_{<3} \stackrel{P_2}{\rightsquigarrow} Z_3 & \ldots & Z_{<k} \stackrel{P_k}{\rightsquigarrow} \{0,1\} \\
& [0,1]^{3+1} \stackrel{f_2}{\longrightarrow} Z_2 & Z_{<3} \stackrel{P_2}{\rightsquigarrow} Z_3 & \ldots & Z_{<k} \stackrel{P_k}{\rightsquigarrow} \{0,1\} \\
& & [0,1]^{7+1} \stackrel{f_2}{\longrightarrow} Z_3 & \ldots & Z_{<k} \stackrel{P_k}{\rightsquigarrow} \{0,1\} \\
& & & \vdots & \\
& & & & [0,1]^{2^k-1} \stackrel{W}{\rightsquigarrow} \{0,1\}
\end{array}$$

where we have written $Z_{<i} := \prod_{a \subset [i]} Z_{|a|}$. When we reach stage $k$ we have a representation of $\mu_k = \mu$ purely in terms of this noise and one last probability kernel; this will be represented by the function $W$ of the Elek-Szegedy Theorem. □

### 4.3. Measures on spaces of isomorphism classes and the work of Razborov

Much of the work of this paper has been to describe those probability measures $\mu$ on a space $K^{(S)}$ of colourings by a $k$-palette $K$ of either the subsets or the directed hyperedges of a vertex set $S$ that are invariant under finitely-supported vertex permutations. Any two points of this space lie in the same orbit of this $\mathrm{Sym}_0(S)$-action if and only if they are isomorphic as (directed) coloured hypergraphs by a finitely-supported vertex permutation, and so we may regard the study of such $\mu$ as a natural alternative to working directly with probability measures on the space of hypergraph colouring isomorphism classes, $K^{(S)}/\mathrm{Sym}_0(S)$.

Indeed, the quotient by $\mathrm{Sym}_0(S)$ destroys all of the nice topological and Borel space structure of the space $K^{(S)}$, and so working directly on this space of isomorphism classes is more difficult. This topological difficulty seems to appear in some guise or other in any attempt to study a limit object for hypergraph colourings, even though it might seem more natural a priori to consider these only up to isomorphism.

Exchangeable random colourings of $\binom{S}{\leq k}$ or $\bigcup_{i \leq k} \mathrm{Inj}([i], S)$ are a convenient alternative requiring only very classical probabilistic ideas. However, there are other ways to circumvent this difficulty. Here we will give an overview of one such, and compare it with the study of exchangeability. This is the approach the underlies recent work of Razborov [54, 55], and rests on a conversion of the necessary data into certain very specialized abstract commutative algebras ('flag algebras') and $\mathbb{R}$-valued homomorphisms on them. In fact, Razborov details his construction in the slightly different lexicon of model theory (much as outlined briefly in Subsection 3.8 above), and we have made a partial translation into that of hypergraph colourings in order to make the comparison with the results of the present note simpler. In addition, Razborov is by no means the first to consider a notion of exchangeable probability measure on a space of models of a simple theory; for example, these appear already in the studies by Gaifman [34] and Krauss [51], and their work is then offered as motivation in Hoover's article [42]. Rather, the novelty of Razborov's approach is his use to which he puts these ideas in the pursuit of certain purely finitary combinatorial results; given the relevance of this to the aims of the present survey, we shall concentrate on Razborov's chosen route through the theory below.

Just as for many of the other constructions described in this survey, the study of these is motivated by their ability to serve as proxies for the leading-order statistics of large graphs or hypergraphs. Razborov's approach emphasizes first an abstractly-defined commutative algebra over $\mathbb{R}$ defined in terms of isomorphism classes of finite hypergraphs, and then those of its linear homomorphisms to $\mathbb{R}$ that enjoy a certain positivity property. It turns out that these algebras correspond in our picture precisely to certain quotients of the algebra of continuous functions on $K^{(S)}$, and that routine soft functional analytic arguments now identify the positive homomorphisms of this algebra with the exchangeable probability measures on this space that are 'effectively concentrated on a single equivalence class' (that is, are ergodic). This is how Razborov's algebras and their homomorphisms give an alternative means for handling the difficulty of defining a probability measure directly onto the space of these equivalence classes.

Here we will discuss the relation between Razborov's formalism (set up in the early sections of [54]) and that of exchangeable random hypergraph colourings. For simplicity we will specialize to the case of undirected $K$-coloured hypergraphs for some $k$-palette $K$, so that the notation $K^{(S)}$ will be used as in Subsection 3.3. The new notation we shall introduce to describe flag algebras follows [54]. Razborov's work actually applies in the quite general setting of models of theories outlined in Subsection 3.8; a more detailed discussion in this greater generality can be found in [9].

Razborov goes on to use his machinery to conduct a variational analysis of a problem in graph theory, an instance of the 'graph copy problem' relating to edge-versus-triangle densities in large graphs (see also Subsection 4.6 below). This leads to a quite specific optimization problem (in some ways similar to that appearing in Fremlin and Talagrand [29]) studied in Sections 4 and 5 of [54], and then analyzed completely in [55]. While they too can presumably be translated into the measure theoretic picture, this actually seems to change the arguments involved very little.

Let us first consider isomorphism classes of those finite $K$-colourings that contain a distinguished copy of some fixed finite $K$-colouring $\sigma$; that is, pairs $F = (M, \theta)$ with $M$ a finite $K$-colouring and $\theta$ an embedding $\sigma \hookrightarrow M$. These isomorphism classes of extensions are examples of the '$\sigma$-flags' of [54], and $\mathcal{F}^\sigma$ is written for the collection of all them, $\mathcal{F}_\ell^\sigma$ for the subcollection of those with vertex set of size $\ell$ (so $\mathcal{F}_\ell^\sigma = \emptyset$ unless $\ell \geq k$), and we can now specify the obvious notions of embedding and isomorphism for $\sigma$-flags.

Let us restrict our attention to $\sigma$-flags $F = (M, \sigma)$ with $M \in K^{([\ell])}$ for some $\ell \geq 1$ and $\sigma$ itself serving as its distinguished copy in $M$ (so we implicit order the vertices of $M$ so that $\sigma = M|_{[r]}$ for some $r$). Let $\Omega := K^{(\mathbb{N})}$, let $\Omega_\sigma$ be the subset of those $\omega \in \Omega$ such that $\omega|_{[r]} = \sigma$, and now associate to each $F \in \mathcal{F}^\sigma$ the finite-dimensional cylinder set

$$A_F := \{\omega \in \Omega : \omega|_{[\ell]} = M\} \subseteq \Omega_\sigma$$

(so $\Omega_\sigma = A_{(\sigma,\sigma)}$; the difference in notation reflects the different rôles of $\sigma$ and $F$).

Let us write $C(\Omega_\sigma)$ and $\mathcal{M}(\Omega_\sigma)$ for the usual Banach spaces of real-valued continuous functions and signed Radon measures on $\Omega_\sigma$ respectively; the Riesz-Kakutani representation identifies $\mathcal{M}(\Omega_\sigma)$ isometrically with $C(\Omega_\sigma)^*$. Moreover, we write $\mathcal{M}^\sigma$ for the subspace of measures supported on $\Omega_\sigma$ and invariant under finitely-supported permutations of the vertices in $\mathbb{N} \setminus [k]$; we will call these $\sigma$-**exchangeable**, and denote the group of these permutations by $\mathrm{Sym}_\sigma$ (so this is just the stabilizer of $1, 2, \ldots, k$ in $\mathrm{Sym}_0(\mathbb{N})$). We write $(\mathcal{M}^\sigma)^\perp$ for the annihilator of this space of continuous linear functionals in $C(\Omega_\sigma)$,

$$(\mathcal{M}^\sigma)^\perp = \{f \in C(\Omega_\sigma) : \langle f, \mu \rangle = 0 \ \forall \mu \in \mathcal{M}^\sigma\};$$

as usual, the dual-of-the-quotient Banach space $\big(C(\Omega_\sigma)/(\mathcal{M}^\sigma)^\perp\big)^*$ can be isometrically identified with $\mathcal{M}^\sigma$. Let $q_\sigma$ be the quotient map $C(\Omega_\sigma) \to C(\Omega_\sigma)/(\mathcal{M}^\sigma)^\perp$.

Given any $f \in C(\Omega_\sigma)$, we must have $f - f \circ g \in (\mathcal{M}^\sigma)^\perp$ for all $g \in \mathrm{Sym}_\sigma$. Write $\mathrm{T}_\sigma$ for the tail $\sigma$-subalgebra $\bigcap_{m \geq k+1} \sigma(\pi_{[k] \cup \{m, m+1, \ldots\}})$, and (with a slight abuse of notation) $L^\infty(\mathrm{T}_\sigma)$ for the bounded $\mathrm{T}_\sigma$-measurable functions under the equivalence relation of "equality $\mu$-a.e. for every $\mu \in \mathcal{M}^\sigma$". Clearly these are invariant under the action of $\mathrm{Sym}_\sigma$. By (for example) the pointwise ergodic theorem for the amenable group $\mathrm{Sym}_\sigma$, the ergodic averages of the compositions $f \circ g$ over $g$ converge to a $\mathrm{T}_\sigma$-measurable function $\bar{f}$ on $\Omega_\sigma$ which is defined $\mu$-almost everywhere for every $\mu \in \mathcal{M}^\sigma$ and is invariant under $\mathrm{Sym}_\sigma$, and hence actually specifies a member of $L^\infty(\mathrm{T}_\sigma)$. Observe that $\bar{f} = \bar{h}$ for $f, h \in C(\Omega_\sigma)$ if and only if $f - h \in (\mathcal{M}^\sigma)^\perp$, and so our map $f \mapsto \bar{f}$ embeds $C(\Omega_\sigma)/(\mathcal{M}^\sigma)^\perp$ as a subspace $V^\sigma$ of $L^\infty(\mathrm{T}_\sigma)$; general nonsense now shows also that this is an isometric embedding, so $V^\sigma$ is closed.

Furthermore, $V^\sigma$ is actually a sub*algebra* of $L^\infty(\mathrm{T}_\sigma)$. To see that it is closed under multiplication, suppose $f, h \in C(\Omega_\sigma)$, and now consider the products $f \cdot (h \circ g)$ for any sequence of permutations $g$ that pushes $h$ 'further and further out', in the following sense: for any $m \geq 1$, there are finite $A, B \subset \mathbb{N} \setminus [r]$ such that $f$ and $h$ are uniformly $(1/m)$-close to functions depending only on the colours of edges above vertices in $A$ and $B$ (respectively), and now we insist that $g$ move $B$ into $\mathbb{N} \setminus (A \cup [r])$. Letting $m \to \infty$ this gives a sequence $g_m$ for which, in terms of their dependence on coordinates, $f$ and $h \circ g_m$ are closer and closer to independent. Now an elementary argument using approximation by step functions shows that the quotients $q_\sigma(f \cdot h \circ g)$ converge in $C(\Omega_\sigma)/(M^\sigma)^\perp$ to a member that depends only on $q_\sigma(f)$ and $q_\sigma(h)$; and it is a routine exercise to check that this actually defines a C*-algebra product on $C(\Omega_\sigma)/(M^\sigma)^\perp$ corresponding exactly to the usual product of functions in $V^\sigma$.

Thus we have identified $C(\Omega_\sigma)/(\mathcal{M}^\sigma)^\perp$ with a closed subalgebra $V^\sigma$ of $L^\infty(\mathrm{T}_\sigma)$ (with a newly-defined product). Let us call functions in $V^\sigma$ **simple** if they are the images of simple (equivalently, finite-dimensional) functions in $C(\Omega_\sigma)$. We can describe the simple functions naturally as follows: to any fixed nonempty cylinder set $A \subseteq \Omega_\sigma$ depending on coordinates in $J \subset \mathbb{N} \setminus [k]$ corresponds a collection of finite models of the theory $T$ on the vertex set $J$ (with some multiplicities), and now the averaged-function $\overline{1_A}(\omega)$ for $\omega \in \Omega_\sigma$ is just the sum of the densities with which each of those finite models appears isomorphically as a submodel of $\omega$ (now summing over the multiplicities). Referring to such a function $\overline{1_A}$ for $A$ corresponding to a single model on $J$ (so that our general $A$ is a disjoint union of such) as a **statistics function**, the simple functions in $V^\sigma$ are now just linear combinations of statistics functions. We write $V_0^\sigma$ for the dense subspace of these.

We now begin our account of flag algebras themselves. Let $\mathbb{R}\mathcal{F}^\sigma$ denote the free $\mathbb{R}$-vectorspace on $\mathcal{F}^\sigma$. We modify the correspondence $F \mapsto A_F \subseteq \Omega_\sigma$ introduced above to associate to a member $\sum_j \lambda_k F_j \in \mathbb{R}\mathcal{F}^\sigma$ a corresponding linear combination of the indicator functions $1_{A_F}$ of these $A_F$, with some suitable renormalizing constraints (which depend on from Razborov's choice of formalism):

$$F \mapsto \frac{1}{|F|!} 1_{A_F};$$
$$\sum_j \lambda_k F_j \mapsto \sum_j \frac{\lambda_k}{|F_j|!} 1_{A_{F_j}}.$$

The step functions that appear above are continuous on $\Omega_\sigma$ since each $A_F$ is clopen, and so this defines a linear operator $\Phi : \mathbb{R}\mathcal{F}^\sigma \to C(\Omega_\sigma)$ with image some peculiar subspace contained within the space of simple functions in $C(\Omega_\sigma)$.

The arbitrariness in our choice of the subsets $A_F$ corresponding to $F$ is reflected in a similar arbitrariness in the linear map $\Phi$ into $C(\Omega_\sigma)$; however, this disappears at the next step, when we define a flag algebra as a quotient of $\mathbb{R}\mathcal{F}^\sigma$. Let $\mathcal{K}^\sigma$ be the subspace of $\mathbb{R}\mathcal{F}^\sigma$ generated by the linear combinations $\tilde{F} - \sum_{F \in \mathcal{F}_\ell^\sigma} p(\tilde{F}, F) F$ for different $\ell \geq |V(\tilde{F})|$, and for certain real numbers $p(\tilde{F}, F)$ (given in Definition 1 of [54]). It turns out that, given our chosen normalization in the definition of $\Phi$ above, the values $p(\tilde{F}, F)$ are such that $\Phi(\mathcal{K}^\sigma)$ is precisely the set of those $a \in \mathbb{R}\mathcal{F}^\sigma$ for which $\Phi(a) \in C(\Omega_\sigma)$ is annihilated by every exchangeable probability measure on $\Omega_\sigma$, and thus by all finite signed measures obeying the vertex-permutation symmetry (this is discussed in a little more detail in [9]).

We now consider the quotient space $\mathfrak{A}^\sigma := \mathbb{R}\mathcal{F}^\sigma/\mathcal{K}^\sigma$, this is Razborov's **flag algebra**. By the above, $q_\sigma \circ \Phi$ factors through this quotient to give an injective map $\Psi : \mathfrak{A}^\sigma \to V^\sigma$. Moreover, since any finite dimensional cylinder set contained in $\Omega_\sigma$ can be identified with some $A_F$ upon a suitable permutation of coordinates in $\mathbb{N} \setminus [k]$, the image of $\Psi$ is actually the subspace $V_0^\sigma$ of *all* simple functions in $V^\sigma$; as such, it is dense.

Under the above identification the image of Razborov's definition of the product of $a, b \in \mathbb{R}\mathcal{F}^\sigma/\mathcal{K}^\sigma$ is now the product of $\Psi(a)$ and $\Psi(b)$ as $L^\infty$-functions on $\Omega_\sigma$; the proof in [54] that this product is well-defined translates into a proof in the exchangeability picture that this product remains in the image of $\Psi$ (that is, in $V_0^\sigma$). This completes our identification of the flag algebra $\mathfrak{A}^\sigma$ with the dense subalgebra $V_0^\sigma$ of $V^\sigma$, which is itself a norm-closed Banach subalgebra of $L^\infty(\mathrm{T}_\sigma)$.

The overall approach during the early stages of [54] is to define flag algebras *first* (Section 2), and then to specify a collection of 'limit objects' for the statistics of large models of a theory as certain homomorphisms of these flag algebras (Section 3). Thus, the next step is to consider multiplicative functions $\phi : \mathfrak{A}^\sigma \to \mathbb{R}$ that are non-negative on the image of any single flag $F \in \mathcal{F}^\sigma$; the set of these is written $\mathrm{Hom}^+(\mathfrak{A}^\sigma, \mathbb{R})$. Now, having identified $\mathfrak{A}^\sigma$ with the dense subalgebra $V_0^\sigma \subseteq V^\sigma$ so that the images of single flags correspond to the single statistics functions, we can easily check that non-negativity on these implies non-negativity on any member of $V_0^\sigma$ which is itself non-negative as a real-valued function. Therefore, given the non-negativity of such a $\phi$, it follows that it must be bounded as a linear functional on $V_0^\sigma$, and so can be extended to a multiplicative linear functional on $V^\sigma$.

We could now, if we wished, apply certain standard representation theorems to this space (perhaps most directly the results of Yosida and Kakutani for M-spaces, as presented in Section XII.5 of [66]). This would identify $\phi$ with a point of the spectrum of $V^\sigma$ (for one or other interpretation of 'spectrum'). In fact, a similar observation is more-or-less implicit in Remark 4 of Subsection 3.2 of [54], although there we still require some of the basic structure of $\mathrm{Hom}^+(\mathfrak{A}^\sigma, \mathbb{R})$ to have been identified.

However, given our identification of $\mathfrak{A}^\sigma$ with $V_0^\sigma$, we have an alternative to the above. Since $V^\sigma \cong C(\Omega_\sigma)/(\mathcal{M}^\sigma)^\perp$, as a linear functional on $V^\sigma$ we can identify $\phi$ with a member of $\mathcal{M}^\sigma$ (uniquely, since $(C(\Omega_\sigma)/(\mathcal{M}^\sigma)^\perp)^* \cong \mathcal{M}^\sigma$): a $\mathrm{Sym}_\sigma$-invariant measure on the space $\Omega_\sigma$ that we started with, rather than a point of some abstractly-produced new space $\mathrm{Spec}\, V^\sigma$. It is now easy to check that those measures in $\mathcal{M}^\sigma$ that are *multiplicative* on $V^\sigma$ are precisely the ergodic $\mathrm{Sym}_\sigma$-invariant probability measures on $\Omega_\sigma$ (since any member of $V^\sigma$ is $\mu$-a.s. constant if $\mu \in \mathcal{M}^\sigma$ is ergodic). Now the order defined on $\mathfrak{A}^\sigma$ in Definition 5 of Section 3 of [54] is precisely the usual pointwise order on $V_0^\sigma$ as a set of real-valued functions; in the setting of abstract flag algebras, where we are unable to define anything 'pointwise', the functionals of $\mathrm{Hom}^+(\mathfrak{A}^\sigma, \mathbb{R})$ are needed as a replacement to formulate this definition. This constitutes the promised correspondence between exchangeable random $K$-coloured hypergraphs (on $\mathbb{N}$, say) and the homomorphisms of Razborov's theory.

Let us finish this subsection by noting that while Razborov's work seems to come closest to developing a theory of 'leading-order statistics' in terms of probability measures or other structures defined directly on a space of isomorphism classes, there are precedents for such an approach among the study of certain other combinatorial structures: again, these can be found in the 'objective method' mentioned in passing in Subsection 4.1 above and recounted in Aldous and Steele [5]. Perhaps the simplest instance of the extraction of a limit object as a measure defined directly on a space of isomorphism classes of infinite models occurs in the work on limit objects for various classes of sparse rooted graph (see, in particular, the papers of Aldous and Steele [5], Aldous and Lyons [4], Benjamini and Schramm [12], Benjamini, Schramm and Shapira [13] and Elek [23]). For example (as in [12]), in the setting of all rooted graphs $(G, o)$ satisfying some fixed maximum-degree bound $d \geq 1$, one can associate to such a graph $G$ a probability measure $\mu_G$ directly on the set of all isomorphism classes of such rooted graphs of maximum degree at most $d$, by defining the probability of the class of a given such graph $(H, v)$ to be the proportion of possible roots $o \in V(G)$ for which $(G, o) \cong (H, v)$. The set of all such isomorphism classes can be endowed with a natural compact metrizable topology by defining a neighbourhood base as those subsets of equivalence classes for which a ball of a given finite radius around the root lies within a given isomorphism class, and we may now take vague limits of these measures $\mu_G$ for suitable sequences of $G$. This ability to study a natural — albeit, perhaps, not easily visualized — compact metrizable topology on the space of isomorphism classes distinguishes this setting from that of dense graphs. On the other hand, in the latter setting the possibility of working in the much larger space of all $k$-uniform hypergraphs on $\mathbb{N}$ (*not* quotiented by the isomorphism relation) is available, and makes contact with the classical theory of exchangeability; it is not clear that the setting of bounded-degree graphs admits a similarly nice class of 'labeled' models on which we could study those probability measures invariant under automorphisms.

### *4.4. Comparison with finitary regularity lemmas*

Since Szemerédi's introduction of his original graph regularity lemma on the way to his proof of Szemerdí's Theorem ([61]), it has become perhaps the single most powerful tool in extremal graph theory (see, for example, Sections IV.5 and IV.6 of Bollobás [16]). More recently, various regularity lemmas have emerged for hypergraphs with enough strength to reproduce many of the results from the graph case: in particular, those of Nagle, Rödl and Schacht [53] and of Gowers [36]. These results can be described loosely as guaranteeing that to an arbitrary very large $k$-uniform hypergraph there must correspond some smaller non-uniform hypergraph, of order bounded only in terms of an a priori fixed error tolerance, and some $[0, 1]$-weights on its hyperedges, such that a certain normalized weighted count of embedded copies of a given small hypergraph in the latter is a good approximation to their normalized count (the 'density') in the former (we shall not make this more precise here). They all rest on one or another notion of 'quasirandomness' for hypergraphs.

It emerges naturally in all the higher-rank regularity lemmas that in order to obtain a good approximation to the densities of small embedded subgraphs in a very large fixed $k$-uniform hypergraph $(V, H)$, one must first introduce a partition of the vertex set that controls only certain statistics pertaining to the $(k-1)$-hyperedges over $V$, and then given this information obtain a much finer partition that controls statistics pertaining to $H$ and to this first partition through the $(k-2)$-hyperedges over $V$, and so on; the final partition of $V$ is obtained after $k$ such steps. This feature is present in different ways in the Nagle-Rödl-Schacht and Gowers regularity lemmas. Gowers' approach, which is tailored very precisely to the needs of proving hypergraph removal lemmas (and thus also Szemerédi's Theorem), rests on the use of 'Gowers norms' to control the count of 'octahedra' in a hypergraph (for which we introduced infinitary analogs in Subsection 3.7 above). On the other hand, Nadle, Rödl and Schacht obtain an iterative control of the correlations of the original $k$-uniform hypergraph with $\ell$-uniform hypergraphs for $\ell \leq k$, starting at $\ell = k$ and working downwards.

A clear discussion of the similarities and differences in these two finitary hypergraph regularity lemmas can be found around Definitions 2.5 and 2.6 and in Section 10 of Gowers [38].

We recall these observations here because this picture of exerting control from the top rank downwards — particulary in the form it takes in Nagle, Rödl and Schacht's work — is reminiscent of the downwards-inductive construction of our ingredients, in which each of the resulting probability kernels must be allowed to depend on type spaces of all ranks beneath it.

### *4.5. Property testing, repairability and joinings*

A more specific class of questions that can be set up roughly in parallel for finite (directed or undirected) hypergraph colourings and for their exchangeable random counterparts is that of testing and local repair of hereditary hypergraph properties. The basic questions of property testing have attracted considerable attention during the last ten years, and we shall discuss them only very incompletely here; see, for example, the papers of Rubinfeld and Sudan [57], Alon and Shapira [8] and Alon, Fischer, Krivelevich and Szegedy [6], and the further references given there.

In this subsection we shall work with coloured directed hypergraphs, as we will in [10]. Suppose that $K$ is a finite $k$-palette and that $\mathcal{P}$ is a property of $K$-coloured directed hypergraphs: that is, it is a subcollection of the collection of all isomorphism classes of finite $K$-coloured directed hypergraphs. In this case, given a vertex set $V$ we shall write $\mathcal{P}^{(V)}$ for the set of all $K$-coloured hypergraphs on $V$ satisfying $\mathcal{P}$ (that is, members of $K^{(V)}$ whose isomorphism class is in $\mathcal{P}$), and shall write that the members of $\mathcal{P}^{(V)}$ **obey** $\mathcal{P}$.

**Definition 4.3** (Hereditary property)**.** A property $\mathcal{P}$ as above is called **hereditary** if whenever $G \in \mathcal{P}^{(V)}$ and $W \subseteq V$ then the induced $K$-coloured hypergraph $G|_W$ on $W$ is in $\mathcal{P}^{(W)}$.

In case $\mathcal{P}$ is hereditary, we may naturally extend its definition to coloured directed hypergraphs $G$ on infinite vertex sets $V$ by specifying that $G$ obeys $\mathcal{P}$ if and only if all of its finite induced subgraphs do.

**Definition 4.4** (Testability). [57] Let $K$ be a finite $k$-palette and let $\mathcal{P}$ be a hereditary $K$-coloured directed hypergraph property. We say that $\mathcal{P}$ is **testable with one-sided error** if for every $\varepsilon > 0$ there are an integer $N \geq 1$ and a real number $\delta > 0$ for which the following holds: if $G \in K^{(V)}$ is a $K$-coloured hypergraph on a vertex set $V$ with with $N \leq |V| < \infty$, and $G$ 'locally almost obeys $\mathcal{P}$' in the sense that

$$\frac{1}{|\binom{V}{N}|}\left|\left\{W \in \binom{V}{N} : G|_W \in \mathcal{P}^{(W)}\right\}\right| \geq 1 - \delta,$$

then there exists $G' \in \mathcal{P}^{(V)}$ which is close to $G$ in the sense that

$$\frac{1}{|\binom{V}{k}|}\left|\left\{e \in \binom{V}{k} : (G_\phi)_{\phi \in \mathrm{Inj}([k],e)} \neq (G'_\phi)_{\phi \in \mathrm{Inj}([k],e)}\right\}\right| \leq \varepsilon.$$

In the case of undirected $\{0,1\}$-coloured graphs or hypergraphs, it has recently been shown in work of Alon and Shapira [7] and Rödl and Schacht [56] that every hereditary property is testable. Their arguments are purely finitary, each resting on several applications of a graph or hypergraph regularity lemma to 'process' a given hypergraph (although the details of their approaches are very different). However, reconsidering these results in the setting of infinitary exchangeable random hypergraph colourings shows some interesting features, and also proves convenient for several further extensions of the combinatorial results to directed, non-uniform hypergraph colourings with larger finite palettes. This approach will be elaborated in detail in [10], where the main positive result is the assertion that all hereditary properties in the full generality of Definition 4.4 are testable.

In [10] we will also introduce certain related notions of 'local repairability', which amount to additional conditions on the manner in which a suitable modification $G'$ is to be found in Definition 4.4. In particular, they demand that there be some randomized procedure for modifying $G$ to $G'$ that modifies the colour of a given edge in $G$ according to the behaviour of $G$ across only that edge and some bounded number of other vertices, these latter chosen at random, and that yields a suitable modification with high probability. These notions of 'local repairability' take various forms, both in terms of finite coloured directed graphs or hypergraphs and with related versions for their exchangeable random cousins. We will not elaborate on the different forms of 'local repairability' here, referring the reader to [10] instead, but will restrict ourselves to an informal statement of the infinitary counterpart of simple testability as illustration of these relationships.

Thus, let $K$ be a finite $k$-palette, let $\mathcal{P}$ be a hereditary property of $K$-coloured directed hypergraphs. We will write that $\mathcal{P}$ is **infinitarily testable with one-sided error** if whenever $\mu$ is an exchangeable random $K$-coloured directed hypergraph on a countably infinite vertex set $S$ whose sample points obey $\mathcal{P}$ almost surely and $\varepsilon > 0$ is an error tolerance then the following holds. If we extract the structural ingredients guaranteed by Theorem 3.23:

- a $k$-palette $Z = (Z_i)_{i=0}^k$ consisting of totally disconnected compact metric spaces,
- continuous maps $\kappa_i : Z_i \to K_i$,
- and an exchangeable random $Z$-coloured hypergraph $\mu'$ under which the clusters of random variables $(\pi_\phi^Z)_{\phi \in \mathrm{Inj}([i],u)}$ for distinct $u \in \binom{S}{i}$ are relatively independent when conditioned on all the random variables $\pi_\psi^Z$ with $\psi \in \mathrm{Inj}(<i,S)$, and $\mu'$ satisfies $(\kappa_0, \kappa_1, \ldots, \kappa_k)_\#^{(S)} \mu' = \mu$, then there exists a $\mathrm{Sym}_0$-covariant probability kernel $T : Z^{(S)} \rightsquigarrow K^{(S)}$ such that

- the map $(z_\emptyset, (z_i)_{i \in e}, \ldots) \mapsto T(z_\emptyset, (z_i)_{i \in e}, \ldots, \cdot)$ is vaguely continuous from $Z^{(S)}$ into $\mathrm{Pr}\, K^{(S)}$,
- for each $W \subseteq S$, the projected measure $(\pi^K_{\bigcup_{i \leq k} \mathrm{Inj}([i],S)})_\# (T(z_\emptyset, (z_i)_{i \in S}, \ldots, \cdot))$ depends only on $(z_\emptyset, (z_i)_{i \in W}, \ldots)$,
- the probability measure $T(z_\emptyset, (z_i)_{i \in S}, \ldots, \cdot)$ is concentrated on $\mathcal{P}^{(S)}$ for every $(z_\emptyset, (z_i)_{i \in S}, \ldots) \in Z^{(S)}$,
- and $T$ is close to $\kappa$ under $\mu'$ in the sense that

$$\int_{Z^{(e)}} T\Big(z, \Big(\pi^K_{\bigcup_{i \leq k} \mathrm{Inj}([i],e)}\Big)^{-1} (K^{(e)} \setminus \{(\kappa_0(z_\emptyset), (\kappa_1(z_i))_{i \in e}, \ldots)\})\Big)$$
$$\times \Big(\pi^Z_{\bigcup_{i \leq k} \mathrm{Inj}([i],e)}\Big)_\# \mu'(\mathrm{d}z) < \varepsilon$$

for any $e \in \binom{S}{k}$ (the choice being irrelevant by exchangeability).

The proof of the positive results of [10] will rely on a 'correspondence principle', which asserts that a hereditary property $\mathcal{P}$ is testable with one-sided error if it is infinitarily testable with one-sided error, and similarly for the stronger finitary and infinitary notions of 'strong local repairability'. Having established this principle, the bulk of the work of [10] will go into proving infinitary testability, using the formalism of exchangeable random colourings and the structure theorem 3.23. Interestingly, we find that both finitary and infinitary strong local repairability can fail for some hereditary properties, but they do always hold in the presence of various additional assumptions (for example, in rank 2 or if the property is actually monotone).

This approach to proving testability results has already appeared in the special case of hypergraph removal lemmas in Tao [62], making only partial use of the full structure theorem for exchangeable random hypergraphs.

## *4.6. Extremal problems*

Szemerédi's Theorem asserts that *any* positive-density subset of a sufficiently long discrete interval contains long arithmetic progressions. By contrast, the analogous search for copies of a given small hypergraph in a large dense hypergraph is typically not a foregone conclusion, and the nature of our question changes: rather than ask whether every positive-density hypergraph must contain embedded or induced copies of a fixed small hypergraph $F$, we instead turn to the classical Turán question of how dense our large hypergraph must be in order to guarantee the presence of such subhypergraphs.

However, structurally this adds a new kind of complexity; for in the typical case in which the Turán density lies in $(0, 1)$, identifying it precisely usually requires also a precise description of the density-extremizers among those hypergraphs that do not contain any embedded copies of the small hypergraph. Thus we must delve into the structure of individual hypergraphs in order to find these extremizers (at least up to their leading-order statistics), and so establish the value of the critical density. To date, this program has been successfully carried out for graphs, but only a handful of results are known for higher-rank hypergraphs; see, for example, the survey by Sidorenko [60] and the construction of conjectured edge-density extremizers among 3-uniform hypergraphs containing no tetrahedron in Kostochka [50] and (for a nice systematization of the Kostochka examples) Fon-Der-Flaass [26].

Such extremal questions also have natural formulations for exchangeable random hypergraphs. Let $H$ be a fixed finite $k$-uniform hypergraph on some finite vertex set $I \subset \mathbb{N}$; identify $H$ with some $\eta \in \{0,1\}^{\binom{I}{k}}$. We wish to know what is the supremal $\delta \geq 0$ (the **Turán density of** $H$) such that there exists an exchangeable random $k$-uniform hypergraph $\mu$ on $\{0,1\}^{\binom{\mathbb{N}}{k}}$ with $\mu(\pi_{\{1,2,\ldots,k\}}^{-1}\{1\}) = \delta$ and

$$\mu\{\omega \in X : \omega|_{\binom{I}{k}} \geq \eta\} = 0.$$

Note that it is clear by compactness and continuity that some $\mu$ attains this supremum. We refer to $\mu(\pi_{\{1,2,\ldots,k\}}^{-1}\{1\})$ as the **edge density of** $\mu$; and will say that $\mu$ is $H$-**free** in case the second $\mu$-probability above is zero. Simply by approximating such random hypergraphs by finite hypergraphs one sees that the resulting Turán density of $H$ is the same in either setting.

One can also ask the induced variant of this question, for which the condition $\omega|_{\binom{I}{2}} \geq \eta$ is replaced by $\omega|_{\binom{I}{2}} = \eta$.

Furstenberg's ergodic-theoretic approach to Szemerédi's Theorem can help to find a proof because we can reduce that theorem to a result about all non-negligible subsets of a probability-preserving system. On the other hand, the problem of understanding the structure of possible extremizers in the hypergraph setting, although easily formulable in both finitary and infinitary settings, is not clearly more vulnerable to attack in the latter than in the former. The structural ingredients of an exchangeable random hypergraph do offer a possibly richer collection of data over which to optimize when seeking to isolate the Turán extremizers, and, moreover, some of which can be varied continuously. However, it is not clear when this additional manoeuvrability might actually make the search any easier. One important instance of a related question in which such an infinitary formalism does seem to be more manageable is that studied by Razborov in [54, 55], where the interest is in what possible pairs of values $(\delta_1, \delta_2) \in [0,1]^2$ can appear (up to $o_N(1)$ corrections, for graphs with very large numbers of vertices $N$) simultaneously as the edge- and triangle-densities of a very large finite graph; Razborov uses a variational analysis on the objects of a different itary formalism (introduced in Subsection 4.3 above) to answer this question completely.

Finally, we note that such a variational approach to the structure of extremizing measures appears already in Section 4 of Fremlin and Talagrand [29]. Indeed, the motivation for their work is the analogous extremal questions for spreadable random graphs (see Subsection 3.9 above) for which the slightly weaker symmetry results in slightly different critical densities. The hypothesis of spreadability does not seem to relate to any simple finitary situation through a correspondence principle as in Subsection 2.3, and this may be partly why the Fremlin-Talagrand analysis seems to be largely unknown in the combinatorial community.

### *4.7. Broader context in ergodic theory*

An interesting viewpoint of the correspondence between statistics of large dense hypergraphs and exchangeable random hypergraph colourings, which we have not explained previously, is that it fills a suggestive lacuna among some of the known approaches to Szemerédi's Theorem.

On the one hand, two separate 'purely combinatorial' proofs of Szemerédi's Theorem are now known, both relying on strong regularity lemmas: Szemerédi's original proof ([61]) led to his introduction of the regularity lemma for graphs, and more recently hypergraph extensions of the regularity lemma have been introduced by Nagle, Rödl and Schacht [53] and by Gowers [36] and used to give a different combinatorial argument.

On the other hand, in 1977 Furstenberg [30] gave another proof, superficially very different, relying on a translation of the problem into the highly infinitary language of ergodic theory via a another correspondence principle, followed by an analysis of the relevant types of behaviour that an arbitrary measure-preserving system can display. This latter analysis relies on a 'structure theory' for probability preserving $\mathbb{Z}$-systems, developed by Furstenberg in [30] and by Zimmer in [69, 68]. An accessible introduction to the extent of this theory needed for a proof of Szemerédi's Theorem can be found in Furstenberg's book [31].

In essence, Furstenberg and Zimmer show how to extract from a probability-preserving $\mathbb{Z}$-system $(X, \Sigma, \mu, \tau)$ a (possibly transfinite) tower of factors

$$X \to \cdots \to Y_{\eta+1} \to Y_\eta \to \cdots \to Y_2 \to Y_1 \to Y_0$$

such that $Y_0$ is the invariant factor, for any limit ordinal $\eta$ the factor $Y_\eta$ is the inverse limit of its predecessors, and for each $\eta$ the extension $Y_{\eta+1} \to Y_\eta$ takes a certain 'primitive' form: it is either 'relatively weakly mixing' or 'relatively compact'. Relative weak mixing tells us that given a pair of functions $f, g \in L^\infty(Y_{\eta+1})$, $f$ and $g \circ \tau^n$ are asymptotically relatively independent over $Y_\eta$ (in a certain weak sense) as $n \to \infty$. Relative compactness is a counterpoint to this, according to which the shifts $g \circ \tau^n$ of the function $g \in L^\infty(Y_{\eta+1})$ are constrained in the forms they can take relative to the subfactor $Y_\eta$ as $n \to \infty$.

During the long period of research since these early papers it has become clear that many of the ideas underlying these two approaches to Szemerédi's Theorem are analogous. It turns out that the exchangeability theory can be seen as making this analogy a little more concrete, by offering an infinitary picture that relates to hypergraphs (via the correspondence principle of Subsection 2.3) as do measure-preserving $\mathbb{Z}$-systems to the additive combinatorics of Szemerédi's Theorem, and in which the structure theorems for exchangeable random hypergraph colourings then take the place of the Furstenberg-Zimmer tower.

However, notwithstanding the analogy outlined above, our structural analysis differs from that for $\mathbb{Z}$-actions in some instructive ways.

Firstly, while the work of Furstenberg and Zimmer applies to an arbitrary probability-preserving group action (although some further refinements are needed to extract the multidimensional Szemerédi Theorem; see [32]), in the setting of exchangeable random hypergraph colourings or set-indexed systems we are considering only a very special restricted class of probability-preserving $\mathrm{Sym}_0(S)$-actions, precisely because we insist that our measures live on the product space $K^T$ and be invariant under coordinate permutations. Such special cases are very far from representing a generic $\mathrm{Sym}_0(S)$-system.

Secondly, and perhaps more curiously, our structural description necessitates a slight extension of the notion of factor that appears in the Furstenberg-Zimmer theory: we work instead with a tower of 'quasifactors' (for which the relevant probability kernel is not directly recoverable from the original system, but only after our trick of introducing additional 'auxiliary vertices'), and the second example of Subsection 3.6 shows that this is really necessary. Each quasifactor in our tower has a simple description in terms of its predecessors, but in general those predecessors may not be recoverable as factors of the original system. While quasifactors certainly do play a rôle in general ergodic theory also (see, for example, Chapter 8 of Glasner [35]), it seems that to date the basic ergodic theoretic analyses of questions of multiple recurrence (or, relatedly, convergence of nonconventional ergodic averages, as in [43, 67]) has always proved reducible to the study of true factors of the original system.

However, perhaps the most instructive difference is in the basic constructions that underly the proofs. Both the exchangeable law structural results and the Furstenberg-Zimmer theory rest on our ability to extract, given a suitable exchangeable or stationary law $\mu$ on the product space $K^T$ (with $T = \mathbb{Z}$ in classical ergodic theory), a $\Gamma$-invariant $\sigma$-subalgebra T of $\Sigma_K^{\otimes T}$ such that, firstly, T enjoys some strengthened symmetry or other additional structure under the law $\mu$ (up to negligible sets), and secondly the coordinate projections $\pi_t$ enjoy some additional regularity (such as relative independence or some weaker notion related to it, such as relative weak mixing) when conditioned on T under $\mu$. Typically a 'composition series' comprising several nested such $\sigma$-subalgebras is called for.

Of course, the details of these additional demands made on the $\sigma$-subalgebra vary among different notions of exchangeability, and are different again in the ergodic-theoretic setting — indeed, much of the ingenuity behind each theory lies in the selection of the right demands to impose. However, in some sense the most basic difference between the methods of exchangeability theory and of Furstenberg-Zimmer theory lies in *how these $\sigma$-subalgebras are specified*.

In the case of exchangeable processes, the decisive observation for solving our problem is that we can embed a copy of the index set $T$ as an infinite-coinfinite subset $T_1$ of itself, so that the subprocess $(\pi_t)_{t \in T_1}$ is effectively just a copy of the original process and so that the remainder process $(\pi_t)_{t \in T \setminus T_1}$ also still enjoys enough symmetry that we can construct a suitable auxiliary process $(\pi_t^Z)_{t \in T_1}$ simply by grouping together suitable further subprocesses of this remainder process. The clusters of random variables obtained from this grouping then specify directly the factor needed to describe the subprocess $(\pi_t)_{t \in T_1}$, and so by our earlier identification they also tell us the law of the original process $(\pi_t)_{t \in T}$. Thus, with the right sleight of hand, we can in this setting specify our auxiliary process directly by grouping into subprocesses random variables from the original process.

In practice, this approach clearly relies on the ability to move the coordinate projections $\pi_t$ about by the action $\Gamma \curvearrowright T$ very flexibly; in particular, in the proofs of the relative-independence assertions of Propositions 3.4 and 3.12 we made crucial use of our ability to move some vertices around almost arbitrarily while leaving certain other vertices fixed.

In ergodic-theoretic problems, however, we must usually do without this 'overwhelming strength' of the action of $\Gamma$ on $T$. In this setting, is is therefore more common to obtain the $\sigma$-subalgebra T (or, equivalently, the random variables $\pi_t^Z$ or the factor they specify) as suitable limits of functions of several of the original projections $\pi_t$, either by averaging over the action of $\Gamma$ (which is usually amenable) or sometimes by taking some more exotic kind of limit, such as the 'IP-limits' that under Furstenberg and Katznelson's analysis of general 'IP systems' in [33]. This approach generally requires more hard analysis, firstly to show that our averages or limits defining the $\pi_t^Z$ exist at all, and then that they specify an auxiliary process that still enjoys the desired properties (since this may not now be so easy to read off from the corresponding properties of $(\pi_t)_{t \in T}$).

Often — particularly in ergodic theory — there is more than one way to extract these auxiliary processes, and our knowledge of structural results has improved over time owing to refinements to the techniques that are available. In particular, the original Furstenberg-Zimmer structural approach — while still the basic foundation of the analysis of multiple recurrence phenomena for $\mathbb{Z}$-actions — has gradually been complemented by a much more precise description of those 'factors' of a probability-preserving $\mathbb{Z}$-system that govern the convergence behaviour of various 'nonconventional ergodic averages' (which are somewhat analogous to individual products of observables in our setting). We direct the reader, in particular, to the two slightly different treatments of these 'characteristic factors' by Host and Kra [43] and by Ziegler [67].

As a result of this difference in extraction-method, the Furstenberg-Zimmer structure theory (and its later relatives) uses factors built from the 'bottom up':

starting with the raw $\mathbb{Z}$-action, a tower of factors is built recursively upwards, each obtained by an extension of its predecessor (and with inverse limits taken up to limit ordinals, if desired). The resulting description of the overall system is much less precise than in the exchangeability theory: the Furstenberg-Zimmer structure theorem yields a large tower of 'relatively compact' extensions (a very special kind), followed by one last 'relatively weakly mixing' extension (about which we know very little, except that it is irrelevant for the estimation of various nonconventional averages).

In contrast, our structure theorem for an exchangeable random hypergraph colouring is proved from the top-down: starting from an invariant random hypergraph colouring $\mu$ of rank $k$, we pass in one step down to a quasifactor that is itself a rank-$(k-1)$ random hypergraph colouring, and show that the original $\mu$ enjoys the desired relative independence of the random variables $\pi_e$ over this quasifactor (a much more precise condition than 'relative weak mixing'). After iterating this argument $k$ times, we have obtained the whole structure of the original $\mu$, leaving only very completely-described degrees of freedom at each step.

Naturally the above distinction is not precise, nor does it completely account for all manners of specifying $\sigma$-subalgebras that appear in this field. In particular, many of the constructions on the ergodic theory side can be recovered by considering instead a suitable invariant $\sigma$-algebra for some measure-preserving action, either within the original system or within some extension of it (the work [43] of Host and Kra providing a particularly striking example of such an approach). However, such a $\sigma$-subalgebra is still precisely that which is generated, up to $\mu$-negligible sets, by the abovementioned ergodic averages; and such averages inevitably then surface elsewhere in the proof that these $\sigma$-algebras have the desired properties, so this difference seems to be rather cosmetic.

## Acknowledgements

My thanks go to Terence Tao for frequent advice and guidance, to David Fremlin for bringing [29] to my attention and for helpful suggestions and to Alexander Razborov for further helpful comments.

## References

[1] D. J. Aldous. Representations for partially exchangeable arrays of random variables. *J. Multivariate Anal.*, 11(4):581–598, 1981. MR0637937

[2] D. J. Aldous. On exchangeability and conditional independence. In *Exchangeability in probability and statistics (Rome, 1981)*, pages 165–170. North-Holland, Amsterdam, 1982. MR0675972

[3] D. J. Aldous. Exchangeability and related topics. In *École d'été de probabilités de Saint-Flour, XIII—1983*, volume 1117 of *Lecture Notes in Math.*, pages 1–198. Springer, Berlin, 1985. MR0883646

[4] D. J. Aldous and R. Lyons. Processes on Unimodular Random Networks. *Electronic J. Probab.*, 12:1454–1508, 2007. MR2354165

[5] D. J. Aldous and J. M. Steele. The objective method: probabilistic combinatorial optimization and local weak convergence. In H. Kesten, editor, *Probability on Discrete Structures*, volume 110 of *Enyclopaedia Math. Sci.*, pages 1–72. Springer, Berlin, 2004. MR2023650

[6] N. Alon, E. Fischer, M. Krivelevich, and M. Szegedy. Efficient testing of large graphs. *Combinatorica*, 20:451–476, 2000. MR1804820

[7] N. Alon and A. Shapira. A Characterization of the (natural) Graph Properties Testable with One-Sided Error. preprint, available online at http://www.math.tau.ac.il/~nogaa/PDFS/heredit2.pdf.

[8] N. Alon and A. Shapira. Every monotone graph property is testable. In *Proc. of the $37^{\text{th}}$ ACM STOC, Baltimore*. ACM Press, 2005. available online at http://www.math.tau.ac.il/~nogaa/PDFS/MonotoneSTOC.pdf. MR2181610

[9] T. Austin. Razborov flag algebras as algebras of measurable functions. manuscript, available online at arXiv.org: 0801.1538, 2007. MR2371204

[10] T. Austin and T. Tao. On the testability and repair of hereditary hypergraph properties. preprint, available online at arXiv.org: 0801.2179, 2008.

[11] H. Becker and A. S. Kechris. *The Descriptive Set Theory of Polish Group Actions*. Cambridge University Press, Cambridge, 1996. MR1425877

[12] I. Benjamini and O. Schramm. Recurrence of distributional limits of finite planar graphs. *Electron. J. Probab.*, 6:no. 23, 13 pp. (electronic), 2001. MR1873300

[13] I. Benjamini, O. Schramm, and A. Shapira. Every Minor-Closed Property of Sparse Graphs is Testable.

[14] Y. Benyamini and J. Lindenstrauss. *Geometric Nonlinear Functional Analysis*. American Mathematical Society, Providence, 2000. MR1727673

[15] V. Bergelson. Ergodic Ramsey Theory – an Update. In M. Pollicott and K. Schmidt, editors, *Ergodic Theory of $\mathbb{Z}^d$-actions: Proceedings of the Warwick Symposium 1993-4*, pages 1–61. Cambridge University Press, Cambridge, 1996. MR1411215

[16] B. Bollobás. *Modern Graph Theory*. Springer, Berlin, 1998. MR1633290

[17] C. Borgs, J. Chayes, L. Lovász, V. T. Sós, B. Szegedy, and K. Vesztergombi. Graph limits and parameter testing. In *STOC'06: Proceedings of the 38th Annual ACM Symposium on Theory of Computing*, pages 261–270, New York, 2006. ACM. MR2277152

[18] B. de Finetti. Fuzione caratteristica di un fenomeno aleatorio. *Mem. R. Acc. Lincei*, 4(6):86–133, 1930.

[19] B. de Finetti. La prévision: ses lois logiques, ses sources subjectives. *Ann. Inst. H. Poincaré*, 7:1–68, 1937. MR1508036

[20] P. Diaconis and S. Janson. Graph limits and exchangeable random graphs. Preprint; available online at arXiv.org: math.PR math.CO/0712.2749, 2007. MR2346812

[21] E. B. Dynkin. Classes of equivalent random quantities. *Uspehi Matem. Nauk (N.S.)*, 8(2(54)):125–130, 1953. MR0055601

[22] G. Elek. On limits of finite graphs. preprint; available online at arXiv.org: math.CO/0505335, 2005. MR2359831
[23] G. Elek. A Regularity Lemma for Bounded Degree Graphs and Its Applications: Parameter Testing and Infinite Volume Limits. preprint; available online at arXiv.org: math.CO/0711.2800, 2007.
[24] G. Elek and B. Szegedy. Limits of Hypergraphs, Removal and Regularity Lemmas. A Nonstandard Approach. preprint; available online at arXiv.org: math.CO/0705.2179, 2007.
[25] P. Erdős and A. Hajnal. Some remarks on set theory, IX. Combinatorial problems in measure theory and set theory. *Michigan Math. J.*, 11:107–127, 1964. MR0171713
[26] D. G. Fon-Der-Flaass. A method for constructing $(3, 4)$-graphs. *Mat. Zametki*, 44(4):546–550, 559, 1988. MR0975195
[27] D. H. Fremlin. list of problems. available online at `http://www.essex.ac.uk/maths/staff/fremlin/problems.htm`.
[28] D. H. Fremlin. Random equivalence relations. preprint, available online at `http://www.essex.ac.uk/maths/staff/fremlin/preprints.htm`, 2006.
[29] D. H. Fremlin and M. Talagrand. Subgraphs of random graphs. *Trans. Amer. Math. Soc.*, 291(2):551–582, 1985. MR0800252
[30] H. Furstenberg. Ergodic behaviour of diagonal measures and a theorem of Szemerédi on arithmetic progressions. *J. d'Analyse Math.*, 31:204–256, 1977. MR0498471
[31] H. Furstenberg. *Recurrence in Ergodic Theory and Combinatorial Number Theory*. Princeton University Press, Princeton, 1981. MR0603625
[32] H. Furstenberg and Y. Katznelson. An ergodic Szemerédi Theorem for commuting transformations. *J. d'Analyse Math.*, 34:275–291, 1978. MR0531279
[33] H. Furstenberg and Y. Katznelson. An ergodic Szemerédi theorem for IP-systems and combinatorial theory. *J. d'Analyse Math.*, 45:117–168, 1985. MR0833409
[34] H. Gaifman. Concerning measures in first-order calculi. *Israel J. Math.*, 2:1–18, 1964. MR0175755
[35] E. Glasner. *Ergodic Theory via Joinings*. American Mathematical Society, Providence, 2003. MR1958753
[36] W. T. Gowers. Hypergraph regularity and the multidimensional Szemerédi Theorem. preprint.
[37] W. T. Gowers. A new proof of Szemerédi's theorem. *Geom. Funct. Anal.*, 11(3):465–588, 2001. MR1844079
[38] W. T. Gowers. Quasirandomness, counting and regularity for 3-uniform hypergraphs. *Combin. Probab. Comput.*, 15(1-2):143–184, 2006. MR2195580
[39] M. Gromov. *Metric Structures for Riemannian and Non-Riemannian Spaces*. Birkhäuser, Boston, 1999. MR1699320
[40] E. Hewitt and L. J. Savage. Symmetric measures on Cartesian products. *Trans. Amer. Math. Soc.*, 80:470–501, 1955. MR0076206
[41] D. N. Hoover. Relations on probability spaces and arrays of random variables. 1979.
[42] D. N. Hoover. Row-columns exchangeability and a generalized model for

exchangeability. In *Exchangeability in probability and statistics (Rome, 1981)*, pages 281–291, Amsterdam, 1982. North-Holland. MR0675982
- [43] B. Host and B. Kra. Nonconventional ergodic averages and nilmanifolds. *Ann. Math.*, 161(1):397–488, 2005. MR2150389
- [44] O. Kallenberg. Symmetries on random arrays and set-indexed processes. *J. Theoret. Probab.*, 5(4):727–765, 1992. MR1182678
- [45] O. Kallenberg. *Foundations of modern probability.* Probability and its Applications (New York). Springer-Verlag, New York, second edition, 2002. MR1876169
- [46] O. Kallenberg. *Probabilistic symmetries and invariance principles.* Probability and its Applications (New York). Springer, New York, 2005. MR2161313
- [47] J. F. C. Kingman. The representation of partition structures. *J. London Math. Soc. (2)*, 18(2):374–380, 1978. MR0509954
- [48] J. F. C. Kingman. Uses of exchangeability. *Ann. Probability*, 6(2):183–197, 1978. MR0494344
- [49] R. Kopperman. *Model Theory and its Applications.* Allyn and Bacon, Boston, 1972. MR0363873
- [50] A. V. Kostochka. A class of constructions for Turán's $(3,4)$-problem. *Combinatorica*, 2:187–192, 1982. MR0685045
- [51] P. H. Krauss. Representation of symmetric probability models. *J. Symbolic Logic*, 34:183–193, 1969. MR0275482
- [52] L. Lovász and B. Szegedy. Limits of dense graph sequences. *J. Combin. Theory Ser. B*, 96(6):933–957, 2006. MR2274085
- [53] B. Nagle, V. Rödl, and M. Schacht. The counting lemma for regular $k$-uniform hypergraphs. *Random Structures and Algorithms*, to appear. MR2198495
- [54] A. Razborov. Flag Algebras. *J. Symbolic Logic*, 72(4):1239–1282, 2007. MR2371204
- [55] A. Razborov. On the minimal density of triangles in graphs. preprint; available online at http://www.mi.ras.ru/~razborov/triangles.pdf, 2007.
- [56] V. Rödl and M. Schacht. Generalizations of the removal lemma. preprint, available online at http://www.informatik.hu-berlin.de/~schacht/pub/preprints/gen_removal.pdf.
- [57] R. Rubinfeld and M. Sudan. Robust characterization of polynomials with applications to program testing. *SIAM Journal on Computing*, 25(2):252–271, 1996. MR1379300
- [58] C. Ryll-Nardzewski. On stationary sequences of random variables and the de finetti's equivalence. *Colloq. Math.*, 4:149–156, 1957. MR0088823
- [59] O. Schramm. Hyperfinite graph limits. preprint, available online at arXiv.org: math.CO/0711.3808, 2007. MR2334202
- [60] A. Sidorenko. What We Know and What We Do not Know about Turán Numbers. *Graphs and Combinatorics*, 11:179–199, 1995. MR1341481
- [61] E. Szemerédi. On sets of integers containing no $k$ elements in arithmetic progression. *Acta Arith.*, 27:199–245, 1975. MR0369312
- [62] T. Tao. A correspondence principle between (hyper)graph theory and

probability theory, and the (hyper)graph removal lemma. *J. d'Analyse Mathematique.* to appear; available online at arXiv.org: math.CO/0602037. MR2373263
[63] T. Tao. A quantitative ergodic theory proof of Szemerédi's theorem. *Electron. J. Combin.*, 13(1):Research Paper 99, 49 pp. (electronic), 2006. MR2274314
[64] T. Tao and V. Vu. *Additive combinatorics.* Cambridge University Press, Cambridge, 2006. MR2289012
[65] P. Wojtaszczyk. *Banach spaces for analysts.* Cambridge University Press, Cambridge, 1991. MR1144277
[66] K. Yosida. *Functional analysis.* Classics in Mathematics. Springer-Verlag, Berlin, 1995. Reprint of the sixth (1980) edition. MR1336382
[67] T. Ziegler. Universal characteristic factors and Furstenberg averages. *J. Amer. Math. Soc.*, 20(1):53–97 (electronic), 2007. MR2257397
[68] R. J. Zimmer. Ergodic actions with generalized discrete spectrum. *Illinois J. Math.*, 20(4):555–588, 1976. MR0414832
[69] R. J. Zimmer. Extensions of ergodic group actions. *Illinois J. Math.*, 20(3):373–409, 1976. MR0409770